\title{Spectral Measures for $Sp(2)$}
\author{
{\sc David E.\ Evans and Mathew Pugh}\\
 {\footnotesize School of Mathematics, Cardiff University,}\\  {\footnotesize Senghennydd Road, Cardiff CF24 4AG, Wales, U.K.}
}
\date{\today}

\documentclass[12pt]{article}
\usepackage{amssymb,amsmath}
\usepackage{graphicx}
\usepackage{color}
\usepackage[usenames,dvipsnames]{xcolor}

\newtheorem{Def}{Definition}[section]

\newtheorem{Cor}[Def]{Corollary}
\newtheorem{Thm}[Def]{Theorem}

\begin{document}
\maketitle

\begin{abstract}
Spectral measures provide invariants for braided subfactors via fusion modules. In this paper we study joint spectral measures associated to the compact connected rank two Lie group $SO(5)$ and its double cover the compact connected, simply-connected rank two Lie group $Sp(2)$, including the McKay graphs for the irreducible representations of $Sp(2)$ and $SO(5)$ and their maximal tori, and fusion modules associated to the $Sp(2)$ modular invariants.
\end{abstract}

{\footnotesize
\tableofcontents
}

\section{Introduction} \label{sect:intro}

Spectral measures associated to the compact Lie groups $SU(2)$, $SU(3)$ and $G_2$, their maximal tori, nimrep graphs associated to the $SU(2)$, $SU(3)$ and $G_2$ modular invariants, and the McKay graphs for finite subgroups of $SU(2)$, $SU(3)$ and $G_2$ were studied in \cite{banica/bisch:2007, evans/pugh:2009v, evans/pugh:2010i, evans/pugh:2012i, evans/pugh:2012ii}.
Spectral measures associated to other compact rank two Lie groups and their maximal tori are studied in \cite{evans/pugh:2012iv}.

For the $SU(2)$ and $SU(3)$ graphs, the spectral measures distill onto very special subsets of the semicircle/circle for $SU(2)$ (which are both one-dimensional spaces) and discoid/torus for $SU(3)$ (which are both two-dimensional spaces), and the theory of nimreps allowed us to compute these measures precisely. Our methods gave an alternative approach to deriving the results of Banica and Bisch \cite{banica/bisch:2007} for $ADE$ graphs and subgroups of $SU(2)$, and explained the connection between their results for affine $ADE$ graphs and the Kostant polynomials.
In the case of $G_2$, the spectral measures distill onto subsets of $\mathbb{R}$ and the torus $\mathbb{T}^2$, which are one-dimensional and two-dimensional respectively, resulting in an infinite family of pullback measures over $\mathbb{T}^2$ for any spectral measure on $\mathbb{R}$. This ambiguity was removed by considering instead joint spectral measures for pairs of graphs corresponding to the two fundamental representations of $G_2$. Such joint spectral measures, which have support in $\mathbb{R}^2$, yield a unique pullback measure over $\mathbb{T}^2$, and the spectral measures are obtained as pushforward measures.

In this paper we study spectral measures for the compact, connected, simply-connected rank two Lie group $Sp(2)$, the group of $4\times4$ unitary symplectic matrices with entries in $\mathbb{C}$. We also study spectral measures for the (non-simply-connected) compact rank two Lie group $SO(5)$, the group of $5\times5$ real orthogonal matrices, whose double cover is $Sp(2)$. In particular we determine the joint spectral measures associated to the Lie groups themselves and their maximal tori, and joint spectral measures for nimrep graphs associated to the $Sp(2)$ modular invariants.

In $C_2 = sp(2)$ (the Lie algebra of $Sp(2)$) conformal field theories, one considers the Verlinde algebra at a finite level $k$, which is represented by a non-degenerately braided system ${}_N \mathcal{X}_N$ of irreducible endomorphisms on a type $\mathrm{III}_1$ factor $N$, whose fusion rules $\{ N_{\lambda \nu}^{\mu} \}$ reproduce exactly those of the positive energy representations of the loop group of $Sp(2)$ at level $k$,
$N_{\lambda} N_{\mu} = \sum_{\nu} N_{\lambda \nu}^{\mu} N_{\nu}$.
The statistics generators $S$, $T$ for the braided tensor category ${}_N \mathcal{X}_N$ match exactly those of the Ka\u{c}-Peterson modular $S$, $T$ matrices which perform the conformal character transformations (see footnote 2 in \cite{bockenhauer/evans:2001}).
The fusion graph for these irreducible endomorphisms are truncated versions of the representation graphs of $Sp(2)$ itself (see Section \ref{sect:measures_AkB2}).
From the Verlinde formula (\ref{eqn:verlinde_formula}) we see that this family $\{ N_{\lambda} \}$ of commuting normal matrices can be simultaneously diagonalised:
\begin{equation} \label{eqn:verlinde_formula}
N_{\lambda} = \sum_{\sigma} \frac{S_{\sigma, \lambda}}{S_{\sigma,0}} S_{\sigma} S_{\sigma}^{\ast},
\end{equation}
where the summation is over each $\sigma \in {}_N \mathcal{X}_N$ and $0$ is the trivial representation. It is intriguing that the eigenvalues $S_{\sigma, \lambda}/S_{\sigma,0}$ and eigenvectors $S_{\sigma} = \{ S_{\sigma, \mu} \}_{\mu}$ are described by the modular $S$ matrix.

A braided subfactor is an inclusion $N \subset M$ where the dual canonical endomorphism decomposes as a finite combination of endomorphisms in ${}_N \mathcal{X}_N$, and yields a modular invariant partition function through the procedure of $\alpha$-induction which allows two extensions of $\lambda$ on $N$, depending on the use of the braiding or its opposite, to endomorphisms $\alpha^{\pm}_{\lambda} \in {}_M \mathcal{X}_M^{\pm}$ of $M$, so that the matrix $Z_{\lambda,\mu} = \langle \alpha_{\lambda}^+, \alpha_{\mu}^- \rangle$ is a modular invariant \cite{bockenhauer/evans/kawahigashi:1999, bockenhauer/evans:2000, evans:2003}.
The systems ${}_M \mathcal{X}_M^{\pm}$ are called the chiral systems, whilst the intersection ${}_M \mathcal{X}_M^0 = {}_M \mathcal{X}_M^+ \cap {}_M \mathcal{X}_M^-$ is the neutral system. Then ${}_M \mathcal{X}_M^0 \subset {}_M \mathcal{X}_M^{\pm} \subset {}_M \mathcal{X}_M$, where ${}_M \mathcal{X}_M \subset \mathrm{End}(M)$ denotes a system of endomorphisms consisting of a choice of representative endomorphisms of each irreducible subsector of sectors of the from $[\iota \lambda \overline{\iota}]$, $\lambda \in {}_N \mathcal{X}_N$, where $\iota: N \hookrightarrow M$ is the inclusion map. Although ${}_N \mathcal{X}_N$ is assumed to be braided, the systems ${}_M \mathcal{X}_M^{\pm}$ or ${}_M \mathcal{X}_M$ are not braided in general.
The action of the $N$-$N$ sectors ${}_N \mathcal{X}_N$ on the $M$-$N$ sectors ${}_M \mathcal{X}_N$ and produces a nimrep (non-negative integer matrix representation of the original Verlinde algebra) $G_{\lambda} = (\langle \xi \lambda, \xi' \rangle)_{\xi,\xi' \in {}_M\mathcal{X}_N}$, i.e.
$G_{\lambda} G_{\mu} = \sum_{\nu} N_{\lambda \nu}^{\mu} G_{\nu}$
whose spectrum reproduces exactly the diagonal part of the modular invariant \cite{bockenhauer/evans/kawahigashi:2000}.
In the case of the trivial embedding of $N$ in itself, the nimrep $G$ is simply the trivial representation $N$.
Since the nimreps are a family of commuting matrices, they can be simultaneously diagonalised and thus the eigenvectors $\psi_{\sigma}$ of $G_{\lambda}$ are the same for each $\lambda \in {}_N \mathcal{X}_N$.
We have
\begin{equation} \label{eqn:verlinde_formulaG}
G_{\lambda} = \sum_{\sigma} \frac{S_{\sigma,\lambda}}{S_{\sigma,0}} \psi_{\sigma} \psi_{\sigma}^{\ast},
\end{equation}
where the summation is over each $\sigma \in {}_N \mathcal{X}_N$ with multiplicity given by the modular invariant, i.e. the spectrum of $G_{\lambda}$ is given by $\{ S_{\sigma, \lambda}/S_{\sigma,0}$ with multiplicity $Z_{\sigma,\sigma} \}$. We call the set $ \{ \mu$ with multiplicity $Z_{\mu,\mu} \}$ the set of exponents of $G$.

Along with the identity invariants for $Sp(2)$, there are orbifold invariants for all levels $k$ \cite{bernard:1987}. There are three exceptional invariants due to conformal embeddings at levels 3, 7, 12 \cite{christe/ravanani:1989}. There is also an exceptional invariant at level 8 \cite{verstegen:1990} which is a twist of the orbifold invariant at level 8, and is analogous to the $E_7$ modular invariant for $SU(2)$ \cite[$\S5.3$]{bockenhauer/evans/kawahigashi:2000} and the Moore-Seiberg $\mathcal{E}_{MS}^{(12)}$ invariant for $SU(3)$ \cite[$\S5.4$]{evans/pugh:2009ii}.
These are all the known $Sp(2)$ modular invariants, and the list is complete for all $k\leq26$ \cite{gannon/ho-kim:1994}.

This paper is organised as follows.
In Section \ref{sect:measure:A8infty(B2)} we describe the representation theory of $Sp(2)$ and $SO(5)$, and their maximal torus $\mathbb{T}^2$, and in particular focus on their fundamental representations. In Sections \ref{sect:measure:A8inftyC2-T2}-\ref{sect:measureA8inftyB2-R} we determine the (joint) spectral measures associated to the (adjacency matrices of the) McKay graphs given by the action of the irreducible characters of $Sp(2)$ on its maximal torus $\mathbb{T}^2$, and the analogous results for $SO(5)$. In Section \ref{sect:measureAinftyB2} we determine the (joint) spectral measures associated to the (adjacency matrices of the) McKay graphs of $Sp(2)$ and $SO(5)$ themselves. In all these cases we focus on the fundamental representations of $Sp(2)$, $SO(5)$ respectively, and determine these (joint) spectral measures over both $\mathbb{T}^2$ and the (joint) spectrum of these adjacency matrices.
Finally in Section \ref{sect:measures_nimrepB2} we determine joint spectral measures over $\mathbb{T}^2$ for nimrep graphs arising from $Sp(2)$ braided subfactors.

\section{Spectral measures for ${}^W \hspace{-2mm} \mathcal{A}_{\infty}(Sp(2))$, ${}^W \hspace{-2mm} \mathcal{A}_{\infty}(SO(5))$} \label{sect:measure:A8infty(B2)}

The irreducible representations $\lambda_{(\mu_1,\mu_2)}$ of $Sp(2)$ are indexed by pairs $(\mu_1,\mu_2) \in \mathbb{N}^2$ such that $\mu_1 \geq \mu_2$.
Let the fundamental representation $\rho_x = \lambda_{(1,0)}$ be the standard representation of $Sp(2)$, $\rho_x(Sp(2)) = Sp(2)$, the group of $4\times4$ unitary symplectic matrices with entries in $\mathbb{C}$.
The maximal torus of $Sp(2)$ is $T = \textrm{diag}(t_1,t_2,t_1^{-1},t_2^{-1})$, for $t_i \in \mathbb{T}$, which is isomorphic to $\mathbb{T}^2$, so that the restriction of $\rho_x$ to $\mathbb{T}^2$ is given by
the $4 \times 4$ diagonal matrix
\begin{equation} \label{eqn:restrict_rho1B2_to_T2}
(\rho_x|_{\mathbb{T}^2})(\omega_1,\omega_2) = \textrm{diag}(\omega_1,\omega_2,\omega_1^{-1},\omega_2^{-1}),
\end{equation}
for $(\omega_1,\omega_2) \in \mathbb{T}^2$.

Let the fundamental representation $\rho_y = \lambda_{(1,1)}$ be the standard representation of $SO(5)$, $\rho_y(Sp(2)) = SO(5)$, the group of $5\times5$ real orthogonal matrices. The restriction of $\rho_y$ to $\mathbb{T}^2$ is given by the $5 \times 5$ diagonal matrix
\begin{equation} \label{eqn:restrict_rho2B2_to_T2}
(\rho_y|_{\mathbb{T}^2})(\omega_1,\omega_2) = \textrm{diag}(\omega_1\omega_2, \omega_1^{-1}\omega_2^{-1}, \omega_1\omega_2^{-1}, \omega_1^{-1}\omega_2,1),
\end{equation}
for $(\omega_1,\omega_2) \in \mathbb{T}^2$.

The irreducible representations of $SO(5)$ are given by the representations $\lambda_{(\mu_1,\mu_2)}$ of $Sp(2)$ for which $\mu_1+\mu_2$ is even. In order to study spectral measures associated to $SO(5)$, we take $\rho_y$ and a second fundamental representation $\rho_z = \lambda_{(2,0)}$ of $SO(5)$, which is the adjoint representation of $Sp(2)$ of dimension 10. The restriction of $\rho_z$ to $\mathbb{T}^2$ is given by the $10 \times 10$ diagonal matrix
\begin{equation} \label{eqn:restrict_rho3B2_to_T2}
(\rho_z|_{\mathbb{T}^2})(\omega_1,\omega_2) = \textrm{diag}(\omega_1^2, \omega_2^2, \omega_1^{-2}, \omega_2^{-2}, \omega_1\omega_2, \omega_1\omega_2^{-1}, \omega_1^{-1}\omega_2, \omega_1^{-1}\omega_2^{-1}, 1, 1),
\end{equation}
for $(\omega_1,\omega_2) \in \mathbb{T}^2$.

Let $\{ \chi_{(\mu_1,\mu_2)} \}_{\mu_1,\mu_2 \in \mathbb{N}:\mu_1 \geq \mu_2}$, $\{ \sigma_{(\mu_1,\mu_2)} \}_{\mu_1,\mu_2 \in \mathbb{Z}}$ be the irreducible characters of $Sp(2)$, $\mathbb{T}^2$ respectively, where $\chi_{(\mu_1,\mu_2)} := \chi_{\lambda_{(\mu_1,\mu_2)}}$. The characters $\chi_{(\mu_1,\mu_2)}$ of $Sp(2)$ are self-conjugate and thus are maps from the torus $\mathbb{T}^2$ to an interval $I_{\mu} := \chi_{\mu}(\mathbb{T}^2) \subset \mathbb{R}$.
For $\omega_i \in \mathbb{T}$, $\mu_i \in \mathbb{Z}$, the characters of $\mathbb{T}^2$ are given by $\sigma_{(\mu_1,\mu_2)}(\omega_1,\omega_2) = \omega_1^{\mu_1}\omega_2^{\mu_2}$, and satisfy $\overline{\sigma_{(\mu_1,\mu_2)}} = \sigma_{(-\mu_1,-\mu_2)}$.

\begin{figure}[tb]
\begin{center}
  \includegraphics[width=140mm]{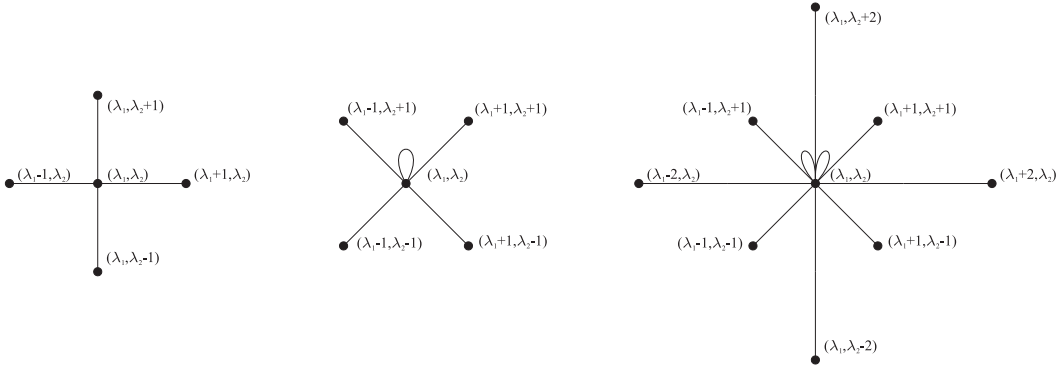}
 \caption{Multiplication by $\chi_{\rho_x}|_{\mathbb{T}^2}$, $\chi_{\rho_y}|_{\mathbb{T}^2}$ and $\chi_{\rho_z}|_{\mathbb{T}^2}$} \label{fig-mult_rhoB2}
\end{center}
\end{figure}

\begin{figure}[tb]
\begin{center}
  \includegraphics[width=155mm]{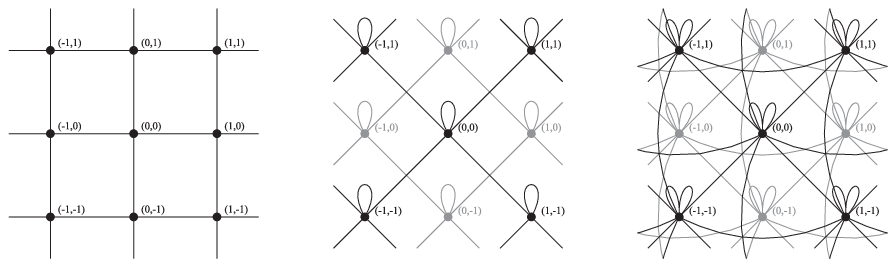}
 \caption{Infinite graphs ${}^W \hspace{-2mm} \mathcal{A}^{\rho_x}_{\infty}(Sp(2))$, ${}^W \hspace{-2mm} \mathcal{A}^{\rho_y}_{\infty}(Sp(2))$ and ${}^W \hspace{-2mm} \mathcal{A}^{\rho_z}_{\infty}(Sp(2))$} \label{fig-A_8infty(B2)all}
\end{center}
\end{figure}

If $\sigma_u$ is the restriction of $\chi_{\rho_u}$ to $\mathbb{T}^2$, $u=x,y,z$, then from (\ref{eqn:restrict_rho1B2_to_T2})-(\ref{eqn:restrict_rho3B2_to_T2})
\begin{align}
\sigma_x & = \chi_{(1,0)}|_{\mathbb{T}^2} = \sigma_{(1,0)} + \sigma_{(-1,0)} + \sigma_{(0,1)} + \sigma_{(0,-1)}, \label{eqn:restriction-chi1} \\
\sigma_y & = \chi_{(1,1)}|_{\mathbb{T}^2} = \sigma_{(0,0)} + \sigma_{(1,1)} + \sigma_{(-1,-1)} + \sigma_{(1,-1)} + \sigma_{(-1,1)}, \label{eqn:restriction-chi2} \\
\sigma_z & = \chi_{(2,0)}|_{\mathbb{T}^2} = 2\sigma_{(0,0)} + \sigma_{(2,0)} + \sigma_{(-2,0)} + \sigma_{(0,2)} + \sigma_{(0,-2)} + \sigma_{(1,1)} + \sigma_{(-1,-1)} + \sigma_{(1,-1)} + \sigma_{(-1,1)}. \label{eqn:restriction-chi3}
\end{align}
Then
\begin{equation} \label{eqn:fusion_rules-B2-1}
\sigma_x \sigma_{(\mu_1,\mu_2)} = \sigma_{(\mu_1+1,\mu_2)} + \sigma_{(\mu_1-1,\mu_2)} + \sigma_{(\mu_1,\mu_2+1)} + \sigma_{(\mu_1,\mu_2-1)},
\end{equation}
for any $\mu_1,\mu_2 \in \mathbb{Z}$, where multiplication by $\sigma_x = \chi_{\rho_x}|_{\mathbb{T}^2}$ corresponds to the edges illustrated in the first diagram in Figure \ref{fig-mult_rhoB2}. The representation graph of $\mathbb{T}^2$ for the first fundamental representation $\rho_x$ is identified with the infinite graph ${}^W \hspace{-2mm} \mathcal{A}^{\rho_x}_{\infty}(Sp(2))$, which is the first figure illustrated in Figure \ref{fig-A_8infty(B2)all}, whose vertices may be labeled by pairs $(\mu_1,\mu_2) \in \mathbb{Z}^2$ such that there is an edge from $(\mu_1,\mu_2)$ to $(\mu_1+1,\mu_2)$, $(\mu_1-1,\mu_2)$, $(\mu_1,\mu_2+1)$ and $(\mu_1,\mu_2-1)$.


Similarly, the representation graph of $\mathbb{T}^2$ for the irreducible representations $\rho_y$, $\rho_z$ are identified with the infinite graphs ${}^W \hspace{-2mm} \mathcal{A}^{\rho_y}_{\infty}(Sp(2))$ and ${}^W \hspace{-2mm} \mathcal{A}^{\rho_z}_{\infty}(Sp(2))$, which are the second, third figures illustrated in Figure \ref{fig-A_8infty(B2)all} respectively, where multiplication by $\sigma_y = \chi_{\rho_y}|_{\mathbb{T}^2}$, $\sigma_z = \chi_{\rho_z}|_{\mathbb{T}^2}$ corresponds to the edges illustrated in the second, third diagram respectively in Figure \ref{fig-mult_rhoB2}. Both these graphs are in fact a disjoint union of two infinite graphs, coloured black, grey respectively, whose vertex sets consists of all $\lambda$ such that $\lambda_1+\lambda_2$ is even, odd respectively.
These graphs ${}^W \hspace{-2mm} \mathcal{A}^{\rho}_{\infty}(Sp(2))$ are essentially $W$-unfolded versions of the graphs $\mathcal{A}^{\rho}_{\infty}(Sp(2))$ (see Figures \ref{fig-A_infty(C2)1}-\ref{fig-A_infty(B2)2}), where $W$ denotes the Weyl group $D_8$ of $Sp(2)$.

We consider first the fixed point algebra of $\bigotimes_{\mathbb{N}}M_4$, $\bigotimes_{\mathbb{N}}M_5$, $\bigotimes_{\mathbb{N}}M_{10}$ under the conjugate action of the torus $\mathbb{T}^2$ given by the restrictions of $\rho_x$, $\rho_y$, $\rho_z$ respectively to $\mathbb{T}^2$ given in (\ref{eqn:restrict_rho1B2_to_T2}), (\ref{eqn:restrict_rho2B2_to_T2}), (\ref{eqn:restrict_rho3B2_to_T2}) respectively. Here $\mathbb{T}^2$ acts by conjugation on each factor in the infinite tensor product.
Thus by \cite[$\S$3.5]{evans/kawahigashi:1998} we have $(\bigotimes_{\mathbb{N}}M_4)^{\mathbb{T}^2} \cong A({}^W \hspace{-2mm} \mathcal{A}^{\rho_x}_{\infty}(Sp(2)))$, $(\bigotimes_{\mathbb{N}}M_5)^{\mathbb{T}^2} \cong A({}^W \hspace{-2mm} \mathcal{A}^{\rho_y}_{\infty}(Sp(2)))$ and $(\bigotimes_{\mathbb{N}}M_{10})^{\mathbb{T}^2} \cong A({}^W \hspace{-2mm} \mathcal{A}^{\rho_z}_{\infty}(Sp(2)))$.
Here $A(\mathcal{G}) = \overline{\bigcup_k A(\mathcal{G})_k}$ is the path algebra of the graph $\mathcal{G}$, where $A(\mathcal{G})_k$ is the algebra generated by pairs $(\eta_1,\eta_2)$ of paths from the distinguished vertex $\ast$ such that the ranges $r(\eta_1)$ and $r(\eta_2)$ are equal, and $|\eta_1|=|\eta_2|=k$, with multiplication defined by $(\eta_1,\eta_2) \cdot (\eta_1',\eta_2') = \delta_{\eta_2,\eta_1'}(\eta_1,\eta_2')$.

We now define commuting self-adjoint operators which may be identified with the adjacency matrix of ${}^W \hspace{-2mm} \mathcal{A}^{\rho_u}_{\infty}(Sp(2))$. We define operators $v_Z^u$ in $\ell^2(\mathbb{Z}) \otimes \ell^2(\mathbb{Z})$, for $u=x,y,x$, by
\begin{align}
v_Z^x & = s \otimes 1 + s^{\ast} \otimes 1 + 1 \otimes s + 1 \otimes s^{\ast}, \\
v_Z^y & = 1 \otimes 1 + s \otimes s + s^{\ast} \otimes s^{\ast} + s \otimes s^{\ast} + s^{\ast} \otimes s, \\
v_Z^z & = 2(1 \otimes 1) + s^2 \otimes 1 + (s^{\ast})^2 \otimes 1 + 1 \otimes s^2 + 1 \otimes (s^{\ast})^2 + s \otimes s + s^{\ast} \otimes s^{\ast} + s \otimes s^{\ast} + s^{\ast} \otimes s,
\end{align}
where $s$ is the bilateral shift on $\ell^2(\mathbb{Z})$.
Let $\Omega$ denote the vector $(\delta_{i,0})_i$. Then $v_Z^u$ is identified with the adjacency matrix of ${}^W \hspace{-2mm} \mathcal{A}^{\rho_u}_{\infty}(Sp(2))$, $u=x,y,z$, where we regard the vector $\Omega \otimes \Omega$ as corresponding to the vertex $(0,0)$ of ${}^W \hspace{-2mm} \mathcal{A}^{\rho_u}_{\infty}(Sp(2))$, and the operators of the form $s^l \otimes s^m$ which appear as terms in $v_Z^u$ as corresponding to the edges on ${}^W \hspace{-2mm} \mathcal{A}^{\rho_u}_{\infty}(Sp(2))$. Then $(s^{\lambda_1} \otimes s^{\lambda_2})(\Omega \otimes \Omega)$  corresponds to the vertex $(\lambda_1,\lambda_2)$ of ${}^W \hspace{-2mm} \mathcal{A}^{\rho_u}_{\infty}(Sp(2))$ for any $\lambda_1,\lambda_2 \in \mathbb{Z}$, and applying $(v_Z^u)^m$ to $\Omega \otimes \Omega$ gives a vector $y=(y_{(\lambda_1,\lambda_2)})$ in $\ell^2({}^W \hspace{-2mm} \mathcal{A}^{\rho_u}_{\infty}(Sp(2)))$, where $y_{(\lambda_1,\lambda_2)}$ gives the number of paths of length $m$ on ${}^W \hspace{-2mm} \mathcal{A}^{\rho_u}_{\infty}(Sp(2))$ from $(0,0)$ to the vertex $(\lambda_1,\lambda_2)$.

We define a state $\varphi$ on $C^{\ast}(v_Z^u)$ by $\varphi( \, \cdot \, ) = \langle \, \cdot \, (\Omega \otimes \Omega), \Omega \otimes \Omega \rangle$. We use the notation $(a_1,a_2,\ldots,a_k)!$ to denote the multinomial coefficient $(\sum_{i=1}^k a_i)!/\prod_{i=1}^k(a_i!)$.
Then
\begin{align*}
\varphi((v_Z^u)^m) & = \sum_{\stackrel{k_i \geq 0}{\scriptscriptstyle{\sum_i k_i \leq m}}} (k_1, \ldots, k_{l(u)}, m - \sum_i k_i)! \; \varphi(s^{r_1^u} \otimes s^{r_2^u}) \\
& = \sum_{\stackrel{k_i \geq 0}{\scriptscriptstyle{\sum_i k_i \leq m}}} (k_1, \ldots, k_{l(u)}, m - \sum_i k_i)! \; \delta_{r_1^u, 0} \; \delta_{r_2^u, 0},
\end{align*}
where $l(u) = 3,4,9$ for $u=x,y,z$ respectively, and
\begin{align}
\hspace{-5mm} r_1^x & = k_1-k_2,                    & r_2^x &= k_1+k_2+2k_3-m, \label{eqn:r1^1,r2^1-B2} \\
\hspace{-5mm} r_1^y & = k_1-k_2+k_3-k_4,            & r_2^y &= k_1-k_2-k_3+k_4, \label{eqn:r1^2,r2^2-B2} \\
\hspace{-5mm} r_1^z & = 2k_1-2k_2+k_5+k_6-k_7-k_8,  & r_2^z &= 2k_3-2k_4+k_5-k_6+k_7-k_8, \label{eqn:r1^3,r2^3-B2}
\end{align}
When $u=x$, we get a non-zero contribution when $k_2 = k_1$ and $k_3 = -k_1+m/2$. So we obtain
\begin{equation} \label{eqn:momentsA_8infty(B2)1}
\varphi((v_Z^x)^m) = \sum_{k_1} (k_1, k_1, -k_1+m/2, -k_1+m/2)!
\end{equation}
where the summation is over all integers $0 \leq k_1 \leq m/2$.
When $u=y$, we get a non-zero contribution when $k_2 = k_1$ and $k_4 = k_3$. So we obtain
\begin{equation} \label{eqn:momentsA_8infty(B2)2}
\varphi((v_Z^y)^m) = \sum_{k_1,k_3} (k_1, k_1, k_3, k_3, m-2k_1-2k_3)!
\end{equation}
where the summation is over all integers $k_1,k_3 \geq 0$ such that $2k_1+2k_3 \leq m$.
When $u=z$, we get a non-zero contribution when $k_7 = k_1-k_2+k_3-k_4+k_5$ and $k_8 = k_1-k_2-k_3+k_4+k_6$. So we obtain
\begin{equation} \label{eqn:momentsA_8infty(B2)3}
\varphi((v_Z^z)^m) = \sum_{k_i} (k_1, k_2, k_3, k_4, k_5, k_6, p_1, p_2, k_9, m-3k_1+k_2-k_3-k_4-2k_5-2k_6-k_9)!
\end{equation}
where $p_1=k_1-k_2+k_3-k_4+k_5$, $p_2=k_1-k_2-k_3+k_4+k_6$, and the summation is over all integers $k_1,k_2,\ldots,k_6, k_9 \geq 0$ such that $3k_1-k_2+k_3+k_4+2k_5+2k_6+k_9 \leq m$.

\subsection{Joint spectral measure for ${}^W \hspace{-2mm} \mathcal{A}_{\infty}(Sp(2))$, ${}^W \hspace{-2mm} \mathcal{A}_{\infty}(SO(5))$ over $\mathbb{T}^2$} \label{sect:measure:A8inftyC2-T2}

The ranges of the restrictions (\ref{eqn:restriction-chi1})-(\ref{eqn:restriction-chi3}) of the characters $\chi_{\rho_u}$ of the irreducible representations $\rho_u$ of $Sp(2)$ to $\mathbb{T}^2$, for $u=x,y,z$, are given by $I_x := \{ 2\mathrm{Re}(\omega_1) + 2\mathrm{Re}(\omega_2) | \, \omega_1,\omega_2 \in \mathbb{T} \} = [-4,4]$, $I_y := \{ 1 + 2\mathrm{Re}(\omega_1\omega_2) + 2\mathrm{Re}(\omega_1\omega_2^{-1}) | \, \omega_1,\omega_2 \in \mathbb{T} \} = [-3,5]$ and $I_z := \{ 2 + 2\mathrm{Re}(\omega_1^2) + 2\mathrm{Re}(\omega_2^2) + 2\mathrm{Re}(\omega_1\omega_2) + 2\mathrm{Re}(\omega_1\omega_2^{-1}) | \, \omega_1,\omega_2 \in \mathbb{T} \} = [-2,10]$:
\begin{align}
\chi_{\rho_x}(\omega_1,\omega_2) &= \omega_1 + \omega_1^{-1} + \omega_2 + \omega_2^{-1} = 2\cos(2\pi\theta_1) + 2\cos(2\pi\theta_2), \label{eqn:Phi_1} \\
\chi_{\rho_y}(\omega_1,\omega_2) &= 1 + \omega_1\omega_2 + \omega_1^{-1}\omega_2^{-1} + \omega_1\omega_2^{-1} + \omega_1^{-1}\omega_2 \nonumber \\
& = 1 + 2\cos(2\pi(\theta_1+\theta_2)) + 2\cos(2\pi(\theta_1-\theta_2)), \label{eqn:Phi_2} \\
\chi_{\rho_z}(\omega_1,\omega_2) &= \chi_{\rho_x}(\omega_1,\omega_2)^2-\chi_{\rho_y}(\omega_1,\omega_2)-1 \nonumber \\
& = 2 + 2\cos(4\pi\theta_1) + 2\cos(4\pi\theta_2) + 2\cos(2\pi(\theta_1+\theta_2)) + 2\cos(2\pi(\theta_1-\theta_2)), \label{eqn:Phi_3}
\end{align}
where $\omega_j = e^{2\pi i \theta_j} \in \mathbb{T}$ for $\theta_j \in [0,1]$, $j=1,2$.
We will write $x,y,z$ for the elements $\chi_{\rho_x}(\omega_1,\omega_2)$, $\chi_{\rho_y}(\omega_1,\omega_2)$, $\chi_{\rho_z}(\omega_1,\omega_2)$ respectively.
Since the spectrum $\sigma(s)$ of $s$ is $\mathbb{T}$, the spectrum $\sigma(v_Z^u)$ of $v_Z^u$ is $I_u$, $u=x,y,z$.

The Weyl group of $Sp(2)$ is the dihedral group $D_8$ of order 8.
If we consider $D_8$ as the subgroup of $GL(2,\mathbb{Z})$ generated by the matrices $T_2$, $T_4$, of orders 2, 4 respectively, given by
\begin{equation} \label{T2,T4}
T_2 = \left( \begin{array}{cc} 0 & 1 \\ 1 & 0 \end{array} \right), \qquad T_4 = \left( \begin{array}{cc} 0 & 1 \\ -1 & 0 \end{array} \right),
\end{equation}
then the action of $D_8$ on $\mathbb{T}^2$ given by $T(\omega_1,\omega_2) = (\omega_1^{a_{11}}\omega_2^{a_{12}},\omega_1^{a_{21}}\omega_2^{a_{22}})$, for $T = (a_{il}) \in D_8$, leaves $\chi_{\rho_u}(\omega_1,\omega_2)$ invariant, for $u=x,y,z$.
Then for $u=x,y,z$, any $D_8$-invariant measure $\varepsilon$ on $\mathbb{T}^2$ produces a probability measure $\mu_u$ on $I_u$ by
\begin{equation} \label{eqn:measures-T2-Ij_B2}
\int_{I_u} \psi(x) \mathrm{d}\mu_u(x) = \int_{\mathbb{T}^2} \psi(\chi_{\rho_u}(\omega_1,\omega_2)) \mathrm{d}\varepsilon(\omega_1,\omega_2),
\end{equation}
for any continuous function $\psi:I_u \rightarrow \mathbb{C}$, where $\mathrm{d}\varepsilon(\omega_1,\omega_2) = \mathrm{d}\varepsilon(g(\omega_1,\omega_2))$ for all $g \in D_8$.
There is a loss of dimension here, in the sense that the integral on the right hand side is over the two-dimensional torus $\mathbb{T}^2$, whereas the spectrum of ${}^W \hspace{-2mm} \mathcal{A}^{\rho_u}_{\infty}(Sp(2))$ is real and lives on the interval $I_u$. We introduce an intermediate probability measure $\nu$ in Section \ref{sect:measureA8inftyB2-D} which lives over the joint spectrum $\mathfrak{D}_{\lambda,\mu} \subset I_{\lambda} \times I_{\mu} \subset \mathbb{R}^2$ for irreducible representations $\lambda$, $\mu$, where there is no loss of dimension.

The spectral measure on $\mathbb{T}^2$ for the graph ${}^W \hspace{-2mm} \mathcal{A}^{\rho_u}_{\infty}(Sp(2))$ is easily seen to be the uniform Lebesgue measure $\mathrm{d}\varepsilon(\omega_1,\omega_2) = \mathrm{d}\omega_1 \, \mathrm{d}\omega_2 / 4 \pi^2$ for $u=x,y,z$, since the $m^{\mathrm{th}}$ moment is given by
\begin{align*}
\frac{1}{4\pi^2} \int_{\mathbb{T}^2} (\chi_{\rho_u}(\omega_1, \omega_2))^m \mathrm{d}\omega_1 \, \mathrm{d}\omega_2 & = \frac{1}{4\pi^2} \sum_{\stackrel{k_i \geq 0}{\scriptscriptstyle{\sum_i k_i \leq m}}} (k_1, k_2, \ldots, k_{l(u)}, m - \sum_i k_i)! \; \int_{\mathbb{T}^2} \omega_1^{r_1^u} \omega_2^{r_2^u} \mathrm{d}\omega_1 \, \mathrm{d}\omega_2 \\
& = \sum_{\stackrel{k_i \geq 0}{\scriptscriptstyle{\sum_i k_i \leq m}}} (k_1, k_2, \ldots, k_{l(u)}, m - \sum_i k_i)! \; \delta_{r_1^u, 0} \; \delta_{r_2^u, 0},
\end{align*}
where $r_1^u$, $r_2^u$ are as in (\ref{eqn:r1^1,r2^1-B2})-(\ref{eqn:r1^3,r2^3-B2}) and $l(u)=3,4,9$ for $u=x,y,z$ respectively, which is equal to $\varphi((v_Z^u)^m)$ given in (\ref{eqn:momentsA_8infty(B2)1})-(\ref{eqn:momentsA_8infty(B2)3}).

\begin{figure}[tb]
\begin{center}
  \includegraphics[width=55mm]{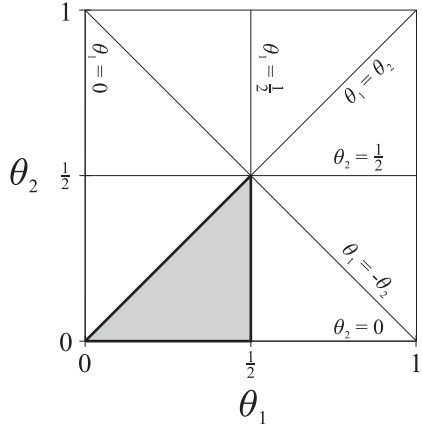}\\
 \caption{A fundamental domain $C$ of $\mathbb{T}^2/D_8$.} \label{fig:fund_domain-B2inT2}
\end{center}
\end{figure}

A fundamental domain $C$ of $\mathbb{T}^2$ under the action of the dihedral group $D_8$ is illustrated in Figure \ref{fig:fund_domain-B2inT2}, where the axes are labelled by the parameters $\theta_1$, $\theta_2$ in $(e^{2 \pi i \theta_1},e^{2 \pi i \theta_2}) \in \mathbb{T}^2$. In Figure \ref{fig:fund_domain-B2inT2}, the lines $\theta_1=0$ and $\theta_2=0$ are also boundaries of copies of the fundamental domain $C$ under the action of $D_8$. The torus $\mathbb{T}^2$ contains 8 copies of $C$, so that
\begin{equation} \label{eqn:measureT2=8C}
\int_{\mathbb{T}^2} \phi(\omega_1,\omega_2) \mathrm{d}\varepsilon(\omega_1,\omega_2) = 8 \int_{C} \phi(\omega_1,\omega_2) \mathrm{d}\varepsilon(\omega_1,\omega_2),
\end{equation}
for any $D_8$-invariant function $\phi:\mathbb{T}^2 \rightarrow \mathbb{C}$. The fixed points of $\mathbb{T}^2$ under the action of $D_8$ are the points $(1,1)$ and $(-1,-1)$, which map to the points 4, $-4$ respectively in the interval $I_x$, whilst both map to the points 5, 10 in the intervals $I_y$, $I_z$ respectively. The point $(-1,1)$ (and its orbit under $D_8$) maps to 0, $-3$, 2 in the intervals $I_x$, $I_y$, $I_z$ respectively.

%
%

\subsection{Joint spectral measure for ${}^W \hspace{-2mm} \mathcal{A}_{\infty}(Sp(2))$, ${}^W \hspace{-2mm} \mathcal{A}_{\infty}(SO(5))$ on $\mathbb{R}^2$} \label{sect:measureA8inftyB2-D}


Let $\Psi_{\lambda,\mu}$ be the map $(\omega_1,\omega_2) \mapsto (x_{\lambda},x_{\mu}) = (\chi_{\lambda}(\omega_1,\omega_2), \chi_{\mu}(\omega_1,\omega_2))$. We denote by $\mathfrak{D}_{\lambda,\mu}$ the image of $\Psi_{\lambda,\mu}(C) \; (= \Psi_{\lambda,\mu}(\mathbb{T}^2))$ in $\mathbb{R}^2$.
Note that we can identify $\mathfrak{D}_{\lambda,\mu}$ with $\mathfrak{D}_{\mu,\lambda}$ by reflecting about the line $x_{\lambda} = x_{\mu}$.
The joint spectral measure $\widetilde{\nu}_{\lambda,\mu}$ is the measure on $\mathfrak{D}_{\lambda,\mu}$ uniquely determined by its cross-moments $\varsigma_{\lambda,\mu}(m,n) = \int_{\mathfrak{D}_{\lambda,\mu}} x_{\lambda}^m x_{\mu}^n \mathrm{d}\nu_{\lambda,\mu}(x_{\lambda},x_{\mu})$.
Then there is a unique $D_8$-invariant pullback measure $\varepsilon_{\lambda,\mu}$ on $\mathbb{T}^2$ such that
\begin{equation} \label{eqn:measures-T2-D_C2}
\int_{\mathfrak{D}_{\lambda,\mu}} \psi(x_{\lambda},x_{\mu}) \mathrm{d}\widetilde{\nu}_{\lambda,\mu}(x_{\lambda},x_{\mu}) = \int_{\mathbb{T}^2} \psi(\chi_{\lambda}(\omega_1,\omega_2),\chi_{\mu}(\omega_1,\omega_2)) \mathrm{d}\varepsilon_{\lambda,\mu}(\omega_1,\omega_2),
\end{equation}
for any continuous function $\psi:\mathfrak{D}_{\lambda,\mu} \rightarrow \mathbb{C}$.

Any probability measure on $\mathfrak{D}_{\lambda,\mu}$ yields a probability measure on the interval $I_{\lambda}$, given by the pushforward $(p_{\lambda})_{\ast}(\widetilde{\nu}_{\lambda,\mu})$ of the joint spectral measure $\widetilde{\nu}_{\lambda,\mu}$ under the orthogonal projection $p_{\lambda}$ onto the spectrum $\sigma(\lambda) = I_{\lambda}$.
In particular, when $\psi(x_{\lambda},x_{\mu}) = \widetilde{\psi}(x_{\lambda})$ is only a function of one variable $x_{\lambda}$, then
$$\int_{\mathfrak{D}_{\lambda,\mu}} \widetilde{\psi}(x_{\lambda}) \mathrm{d}\widetilde{\nu}_{\lambda,\mu}(x_{\lambda},x_{\mu}) = \int_{I_{\lambda}} \widetilde{\psi}(x_{\lambda}) \int_{\mathfrak{D}_{\lambda,\mu}(x_{\lambda})} \mathrm{d}\widetilde{\nu}_{\lambda,\mu}(x_{\lambda},x_{\mu}) = \int_{I_{\lambda}} \widetilde{\psi}(x_{\lambda}) \mathrm{d}\nu_{\lambda}(x_{\lambda})$$
where the measure $\mathrm{d}\nu_{\lambda}(x_{\lambda}) = \int_{x_{\mu} \in \mathfrak{D}_{\lambda,\mu}(x_{\lambda})} \mathrm{d}\widetilde{\nu}_{\lambda,\mu}(x_{\lambda},x_{\mu})$ is given by the integral over $x_{\mu} \in \mathfrak{D}_{\lambda,\mu}(x_{\lambda}) = \{ x_{\mu} \in I_{\mu} | \, (x_{\lambda},x_{\mu}) \in \mathfrak{D}_{\lambda,\mu} \}$.
Since the spectral measure $\nu_{\lambda}$ over $I_{\lambda}$ is also uniquely determined by its (one-dimensional) moments $\widetilde{\varsigma}_{\lambda}(m) = \int_{I_{\lambda}} x_{\lambda}^m \mathrm{d}\nu_{\lambda}(x_{\lambda})$ for all $m \in \mathbb{N}$, one could alternatively consider the moments $\varsigma_{\lambda,\mu}(m,0)$ to determine the measure $\nu_{\lambda}$ over $I_{\lambda}$.

\begin{figure}[tb]
\begin{center}
  \includegraphics[width=135mm]{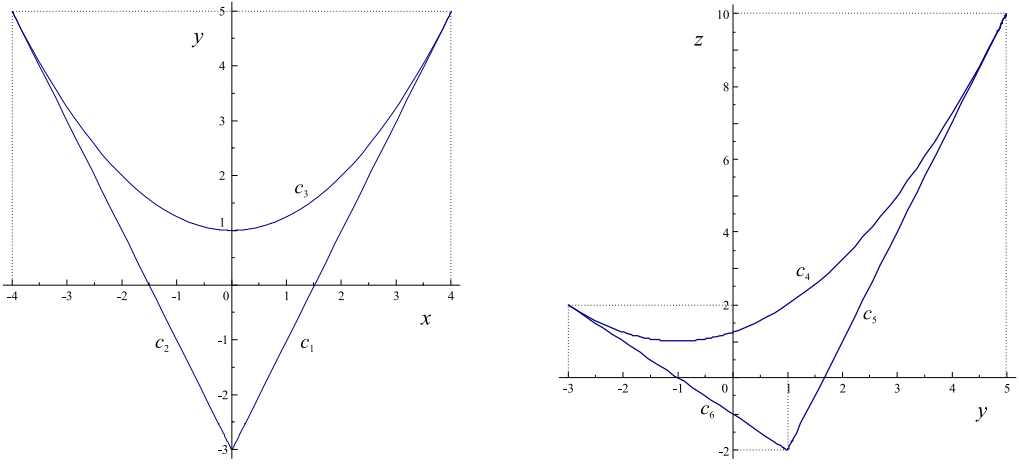}\\
 \caption{The domains $\mathfrak{D}_{x,y}$ and $\mathfrak{D}_{y,z}$ for $Sp(2)$.} \label{fig:DomainD-B2}
\end{center}
\end{figure}

In particular, we will consider the joint spectral measure of the fundamental representations $\rho_x$ and $\rho_y$ of $Sp(2)$ over $\mathfrak{D}_{x,y} := \mathfrak{D}_{\rho_x,\rho_y}$, and the joint spectral measure of the fundamental representations $\rho_y$ and $\rho_z$ of $SO(5)$ over $\mathfrak{D}_{y,z} := \mathfrak{D}_{\rho_y,\rho_z}$, illustrated in Figure \ref{fig:DomainD-B2}.

We first describe $\mathfrak{D}_{x,y}$.
The boundaries of $C$ given by $\theta_2 = 0$, $\theta_1 = 1/2$ respectively, yield the lines $c_1$, $c_2$ respectively, whilst the boundary $\theta_1 = \theta_2$ of $C$ yields the curve $c_3$. These curves are given by given by (c.f. \cite[$\S6.3$]{uhlmann/meinel/wipf:2007})
\begin{equation} \label{eqn:boundaryD-y-B2}
c_1: \;\; y=2x-3, \qquad c_2: \;\; y=-2x-3, \qquad c_3: \;\; 4y=4+x^2.
\end{equation}

For $\mathfrak{D}_{y,z}$, the boundaries of $C$ given by $\theta_2 = 0$ and $\theta_1 = 1/2$ both yield the curve $c_4$, whilst the boundary $\theta_1 = \theta_2$ of $C$ yields the line $c_5$.
Additionally, the line $\theta_2 = 1/2 - \theta_1$ which bisects $C$ yields the third boundary of $\mathfrak{D}_{y,z}$, the line $c_6$. These curves are given by
\begin{equation} \label{eqn:boundaryD-z-B2}
c_4: \;\; 4z=y^2+2y+5, \qquad c_5: \;\; z=3y-5, \qquad c_6: \;\; z=-y-1.
\end{equation}
Note that there is a two-to-one mapping from the fundamental domain $C$ to $\mathfrak{D}_{y,z}$.

\begin{figure}[tb]
\begin{minipage}[t]{8.9cm}
\begin{center}
  \includegraphics[width=80mm]{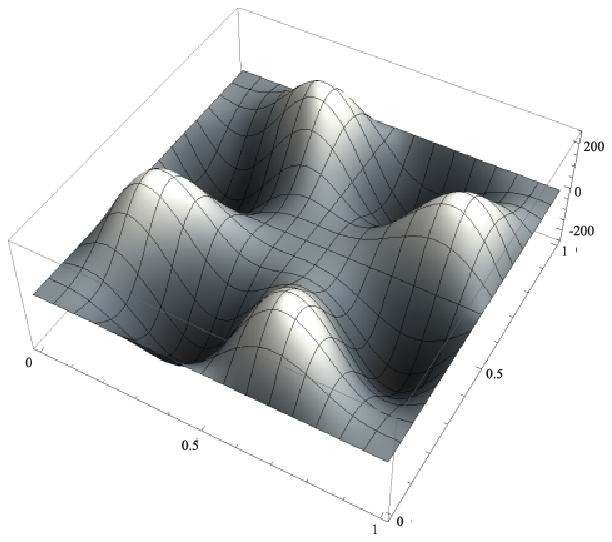}\\
 \caption{The Jacobian $J_{x,y}$ over $\mathbb{T}^2$.} \label{fig:J_overT2-3Dplot-B2}
\end{center}
\end{minipage}
\hfill
\begin{minipage}[t]{6.9cm}
\begin{center}
  \includegraphics[width=55mm]{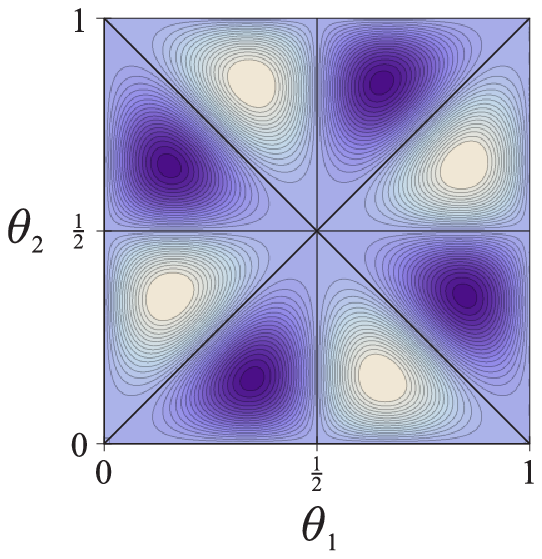}\\
 \caption{Contour plot of $J_{x,y}$ over $\mathbb{T}^2$.} \label{fig:J_overT2-contours-B2}
\end{center}
\end{minipage}
\end{figure}

Under the change of variables $x=\chi_{\rho_x}(\omega_1,\omega_2)$, $y=\chi_{\rho_y}(\omega_1,\omega_2)$, the Jacobian $J_{x,y} = \mathrm{det}(\partial(x,y)/\partial(\theta_1,\theta_2))$ is given by
\begin{align} \label{eqn:J[theta]-B2}
J_{x,y} (\theta_1,\theta_2) & = 8 \pi^2 (\cos(2 \pi (\theta_1 + 2\theta_2)) + \cos(2 \pi (2\theta_1 - \theta_2)) - \cos(2 \pi (2\theta_1 + \theta_2)) \nonumber \\
& \qquad - \cos(2 \pi (\theta_1 - 2\theta_2)).
\end{align}
The Jacobian $J_{x,y}$ is real and is illustrated in Figures \ref{fig:J_overT2-3Dplot-B2}, \ref{fig:J_overT2-contours-B2}, where its values are plotted over the torus $\mathbb{T}^2$.

With $\omega_j = e^{2 \pi i \theta_j}$, $j=1,2$, the Jacobian is given in terms of $\omega_1, \omega_2 \in \mathbb{T}$ by
\begin{align}
J_{x,y} (\omega_1,\omega_2) & = 8 \pi^2 \mathrm{Re}(\omega_1\omega_2^2 + \omega_1^2\omega_2^{-1} - \omega_1^2\omega_2 - \omega_1\omega_2^{-2}) \nonumber \\
& = 4 \pi^2 (\omega_1\omega_2^2 + \omega_1^{-1}\omega_2^{-2} + \omega_1^2\omega_2^{-1} + \omega_1^{-2}\omega_2 - \omega_1^2\omega_2 - \omega_1^{-2}\omega_2^{-1} - \omega_1\omega_2^{-2} - \omega_1^{-1}\omega_2^2). \qquad \label{eqn:J[omega]-B2}
\end{align}
The Jacobian $J_{x,y}$ is invariant under $T_4^2 \in D_8$, but $T(J_{x,y}) = -J_{x,y}$ for $T = T_2,T_4$. Thus $J_{x,y}^2$ is invariant under the action of $D_8$. An expression for $J_{x,y}^2$ in terms of the $D_8$-invariant variables $x$, $y$ may be obtained as a product of the roots appearing as the equations of the boundary of $\mathfrak{D}_{x,y}$ in (\ref{eqn:boundaryD-y-B2}), and is given as (see also \cite{uhlmann/meinel/wipf:2007})
\begin{equation} \label{J2-B2}
J_{x,y}^2(x,y) = 16 \pi^4(y+2x+3)(y-2x+3)(4y-x^2-4),
\end{equation}
for $(x,y) \in \mathfrak{D}_{x,y}$. Thus we see that the Jacobian vanishes only on the boundary of $\mathfrak{D}_{x,y}$, which is equivalent to vanishing only on the boundaries of the images of the fundamental domain in $\mathbb{T}^2$ under $D_8$.

The factorizations of $J_{x,y}$ in (\ref{J2-B2}) and the equations for the boundaries of $\mathfrak{D}_{x,y}$ given in (\ref{eqn:boundaryD-y-B2}) will be used in Sections \ref{sect:measureA8inftyB2-R}, \ref{sect:measureAinftyB2-R} to determine explicit expressions for the weights which appear in the spectral measures $\mu_{v_Z^u}$ over $I_u$ in terms of elliptic integrals.

\begin{figure}[tb]
\begin{minipage}[t]{8.9cm}
\begin{center}
  \includegraphics[width=80mm]{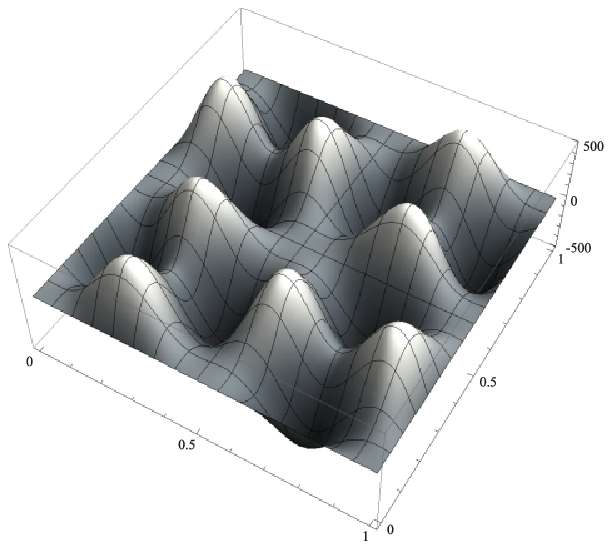}\\
 \caption{The Jacobian $J_{y,z}$ over $\mathbb{T}^2$.} \label{fig:Jyz_overT2-3Dplot-B2}
\end{center}
\end{minipage}
\hfill
\begin{minipage}[t]{6.9cm}
\begin{center}
  \includegraphics[width=55mm]{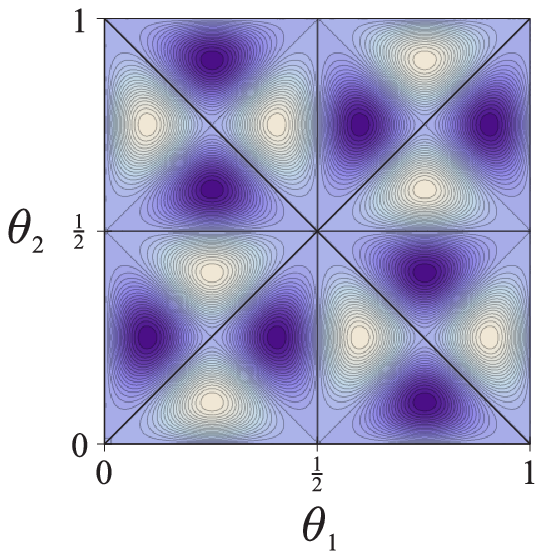}\\
 \caption{Contour plot of $J_{y,z}$ over $\mathbb{T}^2$.} \label{fig:Jyz_overT2-contours-B2}
\end{center}
\end{minipage}
\end{figure}

Similarly, under the change of variables $y=\chi_{\rho_y}(\omega_1,\omega_2)$, $z=\chi_{\rho_z}(\omega_1,\omega_2)$, the Jacobian $J_{y,z} = \mathrm{det}(\partial(y,z)/\partial(\theta_1,\theta_2))$ is given by
\begin{align} \label{eqn:Jyz[theta]-B2}
J_{y,z} (\theta_1,\theta_2) & = 16 \pi^2 (\cos(2 \pi (\theta_1 -3\theta_2)) + \cos(2 \pi (3\theta_1 + \theta_2)) - \cos(2 \pi (\theta_1 + 3\theta_2)) \nonumber \\
& \qquad - \cos(2 \pi (3\theta_1 - \theta_2)).
\end{align}
The Jacobian $J_{y,z}$ is real and is illustrated in Figures \ref{fig:Jyz_overT2-3Dplot-B2}, \ref{fig:Jyz_overT2-contours-B2}, where its values are plotted over the torus $\mathbb{T}^2$.
The Jacobian is given in terms of $\omega_1, \omega_2 \in \mathbb{T}$ by
\begin{equation}
J_{y,z} (\omega_1,\omega_2) = 8 \pi^2 (\omega_1\omega_2^{-3} + \omega_1^{-1}\omega_2^3 + \omega_1^3\omega_2 + \omega_1^{-3}\omega_2^{-1} - \omega_1\omega_2^3 - \omega_1^{-1}\omega_2^{-3} - \omega_1^3\omega_2^{-1} - \omega_1^{-3}\omega_2). \qquad \label{eqn:Jyz[omega]-B2}
\end{equation}
The Jacobian $J_{y,z}$ is again invariant under $T_4^2 \in D_8$, and $T(J_{y,z}) = -J_{y,z}$ for $T = T_2,T_4$.
An expression for $J_{y,z}^2$ in terms of the $D_8$-invariant variables $x$, $y$ may be obtained as a product of the roots appearing as the equations of the boundary of $\mathfrak{D}_{y,z}$ in (\ref{eqn:boundaryD-z-B2}), and is given as
\begin{equation} \label{Jyz2-B2}
J_{y,z}^2(y,z) = 64\pi^4(z-3y+5)(z+y+1)(y^2+2y+5-4z),
\end{equation}
for $(y,z) \in \mathfrak{D}_{y,z}$. Thus we see that the Jacobian vanishes only on the boundary of $\mathfrak{D}_{y,z}$, which is equivalent to vanishing on the boundaries of the images of the fundamental domain in $\mathbb{T}^2$ under $D_8$ as well as on the lines $\theta_2 = 1/2 \pm \theta_1$.
The lines $\theta_2 = 1/2 \pm \theta_1$ denote the lines of reflection of the additional symmetry of $\chi_{\rho_u}$, which corresponds to the fact that there is a two-to-one mapping from the fundamental domain $C$ to $\mathfrak{D}_{y,z}$.

Note that since $z=x^2-y-1$, we find that $|J_{x,y}(y,z)| = 4 \pi^2 \sqrt{(z-3y+5)(y^2+2y+5-4z)}$, and thus $J_{x,y}$ and $J_{y,z}$ are related by
$J_{y,z}(y,z)=2\sqrt{z+y+1}J_{x,y}(y,z)$.
Thus $J_{x,y}(y,z)$ is zero only on the boundaries of $\mathfrak{D}_{y,z}$ given by the curves $c_4$, $c_5$ in (\ref{eqn:boundaryD-z-B2}), but is not zero on the boundary given by $c_6$.

Since $J_{x,y}$, $J_{y,z}$ are real, $J_{x,y}^2, J_{y,z}^2 \geq 0$ and we have the following expressions:
\begin{align*}
J_{x,y} (\theta_1,\theta_2) & = 8 \pi^2 (\cos(2 \pi (\theta_1 + 2\theta_2)) + \cos(2 \pi (2\theta_1 - \theta_2)) - \cos(2 \pi (2\theta_1 + \theta_2)) \nonumber \\
& \qquad - \cos(2 \pi (\theta_1 - 2\theta_2))), \\
J_{x,y} (\omega_1,\omega_2) & = 4 \pi^2 (\omega_1\omega_2^2 + \omega_1^{-1}\omega_2^{-2} + \omega_1^2\omega_2^{-1} + \omega_1^{-2}\omega_2 - \omega_1^2\omega_2 - \omega_1^{-2}\omega_2^{-1} - \omega_1\omega_2^{-2} - \omega_1^{-1}\omega_2^2), \\
|J_{x,y} (x,y)| & = 4 \pi^2 \sqrt{(y+2x+3)(y-2x+3)(4y-x^2-4)}, \\
J_{y,z} (\theta_1,\theta_2) & = 16 \pi^2 (\cos(2 \pi (\theta_1 -3\theta_2)) + \cos(2 \pi (3\theta_1 + \theta_2)) - \cos(2 \pi (\theta_1 + 3\theta_2)) \nonumber \\
& \qquad - \cos(2 \pi (3\theta_1 - \theta_2))), \\
J_{y,z} (\omega_1,\omega_2) & = 8 \pi^2 (\omega_1\omega_2^{-3} + \omega_1^{-1}\omega_2^3 + \omega_1^3\omega_2 + \omega_1^{-3}\omega_2^{-1} - \omega_1\omega_2^3 - \omega_1^{-1}\omega_2^{-3} - \omega_1^3\omega_2^{-1} - \omega_1^{-3}\omega_2), \\
|J_{y,z} (y,z)| & = 8\pi^2 \sqrt{(z-3y+5)(z+y+1)(y^2+2y+5-4z)},
\end{align*}
where $0 \leq \theta_1,\theta_2 < 1$, $\omega_1,\omega_2 \in \mathbb{T}$ and $(x,y) \in \mathfrak{D}_{x,y}$, $(y,z) \in \mathfrak{D}_{y,z}$.

Then
\begin{align}
\int_{C} \psi(\chi_{\rho_x}(\omega_1,\omega_2),\chi_{\rho_y}(\omega_1,\omega_2)) \mathrm{d}\omega_1 \, \mathrm{d}\omega_2 &= \int_{\mathfrak{D}_{x,y}} \psi(x,y) |J_{x,y}(x,y)|^{-1} \mathrm{d}x \, \mathrm{d}y, \label{eqn:integral_overDxy-B2} \\
\int_{C} \psi(\chi_{\rho_y}(\omega_1,\omega_2),\chi_{\rho_z}(\omega_1,\omega_2)) \mathrm{d}\omega_1 \, \mathrm{d}\omega_2 &= 2\int_{\mathfrak{D}_{y,z}} \psi(y,z) |J_{y,z}(y,z)|^{-1} \mathrm{d}y \, \mathrm{d}z, \label{eqn:integral_overDyz-B2}
\end{align}
and from (\ref{eqn:measureT2=8C}) we obtain
\begin{Thm}
The joint spectral measure $\nu_{x,y}$ (over $\mathfrak{D}_{x,y}$) for ${}^W \hspace{-2mm} \mathcal{A}^{\rho_x}_{\infty}(Sp(2))$, ${}^W \hspace{-2mm} \mathcal{A}^{\rho_y}_{\infty}(Sp(2))$ is
$$\mathrm{d}\nu_{x,y}(x,y) = 8 \, |J_{x,y}(x,y)|^{-1} \mathrm{d}x \, \mathrm{d}y,$$
whilst the joint spectral measure $\nu_{y,z}$ (over $\mathfrak{D}_{y,z}$) for ${}^W \hspace{-2mm} \mathcal{A}^{\rho_y}_{\infty}(Sp(2))$, ${}^W \hspace{-2mm} \mathcal{A}^{\rho_z}_{\infty}(Sp(2))$ is
$$\mathrm{d}\nu_{y,z}(y,z) = 16 \, |J_{y,z}(y,z)|^{-1} \mathrm{d}y \, \mathrm{d}z.$$
\end{Thm}

\subsection{Spectral measure for ${}^W \hspace{-2mm} \mathcal{A}_{\infty}(Sp(2))$, ${}^W \hspace{-2mm} \mathcal{A}_{\infty}(SO(5))$ on $\mathbb{R}$} \label{sect:measureA8inftyB2-R}

We now determine the spectral measure $\mu_Z^{u,G} := \mu_{v_Z^{u,G}}$ over $I_u$, $u=x,y,z$, where $G$ is $Sp(2)$ or $SO(5)$, which is determined by its moments $\varphi((v_Z^{u,G})^m) = \int_{I_u} u^m \mathrm{d}\mu_Z^{u,G}(u)$ for all $m \in \mathbb{N}$.

\begin{figure}[tb]
\begin{minipage}[t]{5.9cm}
\begin{center}
  \includegraphics[width=51mm]{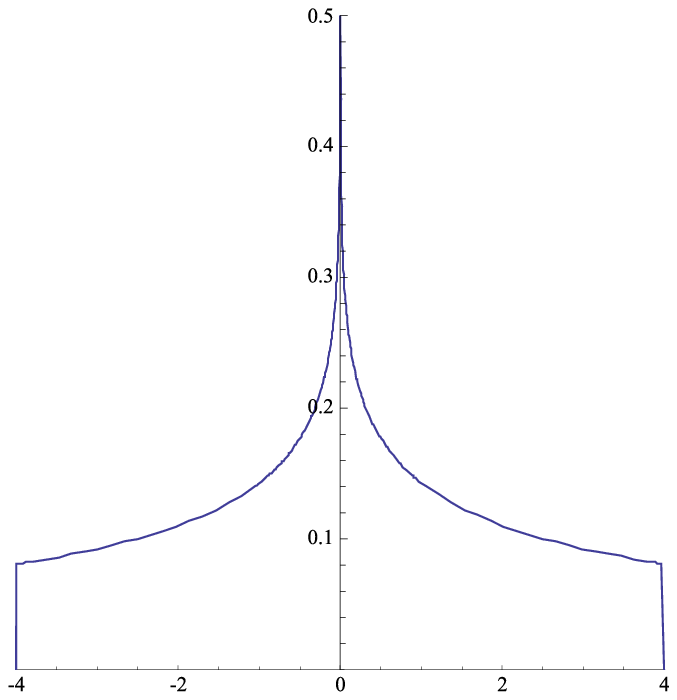}
 \caption{$J_{x}^{\mathbb{T}^2}(x)$} \label{fig-Jx-T2-B2}
\end{center}
\end{minipage}
\hfill
\begin{minipage}[t]{7.9cm}
\begin{center}
  \includegraphics[width=76mm]{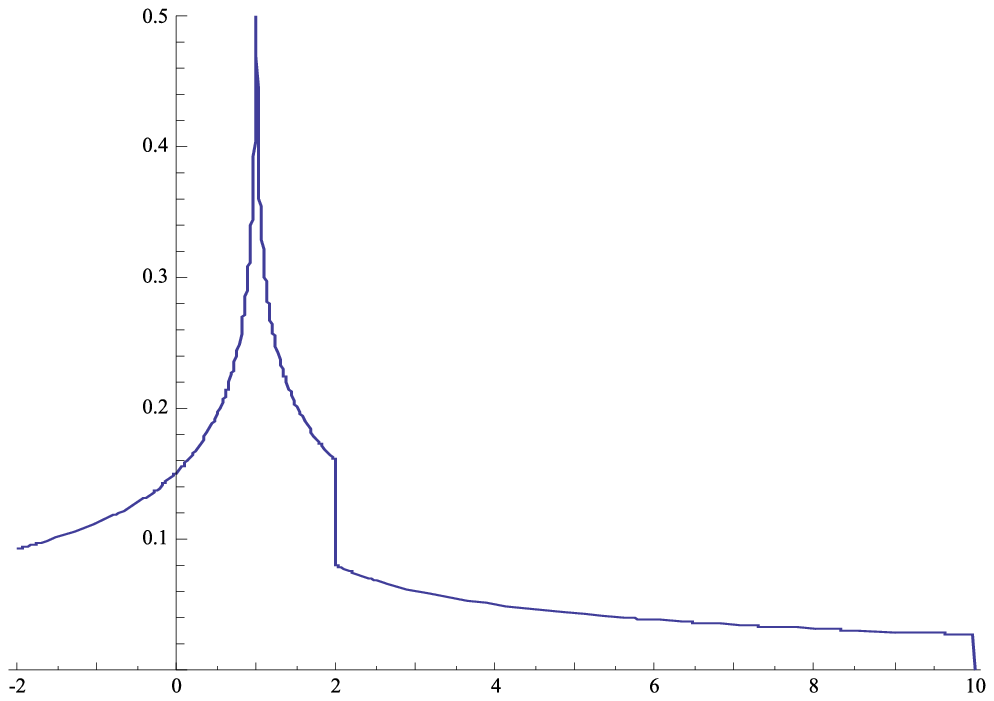}
 \caption{$J_{z}^{\mathbb{T}^2}(z)$} \label{fig-Jz-T2-B2}
\end{center}
\end{minipage}
\end{figure}

Thus for $\mu_Z^{x,Sp(2)}$ we set $\psi(x,y) = x^m$ in (\ref{eqn:integral_overDxy-B2}) and integrate with respect to $y$. Similarly, setting $\psi(x,y) = y^m$ in (\ref{eqn:integral_overDxy-B2}), the measure $\mu_Z^{y,Sp(2)}$ is obtained by integrating with respect to $x$. More explicitly, using the expressions for the boundaries of $\mathfrak{D}$ given in (\ref{eqn:boundaryD-y-B2}),
the spectral measure $\mu_Z^{x,Sp(2)}$ (over $[-4,4]$) for the graph ${}^W \hspace{-2mm} \mathcal{A}^{\rho_x}_{\infty}(Sp(2))$ is $\mathrm{d}\mu_Z^{x,Sp(2)}(x) = J_{x}^{\mathbb{T}^2}(x) \, \mathrm{d}x$, where $J_{x}^{\mathbb{T}^2}(x)$ is given by
$$J_{x}^{\mathbb{T}^2}(x) = \left\{
\begin{array}{cl}
8\int_{-2x-3}^{(x^2+4)/4} |J_{x,y}(x,y)|^{-1} \, \mathrm{d}y & \textrm{ for } x \in [-4,0], \\
8\int_{2x-3}^{(x^2+4)/4} |J_{x,y}(x,y)|^{-1} \, \mathrm{d}y & \textrm{ for } x \in [0,4].
\end{array} \right.$$
The weight $J_{x}^{\mathbb{T}^2}(x)$ is the integral of the reciprocal of the square root of a cubic in $y$, and thus can be written in terms of the complete elliptic integral $K(m)$ of the first kind, $K(m) = \int_0^{\pi/2} (1-m\sin^2\theta)^{-1/2} \mathrm{d}\theta$. Using \cite[Eqn. 235.00]{byrd/friedman:1971}, $J_{x}^{\mathbb{T}^2}(x)$ is given by
$$J_{x}^{\mathbb{T}^2}(x) = \frac{4}{\pi^2 (4-x)} \; K(v(x)) \;\; = \;\; \frac{-4 \, v(x)^{1/2}}{\pi^2 (x+4)} \; K(v(x))$$
for $x \in [-4,0]$, where $v(x) = (x+4)^2/(x-4)^2$, whilst for $x \in [0,4]$, $J_{x}^{\mathbb{T}^2}(x)$ is given by
$$J_{x}^{\mathbb{T}^2}(x) = \frac{4}{\pi^2 (x+4)} \; K(v(x)^{-1}).$$
The weight $J_{x}^{\mathbb{T}^2}(x)$ is illustrated in Figure \ref{fig-Jx-T2-B2}.

The spectral measure $\mu_Z^{y,Sp(2)}$ (over $[-3,5]$) for the graph ${}^W \hspace{-2mm} \mathcal{A}^{\rho_y}_{\infty}(Sp(2))$ is $\mathrm{d}\mu_Z^{y,Sp(2)}(y) = J_{y}^{\mathbb{T}^2}(y) \, \mathrm{d}y$, where $J_{y}^{\mathbb{T}^2}(y)$ is given by
$$J_{y}^{\mathbb{T}^2}(y) = \left\{
\begin{array}{cl}
8\int_{-(y+3)/2}^{(y+3)/2} |J_{x,y}(x,y)|^{-1} \, \mathrm{d}x & \textrm{ for } y \in [-3,1], \\
16\int_{2\sqrt{y-1}}^{(y+3)/2} |J_{x,y}(x,y)|^{-1} \, \mathrm{d}x & \textrm{ for } y \in [1,5],
\end{array} \right. $$
where the value of the square root is taken to be positive. Note that the Jacobian is an even function of $x$.
The weight $J_{y}^{\mathbb{T}^2}(y)$ is the integral of the reciprocal of the square root of a quadratic in $x^2$ and thus can also be written in terms of the complete elliptic integral of the first kind. In fact, using \cite[Eqn. 214.00]{byrd/friedman:1971} for $y \in [-3,1]$ and \cite[Eqn. 218.00]{byrd/friedman:1971} for $y \in [1,5]$, we obtain that $J_{y}^{\mathbb{T}^2}(y) = J_{x}^{\mathbb{T}^2}(y-1)$ for all $y \in [-3,5]$.
This is a surprising result, since there is no obvious symmetry between $x$ and $y$ in the Jacobian $J_{x,y}(x,y)$ -- for one thing $J_{x,y}^2$ is a quartic in $x$ but only a cubic in $y$ -- and yet the integral of $|J_{x,y}(x,y)|^{-1}$ over $x \in \mathbb{D}$ and over $y \in \mathbb{D}$ yields identical weights $J_{x}^{\mathbb{T}^2}$ and $J_{y}^{\mathbb{T}^2}$, up to a shift.

Moving to the case of $SO(5)$, the spectral measure $\mu_Z^{y,SO(5)}$ (over $[-3,5]$) for the graph ${}^W \hspace{-2mm} \mathcal{A}^{\rho_y}_{\infty}(SO(5))$ is $\mathrm{d}\mu_Z^{y,SO(5)}(y) = J_{y}^{\mathbb{T}^2}(y) \, \mathrm{d}y$, where $J_{y}^{\mathbb{T}^2}(y)$ is as above, since ${}^W \hspace{-2mm} \mathcal{A}^{\rho_y}_{\infty}(SO(5))$ is simply the connected component of $(0,0)$ in ${}^W \hspace{-2mm} \mathcal{A}^{\rho_y}_{\infty}(Sp(2))$, thus the moments $\varphi((v_Z^{y,Sp(2)})^m) = \varphi((v_Z^{y,SO(5)})^m)$.
The spectral measure $\mu_Z^{z,SO(5)}$ (over $[-2,10]$) for the graph ${}^W \hspace{-2mm} \mathcal{A}^{\rho_z}_{\infty}(SO(5))$ is $\mathrm{d}\mu_Z^{z,SO(5)}(z) = J_{z}^{\mathbb{T}^2}(z) \, \mathrm{d}z$, where $J_{z}^{\mathbb{T}^2}(z)$ is given by
$$J_{z}^{\mathbb{T}^2}(z) = \left\{
\begin{array}{cl}
16\int_{-z-1}^{(z+5)/3} |J_{y,z}(y,z)|^{-1} \, \mathrm{d}y & \textrm{ for } z \in [-2,1], \\
16\int_{-z-1}^{-1-2\sqrt{z-1}} |J_{y,z}(y,z)|^{-1} \, \mathrm{d}y + 16\int_{-1+2\sqrt{z-1}}^{(z+5)/3} |J_{y,z}(y,z)|^{-1} \, \mathrm{d}y & \textrm{ for } y \in [1,2], \\
16\int_{-1+2\sqrt{z-1}}^{(z+5)/3} |J_{y,z}(y,z)|^{-1} \, \mathrm{d}y & \textrm{ for } y \in [2,10],
\end{array} \right. $$
where the value of the square root is taken to be positive.
A numerical plot of the weight $J_{z}^{\mathbb{T}^2}(z)$ is illustrated in Figure \ref{fig-Jz-T2-B2}.

\section{Spectral measures for $\mathcal{A}_{\infty}(Sp(2))$, $\mathcal{A}_{\infty}(SO(5))$} \label{sect:measureAinftyB2}

We now consider the fixed point algebra of $\bigotimes_{\mathbb{N}}M_4$, $\bigotimes_{\mathbb{N}}M_5$ under the product action of the group $Sp(2)$ given by the fundamental representations $\rho_x$, $\rho_y$ respectively, where $Sp(2)$ acts by conjugation on each factor in the infinite tensor product, and also the fixed point algebra of $\bigotimes_{\mathbb{N}}M_5$, $\bigotimes_{\mathbb{N}}M_{10}$ under the product action of the group $SO(5)$ given by the fundamental representations $\rho_y$, $\rho_z$ respectively.

The characters $\{ \chi_{(\mu_1,\mu_2)} \}_{\mu_1,\mu_2 \in \mathbb{N}:\mu_1 \geq \mu_2}$ of $Sp(2)$ satisfy
\begin{align*}
\chi_{(1,0)} & \chi_{(\mu_1,\mu_2)} = \chi_{(\mu_1+1,\mu_2)} + \chi_{(\mu_1-1,\mu_2)} + \chi_{(\mu_1,\mu_2+1)} + \chi_{(\mu_1,\mu_2-1)}, \\
\chi_{(1,1)} & \chi_{(\mu_1,\mu_2)} = \\
& \left\{ \begin{array}{ll}
\chi_{(\mu_1+1,\mu_2+1)} + \chi_{(\mu_1-1,\mu_2-1)} + \chi_{(\mu_1+1,\mu_2-1)} + \chi_{(\mu_1-1,\mu_2+1)} & \textrm{ if } \mu_1 = \mu_2, \\
\chi_{(\mu_1,\mu_2)} + \chi_{(\mu_1+1,\mu_2+1)} + \chi_{(\mu_1-1,\mu_2-1)} + \chi_{(\mu_1+1,\mu_2-1)} + \chi_{(\mu_1-1,\mu_2+1)} & \textrm{ otherwise},
\end{array} \right. \\
\chi_{(2,0)} & \chi_{(\mu_1,\mu_2)} = \\
& \left\{ \begin{array}{ll}
\chi_{(\mu_1,\mu_2)} + \chi_{(\mu_1-2,\mu_2)} + \chi_{(\mu_1+2,\mu_2)} + \chi_{(\mu_1-1,\mu_2+1)} + \chi_{(\mu_1+1,\mu_2+1)} & \textrm{ if } \mu_2 = 0, \\
\chi_{(\mu_1,\mu_2)} + \chi_{(\mu_1+2,\mu_2)} + \chi_{(\mu_1,\mu_2-2)} + \chi_{(\mu_1+1,\mu_2-1)} & \textrm{ if } \mu_1 = \mu_2 \neq 0, \\
2\chi_{(\mu_1,\mu_2)} + \chi_{(\mu_1-2,\mu_2)} + \chi_{(\mu_1+2,\mu_2)} + \chi_{(\mu_1,\mu_2-2)} + \chi_{(\mu_1,\mu_2+2)} \\
\quad + \chi_{(\mu_1-1,\mu_2-1)} + \chi_{(\mu_1-1,\mu_2+1)} + \chi_{(\mu_1+1,\mu_2-1)} + \chi_{(\mu_1+1,\mu_2+1)} & \textrm{ otherwise},
\end{array} \right.
\end{align*}
where $\chi_{(\mu_1,\mu_2)} = 0$ if $\mu_2 < 0$ or $\mu_1 < \mu_2$.

\begin{figure}[tb]
\begin{minipage}[t]{7.9cm}
\begin{center}
  \includegraphics[width=55mm]{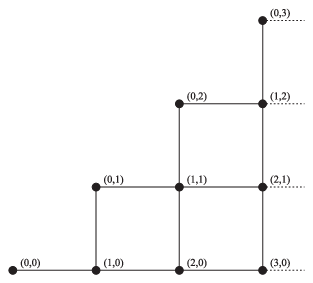}
 \caption{Infinite graph $\mathcal{A}^{\rho_x}_{\infty}(Sp(2))$} \label{fig-A_infty(C2)1}
\end{center}
\end{minipage}
\hfill
\begin{minipage}[t]{7.9cm}
\begin{center}
  \includegraphics[width=55mm]{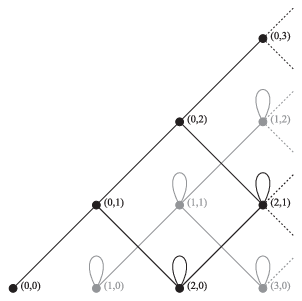}
 \caption{Infinite graph $\mathcal{A}^{\rho_y}_{\infty}(Sp(2))$} \label{fig-A_infty(C2)2}
\end{center}
\end{minipage}
\end{figure}

The representation graph of $Sp(2)$ for the first fundamental representation $\rho_x$ is identified with the infinite graph $\mathcal{A}^{\rho_x}_{\infty}(Sp(2))$, illustrated in Figure \ref{fig-A_infty(C2)1}, where we have made a change of labeling to the Dynkin labels $(\lambda_1,\lambda_2) = (\mu_1-\mu_2,\mu_2)$. This labeling is more convenient in order to be able to define self-adjoint operators $v_N^x$, $v_N^y$ in $\ell^2(\mathbb{N}) \otimes \ell^2(\mathbb{N})$ below. The dashed lines in Figure \ref{fig-A_infty(C2)1} indicate edges that are removed when one restricts to the graph $\mathcal{A}_k(Sp(2))$ at finite level $k$, c.f. Section \ref{sect:measures_AkB2}.

Similarly, the representation graph of $Sp(2)$ for the second fundamental representation $\rho_y$ is identified with the infinite graph $\mathcal{A}^{\rho_y}_{\infty}(Sp(2))$, illustrated in Figure \ref{fig-A_infty(C2)2}, again using the Dynkin labels $(\lambda_1,\lambda_2) = (\mu_1-\mu_2,\mu_2)$. Note that as with the infinite graph ${}^W \hspace{-2mm} \mathcal{A}^{\rho_y}_{\infty}(Sp(2))$, the graph $\mathcal{A}^{\rho_y}_{\infty}(Sp(2))$ is a disjoint union of two infinite graphs.

By \cite[$\S$3.5]{evans/kawahigashi:1998} we have $(\bigotimes_{\mathbb{N}}M_4)^{Sp(2)} \cong A(\mathcal{A}^{\rho_x}_{\infty}(Sp(2)))$ and $(\bigotimes_{\mathbb{N}}M_5)^{Sp(2)} \cong A(\mathcal{A}^{\rho_y}_{\infty}(Sp(2)))$.

We define self-adjoint operators $v_N^{x,Sp(2)}$, $v_N^{y,Sp(2)}$ in $\ell^2(\mathbb{N}) \otimes \ell^2(\mathbb{N})$ by
\begin{align}
v_N^{x,Sp(2)} & = l \otimes 1 + l^{\ast} \otimes 1 + l^{\ast} \otimes l + l \otimes l^{\ast}, \\
v_N^{y,Sp(2)} & = ll^{\ast} \otimes 1 + 1 \otimes l + 1 \otimes l^{\ast} + l^2 \otimes l^{\ast} + (l^{\ast})^2 \otimes l,
\end{align}
identified with the adjacency matrix of $\mathcal{A}^{\rho_u}_{\infty}(Sp(2))$, $u=x,y$, where $l$ is the unilateral shift to the right on $\ell^2(\mathbb{N})$.

Let $\Omega$ denote the vector $(\delta_{i,0})_i$.
The vector $\Omega \otimes \Omega$ is cyclic in $\ell^2(\mathbb{N}) \otimes \ell^2(\mathbb{N})$ since any vector $l^{p_1} \Omega \otimes l^{p_2} \Omega \in \ell^2(\mathbb{N}) \otimes \ell^2(\mathbb{N})$ can be written as a linear combination of elements of the form $(v_N^{x,Sp(2)})^{m_1} (v_N^{y,Sp(2)})^{m_2} (\Omega \otimes \Omega)$ so that $\overline{C^{\ast}(v_N^{x,Sp(2)}, v_N^{y,Sp(2)}) (\Omega \otimes \Omega)} = \ell^2(\mathbb{N}) \otimes \ell^2(\mathbb{N})$.
We define a state $\varphi$ on $C^{\ast}(v_N^{x,Sp(2)}, v_N^{y,Sp(2)})$ by $\varphi( \, \cdot \, ) = \langle \, \cdot \, (\Omega \otimes \Omega), \Omega \otimes \Omega \rangle$.
Since $C^{\ast}(v_N^{x,Sp(2)}, v_N^{y,Sp(2)})$ is abelian and $\Omega \otimes \Omega$ is cyclic, we have that $\varphi$ is a faithful state on $C^{\ast}(v_N^{x,Sp(2)}, v_N^{y,Sp(2)})$.

The moments $\varphi((v_N^{u,Sp(2)})^m)$ count the number of closed paths of length $m$ on the graph $\mathcal{A}^{\rho_u}_{\infty}(Sp(2))$ which start and end at the apex vertex $(0,0)$.

\begin{figure}[tb]
\begin{minipage}[t]{7.9cm}
\begin{center}
  \includegraphics[width=55mm]{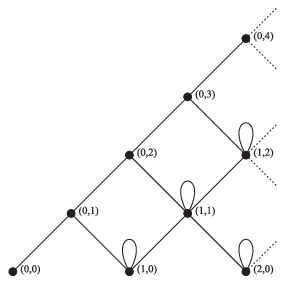}
 \caption{Infinite graph $\mathcal{A}^{\rho_y}_{\infty}(SO(5))$} \label{fig-A_infty(B2)1}
\end{center}
\end{minipage}
\hfill
\begin{minipage}[t]{7.9cm}
\begin{center}
  \includegraphics[width=55mm]{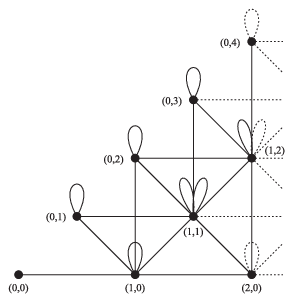}
 \caption{Infinite graph $\mathcal{A}^{\rho_z}_{\infty}(SO(5))$} \label{fig-A_infty(B2)2}
\end{center}
\end{minipage}
\end{figure}

Turning our attention to $SO(5)$, the representation graph $\mathcal{A}^{\rho_y}_{\infty}(SO(5))$ of $SO(5)$ for the first fundamental representation $\rho_y$ of $SO(5)$ is identified with the connected component of the apex vertex $(0,0)$ in the infinite graph $\mathcal{A}^{\rho_y}_{\infty}(Sp(2))$. The graph $\mathcal{A}^{\rho_y}_{\infty}(SO(5))$ is illustrated in Figure \ref{fig-A_infty(B2)1}, where we now use the Dynkin labels for $SO(5)$, $(\lambda_1,\lambda_2) = ((\mu_1-\mu_2)/2,\mu_2)$, where $(\mu_1,\mu_2)$ label the irreducible representations of $Sp(2)$ as in Section \ref{sect:measure:A8infty(B2)}.
The representation graph of $SO(5)$ for the second fundamental representation $\rho_z$ is identified with the infinite graph $\mathcal{A}^{\rho_z}_{\infty}(SO(5))$, illustrated in Figure \ref{fig-A_infty(B2)2}, again using the Dynkin labels for $SO(5)$.
Again, the dashed lines in Figures \ref{fig-A_infty(B2)1}, \ref{fig-A_infty(B2)2} indicate edges that are removed when one restricts to the graph $\mathcal{A}_k(SO(5))$ at finite level $k$, c.f. Section \ref{sect:measures_AkB2}.

Again we have $(\bigotimes_{\mathbb{N}}M_5)^{SO(5)} \cong A(\mathcal{A}^{\rho_y}_{\infty}(SO(5)))$ and $(\bigotimes_{\mathbb{N}}M_5)^{SO(5)} \cong A(\mathcal{A}^{\rho_z}_{\infty}(SO(5)))$.\footnote{Is this actually correct? If the fixed point algebra doesn't see the centre, then it should instead be $(\bigotimes_{\mathbb{N}}M_5)^{SO(5)} \cong A(\mathcal{A}^{\rho_y}_{\infty}(Sp(2)))$ and $(\bigotimes_{\mathbb{N}}M_5)^{SO(5)} \cong A(\mathcal{A}^{\rho_z}_{\infty}(Sp(2)))$?}
We define self-adjoint operators $v_N^{y,SO(5)}$, $v_N^{z,SO(5)}$ in $\ell^2(\mathbb{N}) \otimes \ell^2(\mathbb{N})$ by
\begin{align}
v_N^{y,SO(5)} & = ll^{\ast} \otimes 1 + 1 \otimes l + 1 \otimes l^{\ast} + l \otimes l^{\ast} + l^{\ast} \otimes l, \\
v_N^{z,SO(5)} & = ll^{\ast} \otimes 1 + 1 \otimes ll^{\ast} + l \otimes 1 + l^{\ast} \otimes 1 + l \otimes l^{\ast} + l^{\ast} \otimes l + ll^{\ast} \otimes l + ll^{\ast} \otimes l^{\ast} + l^{\ast} \otimes l^2 + l \otimes (l^{\ast})^2,
\end{align}
identified with the adjacency matrix of $\mathcal{A}^{\rho_u}_{\infty}(SO(5))$, $u=y,z$.

\subsection{Joint spectral measure for $\mathcal{A}_{\infty}(Sp(2))$, $\mathcal{A}_{\infty}(SO(5))$ over $\mathbb{T}^2$} \label{sect:measureAinftyB2-T2}

We will prove in Section \ref{sect:measures_AkB2} that the joint spectral measure over $\mathbb{T}^2$ of $v_N^{x,Sp(2)}$, $v_N^{y,Sp(2)}$ is the measure $\varepsilon$ given by
$$\mathrm{d}\varepsilon(\omega_1,\omega_2) = \frac{1}{128 \pi^4} J_{x,y}(\omega_1,\omega_2)^2 \mathrm{d}\omega_1 \, \mathrm{d}\omega_2,$$
where $\mathrm{d}\omega_l$ is the uniform Lebesgue measure on $\mathbb{T}$, $l=1,2$, and that the joint spectral measure over $\mathbb{T}^2$ of $v_N^{y,SO(5)}$, $v_N^{z,SO(5)}$ is also $\varepsilon$.

\subsection{Spectral measure for $\mathcal{A}_{\infty}(Sp(2))$ on $\mathbb{R}$} \label{sect:measureAinftyB2-R}

We now determine the spectral measure $\mu_N^{u,G} := \mu_{v_N^{u,G}}$ over $I_u$, where $u=x,y$ for $G=Sp(2)$ and $u=y,z$ for $G=SO(5)$.
We first consider the case of $Sp(2)$.
From (\ref{eqn:measureT2=8C}) and (\ref{eqn:integral_overDxy-B2}), with the measure given in Section \ref{sect:measureAinftyB2-T2}, we have that
\begin{equation} \label{eqn:integral_overDj-B2-2}
\frac{1}{128 \pi^4} \int_{C} \psi(\chi_{\rho_u}(\omega_1,\omega_2)) J_{x,y}(\omega_1,\omega_2)^2 \mathrm{d}\omega_1 \, \mathrm{d}\omega_2 = \frac{1}{16 \pi^4} \int_{\mathfrak{D}_{x,y}} \psi(u) |J_{x,y}(x,y)| \mathrm{d}x \, \mathrm{d}y,
\end{equation}
where $C$ is a fundamental domain of $\mathbb{T}^2/D_8$ and $\mathfrak{D}_{x,y}$ is as in Section \ref{sect:measureA8inftyB2-R}.
Thus the joint spectral measure over $\mathfrak{D}_{x,y}$ is $|J_{x,y}(x,y)| \mathrm{d}x \, \mathrm{d}y/16\pi^4$, which is the reduced Haar measure on $Sp(2)$ \cite[$\S$6.2]{uhlmann/meinel/wipf:2007}.
The measure $\mu_N^{x,Sp(2)}$ over $I_x$ is obtained by integrating with respect to $y$ in (\ref{eqn:integral_overDj-B2-2}), whilst the measure $\mu_N^{y,Sp(2)}$ over $I_y$ is obtained by integrating with respect to $x$ in (\ref{eqn:integral_overDj-B2-2}).
More explicitly, using the expressions for the boundaries of $\mathfrak{D}_{x,y}$ given in (\ref{eqn:boundaryD-y-B2}), the spectral measure $\mu_N^{x,Sp(2)}$ (over $[-4,4]$) for the graph $\mathcal{A}^{\rho_x}_{\infty}(Sp(2))$ is $\mathrm{d}\mu_N^{x,Sp(2)}(x) = J_x^{Sp(2)}(x) \, \mathrm{d}x / 16 \pi^4$, where $J_x^{Sp(2)}(x)$ is given by
\begin{align*}
\int_{-2x-3}^{(x^2+4)/4} |J_{x,y}(x,y)| \, \mathrm{d}y & \quad \textrm{ for } x \in [-4,0], \\
\int_{2x-3}^{(x^2+4)/4} |J_{x,y}(x,y)| \, \mathrm{d}y & \quad \textrm{ for } x \in [0,4].
\end{align*}
The spectral measure $\mu_N^{y,Sp(2)}$ (over $[-3,5]$) for the graph $\mathcal{A}^{\rho_y}_{\infty}(Sp(2))$ is $\mathrm{d}\mu_N^{y,Sp(2)}(y) = J_y^{Sp(2)}(y) \, \mathrm{d}y / 16 \pi^4$, where $J_y^{Sp(2)}(y)$ is given by
\begin{align*}
\int_{-(y+3)/2}^{(y+3)/2} |J_{x,y}(x,y)| \, \mathrm{d}x \hspace{35mm} & \quad \textrm{ for } y \in [-3,1], \\
\int_{-(y+3)/2}^{-2\sqrt{y-1}} |J_{x,y}(x,y)| \, \mathrm{d}x + \int_{2\sqrt{y-1}}^{(y+3)/2} |J_{x,y}(x,y)| \, \mathrm{d}x = 2\int_{2\sqrt{y-1}}^{(y+3)/2} |J_{x,y}(x,y)| \, \mathrm{d}x & \quad \textrm{ for } y \in [1,5],
\end{align*}

\begin{figure}[tb]
\begin{minipage}[t]{7.5cm}
\begin{center}
  \includegraphics[width=70mm]{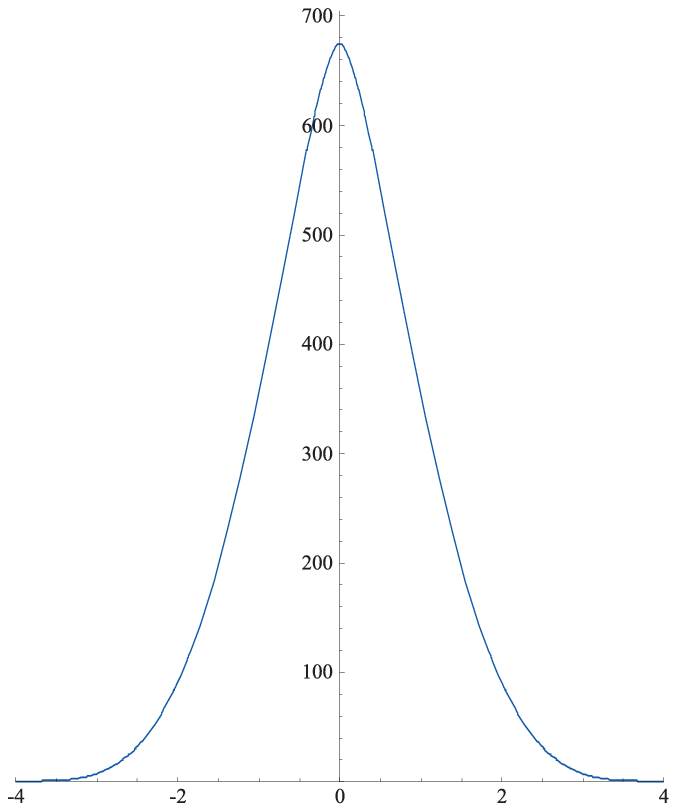}
 \caption{$J_x^{Sp(2)}(x)$} \label{fig-Jx-C2}
\end{center}
\end{minipage}
\hfill
\begin{minipage}[t]{7.5cm}
\begin{center}
  \includegraphics[width=70mm]{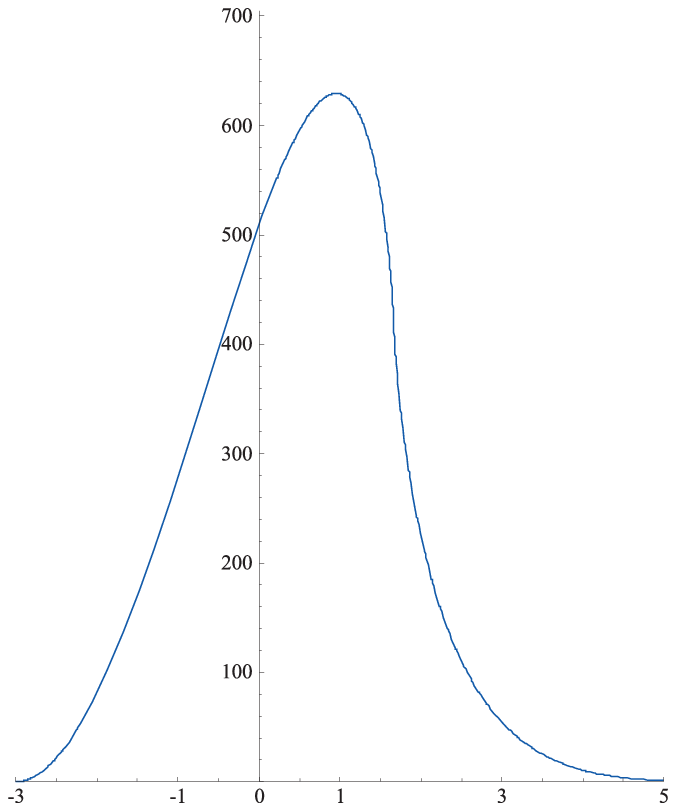}
 \caption{$J_y^{Sp(2)}(y)$} \label{fig-Jy-C2}
\end{center}
\end{minipage}
\end{figure}

The weight $J_x^{Sp(2)}(x)$ is the integral of the square root of a cubic in $y$, and thus can be written in terms of the complete elliptic integrals $K(m)$, $E(m)$ of the first, second kind respectively, where $K(m) = \int_0^{\pi/2} (1-m\sin^2\theta)^{-1/2} \mathrm{d}\theta$ and $E(m) = \int_0^{\pi/2} (1-m\sin^2\theta)^{1/2} \mathrm{d}\theta$. Using \cite[Eqn. 235.14]{byrd/friedman:1971}, $J_x^{Sp(2)}(x)$ is given by
$$\frac{\pi^2}{15} (4-x) \bigg[ (x^4 + 224 x^2 + 256) \; E(v(x)) + 8x(x^2 - 24x + 12) \; K(v(x)) \bigg],$$
for $x \in [-4,0]$, where $v(x) = (x+4)^2/(x-4)^2$, whilst for $x \in [0,4]$, $J_x^{Sp(2)}(x)$ is given by
$$\frac{\pi^2}{15} (x+4) \bigg[ (x^4 + 224 x^2 + 256) \; E(v(x)^{-1}) - 8x(x^2 + 24x + 12) \; K(v(x)^{-1}) \bigg],$$
The weight $J_x^{Sp(2)}(x)$ is illustrated in Figure \ref{fig-Jx-C2}.

Similarly, the weight $J_y^{Sp(2)}(y)$ is the integral of the square root of a quadratic in $x^2$, and can also be written in terms of the complete elliptic integrals of the first and second kinds. Using \cite[Eqn. 214.12]{byrd/friedman:1971}, $J_y^{Sp(2)}(y)$ is given by
$$\frac{2\pi^2}{3} (5-y) \bigg[ 16(1-y) \; K(v(y-1)) +(y^2+22y-7) \; E(v(y-1)) \bigg],$$
for $y \in [-3,1]$, whilst for $y \in [1,5]$, $J_y^{Sp(2)}(y)$ is given by
$$\frac{2\pi^2}{3} (y+3) \bigg[ 32(1-y) \; K(v(y-1)^{-1}) +(y^2+22y-7) \; E(v(y-1)^{-1}) \bigg],$$
using \cite[Eqn. 217.09]{byrd/friedman:1971}.
The weight $J_y^{Sp(2)}(y)$ is illustrated in Figure \ref{fig-Jy-C2}.

\begin{figure}[tb]
\begin{center}
  \includegraphics[width=105mm]{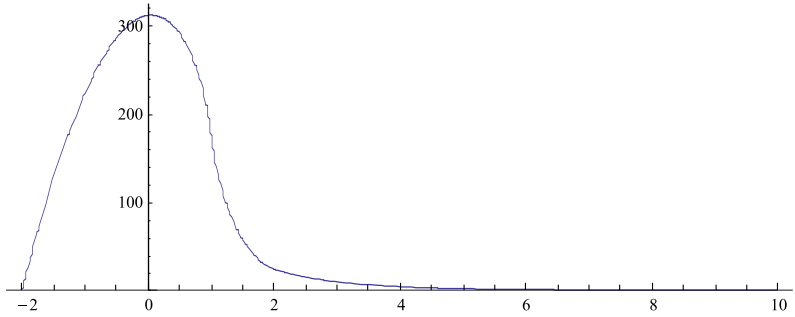}
 \caption{$J_z^{Sp(2)}(z)$} \label{fig-Jz-B2}
\end{center}
\end{figure}

We now consider the measures for $SO(5)$.
From (\ref{eqn:measureT2=8C}) and (\ref{eqn:integral_overDyz-B2}), with the measure given in Section \ref{sect:measureAinftyB2-T2}, we have that
\begin{align*}
\frac{1}{128 \pi^4} \int_{C} \psi(\chi_{\rho_u}(\omega_1,\omega_2)) J_{x,y}(\omega_1,\omega_2)^2 \mathrm{d}\omega_1 \, \mathrm{d}\omega_2 &= \frac{1}{8 \pi^4} \int_{\mathfrak{D}_{y,z}} \psi(u) J_{x,y}(y,z)^2 \, |J_{y,z}(y,z)|^{-1} \, \mathrm{d}y \, \mathrm{d}z \\
&= \frac{1}{16 \pi^4} \int_{\mathfrak{D}_{y,z}} \psi(u) |J_{x,y}(y,z)| \, (z+y+1)^{-1/2} \, \mathrm{d}y \, \mathrm{d}z.
\end{align*}
The spectral measure $\mu_N^{y,SO(5)}$ (over $[-3,5]$) for the graph $\mathcal{A}^{\rho_y}_{\infty}(SO(5))$ is $\mathrm{d}\mu_N^{y,SO(5)}(y) = J_{y}^{Sp(2)}(y) \, \mathrm{d}y/ 16 \pi^4$, where $J_{y}^{Sp(2)}(y)$ is as above.
The spectral measure $\mu_N^{z,SO(5)}$ (over $[-2,10]$) for the graph $\mathcal{A}^{\rho_z}_{\infty}(SO(5))$ is $\mathrm{d}\mu_N^{z,SO(5)}(z) = J_{z}^{Sp(2)}(z) \, \mathrm{d}z/ 8 \pi^4$, where $J_{z}^{Sp(2)}(z)$ is given by
\begin{align*}
\int_{-z-1}^{(z+5)/3} |J_{x,y}(y,z)| \, (z+y+1)^{-1/2} \, \mathrm{d}y \hspace{35mm} & \quad \textrm{ for } z \in [-2,1], \\
\int_{-z-1}^{-1-2\sqrt{z-1}} |J_{x,y}(y,z)| \, (z+y+1)^{-1/2} \, \mathrm{d}y + \int_{-1+2\sqrt{z-1}}^{(z+5)/3} |J_{x,y}(y,z)| \, (z+y+1)^{-1/2} \, \mathrm{d}y & \quad \textrm{ for } z \in [1,2], \\
\int_{-1+2\sqrt{z-1}}^{(z+5)/3} |J_{x,y}(y,z)| \, (z+y+1)^{-1/2} \, \mathrm{d}y \hspace{35mm} & \quad \textrm{ for } z \in [2,10],
\end{align*}
where the value of the square root is taken to be positive.
A numerical plot of the weight $J_{z}^{Sp(2)}(z)$ is illustrated in Figure \ref{fig-Jz-B2}.

\section{Spectral Measures for Nimrep Graphs associated to $Sp(2)$ and $SO(5)$ Modular Invariants} \label{sect:measures_nimrepB2}

We now determine joint spectral measures for nimrep graphs associated to known $Sp(2)$ modular invariants, where we will focus in particular on the nimrep graphs for the fundamental generators $\rho_j$, $j=1,2$,which have quantum dimensions $[2][6]/[3]$, $[5][6]/[2][3]$ respectively, where $[m]$ denotes the quantum integer $[m] = (q^m-q^{-m})/(q-q^{-1})$ for $q=e^{i\pi/2(k+3)}$.
The nimrep graphs $G_{\rho_j}$ were found in \cite{coquereaux/rais/tahri:2010} for the conformal embeddings at levels 3, 7, 12.
We also determine joint spectral measures for the nimrep graphs associated to the trivial $SO(5)$ modular invariants. It is not clear to us which of the other $Sp(2)$ invariants restrict to $SO(5)$ invariants in some way, or indeed whether there are other $SO(5)$ modular invariants for which there are no $Sp(2)$ counterpart.
The realisation of modular invariants for $Sp(2)$ and $SO(5)$ by braided subfactors is parallel to the realisation of $SU(2)$ and $SU(3)$ modular invariants by $\alpha$-induction for a suitable braided subfactors \cite{ocneanu:2000i, ocneanu:2002, xu:1998, bockenhauer/evans:1999i, bockenhauer/evans:1999ii, bockenhauer/evans/kawahigashi:1999, bockenhauer/evans/kawahigashi:2000}, \cite{ocneanu:2000ii, ocneanu:2002, xu:1998, bockenhauer/evans:1999i, bockenhauer/evans:1999ii, bockenhauer/evans/kawahigashi:1999, bockenhauer/evans:2001, bockenhauer/evans:2002, evans/pugh:2009i, evans/pugh:2009ii} respectively. The realisation of modular invariants for $G_2$ is also under way \cite{evans/pugh:2012i}.

Let $G$ be the nimrep associated to a braided subfactor $N \subset M$. Then the graphs $G_{\lambda}$, $\lambda \in {}_N \mathcal{X}_N$ are finite (undirected) graphs which share the same set of vertices ${}_N \mathcal{X}_M$. Their adjacency matrices (which we also denote by $G_{\lambda}$) are clearly self-adjoint. The $m,n^{\mathrm{th}}$ moment $\int_{\mathfrak{D}_{\lambda,\zeta}} x_{\lambda}^m x_{\zeta}^n \mathrm{d}\mu_{\lambda,\zeta}(x_{\lambda},x_{\zeta})$ is given by $\langle G_{\lambda}^m G_{\zeta}^n e_1, e_1 \rangle$, where $e_1$ is the basis vector in $\ell^2(G_{\lambda}) \;\; (=\ell^2(G_{\zeta}))$ corresponding to the distinguished vertex $\ast$ of $G_{\lambda}$ with lowest Perron-Frobenius weight.

Let $\beta^{\nu}_{\lambda}$ be the eigenvalues of $G_{\lambda}$, indexed by $\nu \in \mathrm{Exp}(G)$, which are ratios of the $S$-matrix given by $\beta^{\nu}_{\lambda} = S_{\lambda \nu}/S_{0 \nu}$, with corresponding eigenvectors $(\psi^{\nu}_{\zeta})_{\zeta \in \mathrm{Exp}(G)}$ (note that as the nimreps are a family of commuting matrices they can be simultaneously diagonalised, and thus the eigenvectors of $G_{\lambda}$ are the same for all $\lambda$). Then $G_{\lambda}^m G_{\zeta}^n = \mathcal{U} \Lambda_{\lambda}^m \Lambda_{\zeta}^n \mathcal{U}^{\ast}$, where $\Lambda_{\nu}$ is the diagonal matrix $\Lambda_{\lambda} = \mathrm{diag}(\beta^{\nu_1}_{\lambda}, \beta^{\nu_2}_{\lambda}, \ldots, \beta^{\nu_s}_{\lambda})$ and $\mathcal{U}$ is the unitary matrix $\mathcal{U} = (\psi^{\nu_1}, \psi^{\nu_s}, \ldots, \psi^{\nu_s})$, for $\nu_i \in \mathrm{Exp}(G)$, so that
\begin{align}
\int_{\mathbb{T}^2} (\chi_{\lambda}(\omega_1,\omega_2))^m (\chi_{\zeta}(\omega_1,\omega_2))^n \mathrm{d}\varepsilon_{\lambda,\zeta}(\omega_1,\omega_2) & = \langle \mathcal{U} \Lambda_{\lambda}^m \Lambda_{\zeta}^n \mathcal{U}^{\ast} e_1, e_1 \rangle = \langle \Lambda_{\lambda}^m \Lambda_{\zeta}^n \mathcal{U}^{\ast} e_1, \mathcal{U}^{\ast} e_1 \rangle \nonumber \\
& = \sum_{\nu \in \mathrm{Exp}(G)} (\beta^{\nu}_{\lambda})^m (\beta^{\nu}_{\zeta})^n |\psi^{\nu}_{\ast}|^2. \label{eqn:moments-nimrep-B2}
\end{align}

We denote by ${}_N \mathcal{X}_N^{SO(5)}$ the even subsystem of ${}_N \mathcal{X}_N$ given by the set of all $\lambda \in {}_N \mathcal{X}_N$ which are in the connected component of $\lambda_{(0,0)}$ in $G_{\lambda_{(0,1)}}$, i.e. the Verlinde algebra for $SO(5)$ \footnote{Is this correct? Or is the Verlinde algebra for $SO(5)$ some ($\mathbb{Z}_2$-quotient of the even subsystem? (c.f. the case of $SO(3)$)}. The eigenvalues $\beta^{\nu}_{\lambda}$ of $G_{\lambda}$ for $\lambda \in {}_N \mathcal{X}_N^{SO(5)}$ satisfy $\beta^{\lambda_{(i,j)}}_{\lambda} = \beta^{\lambda_{(i,k-i-j)}}_{\lambda}$. Thus the eigenvalues of the connected component of $\lambda_{(0,0)}$ in $G_{\lambda}$, $\lambda \in {}_N \mathcal{X}_N^{SO(5)}$, are given by the ratio $\beta^{\nu}_{\lambda} = S_{\lambda \nu}/S_{0 \nu}$ as above, where now the exponents $\nu$ essentially correspond to a fundamental domain of the exponents of $G$ under the $\mathbb{Z}_2$ action $\lambda_{(i,j)} \leftrightarrow \lambda_{(i,k-i-j)}$ (in the case of the $\mathcal{D}_{2l+1}$ the exponents which are fixed under this $\mathbb{Z}_2$ action are excluded, whilst for $\mathcal{E}_{12}$ only two copies of the exponent $(4,4)$, which is fixed under this $\mathbb{Z}_2$ action, are included). One could then conjecture what the corresponding partition functions for $SO(5)$ would be, but it would need to be verified that these are modular invariant.

The following $D_8$-invariant measure on $\mathbb{T}^2$ will be useful in what follows.
\begin{Def} \label{def:B2measure}
We denote by $\mathrm{d}^{(\theta_1,\theta_2)}$ the uniform Dirac measure on the $D_8$-orbit of the point $(e^{2\pi i\theta_1},e^{2\pi i\theta_2}), (e^{\pi i (2-\theta_2)},e^{\pi i (2-\theta_1)}) \in C \subset \mathbb{T}^2$.
\end{Def}

\begin{figure}[tb]
\begin{center}
  \includegraphics[width=55mm]{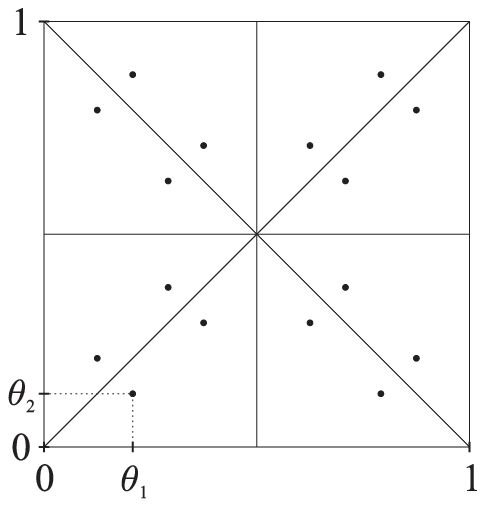}\\
 \caption{$\mathrm{Supp}(\mathrm{d}^{(\theta_1,\theta_2)})$} \label{fig:Supp-dB2}
\end{center}
\end{figure}

The set of points $(\theta_1,\theta_2) \in [0,1]^2$ such that $(e^{2 \pi i \theta_1}, e^{2 \pi i \theta_2})$ is in the support of the measure $\mathrm{d}^{(\theta_1,\theta_2)}$ is illustrated in Figure \ref{fig:Supp-dB2}.
For $(\theta_1,\theta_2) \not \in \partial C$, $|\mathrm{Supp}(\mathrm{d}^{(\theta_1,\theta_2)})| = 16$, whilst $|\mathrm{Supp}(\mathrm{d}^{(0,0)})| = |\mathrm{Supp}(\mathrm{d}^{(1/2,1/2)})| = 2$ and $|\mathrm{Supp}(\mathrm{d}^{(\theta_1,1/2-\theta_1)})| = 8$. For all other $(\theta_1,\theta_2) \in \partial C$, $|\mathrm{Supp}(\mathrm{d}^{(\theta_1,\theta_2)})| = 4$.

\subsection{Graphs $\mathcal{A}_k(Sp(2))$, $k \leq \infty$} \label{sect:measures_AkB2}

The graphs $\mathcal{A}^{\rho_u}_k(Sp(2))$, $u=x,y$, are associated to the trivial inclusion $N \rightarrow N$, and are the trivial nimrep graphs $\mathcal{G}_{\lambda} = N_{\lambda}$, where $\lambda \in {}_N \mathcal{X}_N$ for $Sp(2)$ at level $k$. The graphs $\mathcal{A}^{\rho_u}_k(Sp(2))$ are illustrated in Figures \ref{fig-A_infty(C2)1}, \ref{fig-A_infty(C2)2}, where the set of vertices ${}_N \mathcal{X}_N = P^{k,Sp(2)}_{+} := \{ (\lambda_1,\lambda_2) | \, \lambda_1,\lambda_2 \geq 0; \lambda_1 + \lambda_2 \leq k \}$, and the set of edges is given by the edges between these vertices.  Similarly, $\mathcal{A}^{\rho_u}_k(SO(5))$, $u=y,z$, associated to the trivial inclusion $N \rightarrow N$ are the trivial nimrep graphs $\mathcal{G}_{\lambda} = N_{\lambda}$, where now $\lambda \in {}_N \mathcal{X}_N$ for $SO(5)$ at level $k$, ${}_N \mathcal{X}_N = P^{k,SO(5)}_{+} := \{ (\lambda_1,\lambda_2) | \, \lambda_1,\lambda_2 \geq 0; 2\lambda_1 + \lambda_2 \leq k \}$. The graphs $\mathcal{A}^{\rho_u}_k(SO(5))$ are illustrated in Figures \ref{fig-A_infty(B2)1}, \ref{fig-A_infty(B2)2}.

The eigenvalues $\beta^{\nu,G}_{\rho_u}$ of $\mathcal{A}^{\rho_u}_k(G)$, where $u=x,y$ for $G=Sp(2)$ and $u=y,z$ for $G=SO(5)$, are given by the ratio $S_{\rho_u\nu}/S_{0\nu}$ with corresponding eigenvectors $\psi^{\nu}_{\mu} = S_{\nu,\mu}$ for $\mathcal{A}_k(Sp(2))$ with exponents $\mathrm{Exp}(\mathcal{A}_k(Sp(2))) = P^{k,Sp(2)}_{+}$, and $\psi^{\nu}_{\mu} = \sqrt{2}S_{\nu,\mu}$ for $\mathcal{A}_k(SO(5))$ with exponents $\mathrm{Exp}(\mathcal{A}_k(SO(5))) = \{ \lambda \in \mathrm{Exp}(\mathcal{A}_k(Sp(2))) | \, \lambda_1+2\lambda_2 \leq k \}$. Note that the exponents for $SO(5)$ are not the exponents for $Sp(2)$ which correspond to the labels of the irreducible representations of $(Sp(2))_k$ which are also representations of $(SO(5))_k$, namely $\lambda \in \mathrm{Exp}(\mathcal{A}_k(Sp(2)))$ such that $\lambda_1$ is even. Rather, the exponents for $SO(5)$ are a fundamental domain of the exponents for $Sp(2)$ under the $\mathbb{Z}_2$ action which maps $(\lambda_1,\lambda_2) \leftrightarrow (\lambda_1,k-\lambda_1-\lambda_2)$.
The eigenvalues $\beta^{\lambda,G}_{\rho_u}$, $\lambda \in \mathrm{Exp}(\mathcal{A}_k(G))$, are given by $\beta^{\lambda,G}_{\rho_u} = \chi_{\rho_u}(\omega_1,\omega_2)$, where $\omega_j=\exp^{2\pi i \theta_j}$, $j=1,2$ are related to $\lambda \in \mathrm{Exp}(\mathcal{A}_k(Sp(2)))$ by
\begin{equation} \label{eqn:theta-lambda_C2}
\theta_1 = \hat{\lambda}_2/2\kappa, \quad \theta_2 = (\hat{\lambda}_1 + \hat{\lambda}_2)/2\kappa \qquad \Leftrightarrow \qquad \hat{\lambda}_1 = 2\kappa(\theta_2-\theta_1), \quad \hat{\lambda}_2 = 2\kappa\theta_1.
\end{equation}

The $S$-matrix at level $k$, indexed by $\lambda \in P^{k,Sp(2)}_{+}$, is given by \cite{gannon:2001}:
\begin{align*}
S_{\lambda,\mu} & = \frac{1}{\kappa} \bigg[
\cos(\xi((\hat{\lambda}_1+2\hat{\lambda}_2)(\hat{\mu}_1+2\hat{\mu}_2) + \hat{\lambda}_1\hat{\mu}_1 )) - \cos(\xi((\hat{\lambda}_1+2\hat{\lambda}_2)(\hat{\mu}_1+2\hat{\mu}_2) - \hat{\lambda}_1\hat{\mu}_1 )) \\
& \qquad + \cos(\xi((\hat{\lambda}_1+2\hat{\lambda}_2)\hat{\mu}_1 - \hat{\lambda}_1(\hat{\mu}_1+2\hat{\mu}_2))) - \cos(\xi((\hat{\lambda}_1+2\hat{\lambda}_2)\hat{\mu}_1 + \hat{\lambda}_1(\hat{\mu}_1+2\hat{\mu}_2))) \bigg]
\end{align*}
where $\xi=\pi/2\kappa$, $\kappa=k+3$, $\lambda=(\lambda_1,\lambda_2)$, $\mu=(\mu_1,\mu_2)$, and $\hat{\lambda}_i=\lambda_i+1$, $\hat{\mu}_i=\mu_i+1$ for $i=1,2$.
Then for $\mu$ the distinguished vertex $\ast = (0,0)$, we obtain
\begin{align}
\psi^{\lambda}_{(0,0)} & = \frac{1}{\kappa} \bigg[ \cos(2\xi(2\hat{\lambda}_1+3\hat{\lambda}_2)) + \cos(2\xi(\hat{\lambda}_1-\hat{\lambda}_2)) - \cos(2\xi(\hat{\lambda}_1+3\hat{\lambda}_2)) - \cos(2\xi(2\hat{\lambda}_1+\hat{\lambda}_2)) \bigg] \qquad \label{eqn:PF-evector1-B2} \\
& = - \frac{1}{8 \kappa \pi^2} J_{x,y} \left( \hat{\lambda}_2/2\kappa, (\hat{\lambda}_1+\hat{\lambda}_2)/2\kappa \right), \label{eqn:PF-evector=J-C2} 
\end{align}
where in (\ref{eqn:PF-evector=J-C2}) we have $J_{x,y}(\theta_1,\theta_2)$ with $(\theta_1,\theta_2)$ related to $\lambda \in \mathrm{Exp}(\mathcal{A}_k(Sp(2)))$ by (\ref{eqn:theta-lambda_C2}).

Since the $S$-matrix is unitary, the eigenvector $\psi^{\ast}$ defined by (\ref{eqn:PF-evector1-B2}) has norm 1. Recall that the Perron-Frobenius eigenvector for $\mathcal{A}_k(Sp(2))$ can also be written in the Kac-Weyl factorized form \cite{coquereaux/rais/tahri:2010}:
\begin{equation} \label{eqn:PF-evector2-B2}
\phi^{(0,0)}_{\lambda} = \frac{ \sin(\hat{\lambda}_1\xi) \sin(2\hat{\lambda}_2\xi) \sin((\hat{\lambda}_1+2\hat{\lambda}_2)\xi) \sin((2\hat{\lambda}_1+2\hat{\lambda}_2)\xi) }{ \sin(\xi) \sin(2\xi) \sin(3\xi) \sin(4\xi) }.
\end{equation}
Now $\phi^{\ast}_{\ast} = 1$ whilst $\psi_{\ast}^{\ast} = 16\sin(\xi) \sin(2\xi) \sin(3\xi) \sin(4\xi)/\kappa$, and thus we have $\kappa \psi_{\ast} = 16\sin(\xi) \sin(2\xi) \sin(3\xi) \sin(4\xi)\phi^{\ast}$. Then from (\ref{eqn:PF-evector=J-C2}) we have
\begin{align*}
J_{x,y}(\theta_1,\theta_2) &= - 8 \kappa \pi^2 \; \psi_{\ast}^{(\hat{\lambda}_2/2\kappa,(\hat{\lambda}_1 + \hat{\lambda}_2)/2\kappa)} \\
& = - 128 \kappa \pi^2 \sin(\xi) \sin(2\xi) \sin(3\xi) \sin(4\xi) \; \phi^{\ast}_{(\hat{\lambda}_2/2\kappa,(\hat{\lambda}_1 + \hat{\lambda}_2)/2\kappa)} \\
& =  128 \kappa \pi^2 \sin(2\pi\theta_1) \sin(2\pi\theta_2) \sin(\pi(\theta_1+\theta_2)) \sin(\pi(\theta_1-\theta_2)),
\end{align*}
so that the Jacobian $J_{x,y}(\theta_1,\theta_2)$ can also be written as a product of sine functions.
A similar argument show that
\begin{align*}
J_{y,z}(\theta_1,\theta_2) &= - 4 \kappa \pi^2 \; \psi_{\ast}^{(\hat{\lambda}_1+2\hat{\lambda}_2)/4\kappa, -\hat{\lambda}_1/4\kappa} \\
& = - 64 \kappa \pi^2 \sin(\xi) \sin(2\xi) \sin(3\xi) \sin(4\xi) \; \phi^{\ast}_{(\hat{\lambda}_1+2\hat{\lambda}_2)/4\kappa, -\hat{\lambda}_1/4\kappa} \\
& = 64 \kappa \pi^2 \sin(2\pi\theta_1) \sin(2\pi\theta_2) \sin(2\pi(\theta_1+\theta_2)) \sin(2\pi(\theta_1-\theta_2)),
\end{align*}
so that the Jacobian $J_{y,z}(\theta_1,\theta_2)$ can also be written as a product of sine functions.

We now compute the joint spectral measure for $\mathcal{A}^{\rho_x}_k(Sp(2))$, $\mathcal{A}^{\rho_y}_k(Sp(2))$.
Summing over all $(\lambda_1,\lambda_2) \in \mathrm{Exp}(\mathcal{A}_k(Sp(2)))$ corresponds to summing over all $(\theta_1, \theta_2) \in \{ (\hat{\lambda}_2/2\kappa,(\hat{\lambda}_1+\hat{\lambda}_2)/2\kappa) | \; \hat{\lambda}_1,\hat{\lambda}_2 \geq 1, \hat{\lambda}_1+\hat{\lambda}_2 \leq \kappa-1 \}$, or equivalently, over all $(\theta_1, \theta_2) \in M_k^{Sp(2)} = \{ (q_1/2\kappa, q_2/2\kappa) | \; q_1, q_2 = 0,1,\ldots,2\kappa-1 \}$ such that
\begin{align}
\theta_1 & = \hat{\lambda}_2/2\kappa \geq 1/2\kappa, \qquad
\theta_1 - \theta_2 = -\hat{\lambda}_1/2\kappa \leq -1/2\kappa \label{eqn:C2theta_restrictions-1} \\
\theta_2 & = (\hat{\lambda}_1 + \hat{\lambda}_2)/2\kappa \leq (\kappa-1)/2\kappa = 1/2 - 1/2\kappa. \label{eqn:C2theta_restrictions-2}
\end{align}

Denote by $\mathcal{C}_k^{Sp(2)}$ the set of all $(\omega_1,\omega_2) \in \mathbb{T}^2$ such that $(\theta_1,\theta_2) \in M_k^{Sp(2)}$ satisfies these conditions. Then from (\ref{eqn:moments-nimrep-B2}) and (\ref{eqn:PF-evector=J-C2}) we obtain
\begin{align}
\int_{\mathbb{T}^2} (\chi_{\rho_x}(\omega_1,\omega_2))^m & (\chi_{\rho_y}(\omega_1,\omega_2))^n \mathrm{d}\varepsilon_{x,y}(\omega_1,\omega_2) \nonumber \\
& = \frac{1}{64 \kappa^2 \pi^4} \sum_{\lambda \in \mathrm{Exp}(\mathcal{A}_k(Sp(2)))} (\beta^{\lambda,Sp(2)}_{\rho_x})^m (\beta^{\lambda,Sp(2)}_{\rho_y})^n J_{x,y} \left( \hat{\lambda}_2/2\kappa, (\hat{\lambda}_1+\hat{\lambda}_2)/2\kappa \right)^2 \nonumber \\
& = \frac{1}{64 \kappa^2 \pi^4} \sum_{(\omega_1,\omega_2) \in \mathcal{C}_k^{Sp(2)}} (\chi_{\rho_x}(\omega_1,\omega_2))^m (\chi_{\rho_y}(\omega_1,\omega_2))^n J_{x,y}(\omega_1,\omega_2)^2 \label{eqn:sum_for_measure-A(C2)}
\end{align}

\begin{figure}[tb]
\begin{center}
  \includegraphics[width=55mm]{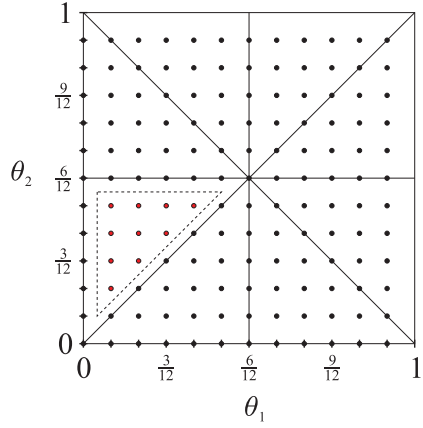}\\
 \caption{The points $(\theta_1,\theta_2)$ such that $(e^{2 \pi i \theta_1}, e^{2 \pi i \theta_2}) \in \mathcal{C}_3^W$.} \label{fig:C3W-C2}
\end{center}
\end{figure}

If we let $C^{Sp(2)}$ be the limit of $\mathcal{C}_k^{Sp(2)}$ as $k \rightarrow \infty$, then $C^{Sp(2)}$ is identified with the fundamental domain $C$ of $\mathbb{T}^2$ under the action of the group $D_8$, illustrated in Figure \ref{fig:fund_domain-B2inT2}. Since $J_{x,y}=0$ along the boundary of $C$, which is mapped to the boundary of $\mathfrak{D}_{x,y}$ under the map $\Psi_{x,y}:\mathbb{T}^2 \rightarrow \mathfrak{D}_{x,y}$, we can include points on the boundary of $C$ in the summation in (\ref{eqn:sum_for_measure-A(C2)}). Since $J_{x,y}^2$ is invariant under the action of $D_8$, we have
\begin{align}
\int_{\mathbb{T}^2} (\chi_{\rho_x}(\omega_1,\omega_2))^m & (\chi_{\rho_y}(\omega_1,\omega_2))^n \mathrm{d}\varepsilon_{x,y}(\omega_1,\omega_2) \nonumber \\
& = \frac{1}{8} \frac{1}{64 \kappa^2 \pi^4} \sum_{(\omega_1,\omega_2) \in \mathcal{C}_k^{W,Sp(2)}} (\chi_{\rho_x}(\omega_1,\omega_2))^m (\chi_{\rho_y}(\omega_1,\omega_2))^n J_{x,y}(\omega_1,\omega_2)^2 \label{eqn:sum_for_measure-A(C2)2}
\end{align}
where
\begin{equation} \label{def:Dl}
\mathcal{C}_k^{W,Sp(2)} = \{ (e^{2 \pi i q_1/2\kappa}, e^{2 \pi i q_2/2\kappa}) \in \mathbb{T}^2 | \; q_1,q_2 = 0, 1, \ldots, 2\kappa-1 \},
\end{equation}
whose intersection with the complement of the image of the boundary of the fundamental domain $C$ is the image of $\mathcal{C}_k^{Sp(2)}$ under the action of the Weyl group $W=D_8$. We illustrate the points $(\theta_1,\theta_2)$ such that $(e^{2 \pi i \theta_1}, e^{2 \pi i \theta_2}) \in \mathcal{C}_3^{W,Sp(2)}$ in Figure \ref{fig:C3W-C2}. The points in the interior of the fundamental domain $C$, those enclosed by the dashed line, correspond to the vertices of the graph $\mathcal{A}_3(Sp(2))$.

Clearly $|\mathcal{C}_k^{W,Sp(2)}| = 4(k+3)^2 = 4\kappa^2$.
Thus from (\ref{eqn:sum_for_measure-A(C2)2}), we obtain (c.f. \cite{evans/pugh:2009v}):

\begin{Thm} \label{thm:C2_k-measureT2}
The joint spectral measure of $\mathcal{A}^{\rho_x}_k(Sp(2))$, $\mathcal{A}^{\rho_y}_k(Sp(2))$, (over $\mathbb{T}^2$) is given by
\begin{equation} \label{eqn:C2_k-measureT2}
\mathrm{d}\varepsilon_{x,y}(\omega_1,\omega_2) = \frac{1}{128 \pi^4} J_{x,y}(\omega_1,\omega_2)^2 \, \mathrm{d}_{2(k+3)} \, \omega_1 \, \mathrm{d}_{2(k+3)} \, \omega_2,
\end{equation}
where $\mathrm{d}_m$ is the uniform Dirac measure over the $m^{\mathrm{th}}$ roots of unity.
\end{Thm}

In fact, $\varepsilon_{Sp(2)}: = \varepsilon_{x,y}$ is the joint spectral measure over $\mathbb{T}^2$ for any $\mathcal{A}^{\lambda}_k(Sp(2))$, $\mathcal{A}^{\mu}_k(Sp(2))$.

We now compute the joint spectral measure for $\mathcal{A}^{\rho_y}_k(SO(5))$, $\mathcal{A}^{\rho_z}_k(SO(5))$.
Summing over all $(\lambda_1,\lambda_2) \in \mathrm{Exp}(\mathcal{A}_k(SO(5)))$ corresponds to summing over all $(\theta_1, \theta_2) \in \{ (\hat{\lambda}_2/2\kappa,(\hat{\lambda}_1+\hat{\lambda}_2)/2\kappa) | \; \lambda_1,\lambda_2 \geq 0, \lambda_1+2\lambda_2 \leq k \}$, or equivalently, over all $(\theta_1, \theta_2) \in M_k^{Sp(2)}$ such that $\theta_1$, $\theta_2$ satisfy (\ref{eqn:C2theta_restrictions-1}) and $\theta_1+\theta_2 \leq 1/2$.

Denote by $\mathcal{C}_k^{SO(5)}$ the set of all $(\omega_1,\omega_2) \in \mathbb{T}^2$ such that $(\theta_1,\theta_2) \in M_k^{Sp(2)}$ satisfies these conditions. Then from (\ref{eqn:moments-nimrep-B2}) and (\ref{eqn:PF-evector=J-C2}), and since the eigenvectors $\psi^{\lambda,SO(5)}$ for $SO(5)$ satisfy $\psi^{\lambda,SO(5)}_{\nu} = \sqrt{2}\psi^{\lambda,Sp(2)}_{\nu}$ for $\nu \in \mathrm{Exp}(\mathcal{A}_k(SO(5)))$, we obtain
\begin{align}
\int_{\mathbb{T}^2} (\chi_{\rho_y}(\omega_1,\omega_2))^m & (\chi_{\rho_z}(\omega_1,\omega_2))^n \mathrm{d}\varepsilon_{y,z}(\omega_1,\omega_2) \nonumber \\
& = \frac{2}{64 \kappa^2 \pi^4} \sum_{\lambda \in \mathrm{Exp}(\mathcal{A}_k(SO(5)))} (\beta^{\lambda,SO(5)}_{\rho_y})^m (\beta^{\lambda,SO(5)}_{\rho_z})^n J_{x,y} \left( \hat{\lambda}_2/2\kappa,(\hat{\lambda}_1+\hat{\lambda}_2)/2\kappa \right)^2 \nonumber \\
& = \frac{1}{32 \kappa^2 \pi^4} \sum_{(\omega_1,\omega_2) \in \mathcal{C}_k^{SO(5)}} (\chi_{\rho_y}(\omega_1,\omega_2))^m (\chi_{\rho_z}(\omega_1,\omega_2))^n J_{x,y}(\omega_1,\omega_2)^2 \label{eqn:sum_for_measure-A(B2)} \\
\end{align}

As discussed in Section \ref{sect:measureA8inftyB2-D}, $(\chi_y(\omega_1,\omega_2),\chi_z(\omega_1,\omega_2)) = (\chi_y(\omega_1',\omega_2'),\chi_z(\omega_1',\omega_2'))$ where $(\omega_1',\omega_2')$ is given by the reflection of $(\omega_1,\omega_2)$ about either of the lines $\theta_1+\theta_2=1/2$ or $\theta_1-\theta_2=1/2$.
Thus we can write
\begin{align*}
\int_{\mathbb{T}^2} &(\chi_{\rho_y}(\omega_1,\omega_2))^m (\chi_{\rho_z}(\omega_1,\omega_2))^n \mathrm{d}\varepsilon_{y,z}(\omega_1,\omega_2) \\
&\qquad = \frac{1}{2} \frac{1}{32 \kappa^2 \pi^4} \sum_{(\omega_1,\omega_2) \in \mathcal{C}_k^{Sp(2)}} (\chi_{\rho_y}(\omega_1,\omega_2))^m (\chi_{\rho_z}(\omega_1,\omega_2))^n J_{x,y}(\omega_1,\omega_2)^2 \\
& \qquad \qquad + \frac{1}{2} \frac{1}{32 \kappa^2 \pi^4} \sum_{\stackrel{(\omega_1,\omega_2) \in \mathcal{C}_k^{Sp(2)}:}{\scriptscriptstyle{\omega_1=-\omega_2 \mbox{\tiny{ or }} \omega_1=-\omega_2^{-1}}}} (\chi_{\rho_y}(\omega_1,\omega_2))^m (\chi_{\rho_z}(\omega_1,\omega_2))^n J_{x,y}(\omega_1,\omega_2)^2
\end{align*}
where the second summation is a correction term since the points along the lines $\theta_1+\theta_2=1/2$ and $\theta_1-\theta_2=1/2$, which are fixed under the $\mathbb{Z}_2$ action, are scaled by a half in the first summation, where they are only counted once.
Then, since $J_{x,y}^2$ is invariant under the action of $D_8$, we have
\begin{align*}
\int_{\mathbb{T}^2} & (\chi_{\rho_y}(\omega_1,\omega_2))^m (\chi_{\rho_z}(\omega_1,\omega_2))^n \mathrm{d}\varepsilon_{y,z}(\omega_1,\omega_2) \\
& \qquad = \frac{1}{8} \frac{1}{64 \kappa^2 \pi^4} \sum_{(\omega_1,\omega_2) \in \mathcal{C}_k^{W,Sp(2)}} (\chi_{\rho_y}(\omega_1,\omega_2))^m (\chi_{\rho_z}(\omega_1,\omega_2))^n J_{x,y}(\omega_1,\omega_2)^2 \\
& \qquad \qquad + \frac{1}{8} \frac{1}{64 \kappa^2 \pi^4} \sum_{\stackrel{(\omega_1,\omega_2) \in \mathcal{C}_k^{W,Sp(2)}:}{\scriptscriptstyle{\omega_1=-\omega_2 \mbox{\tiny{ or }} \omega_1=-\omega_2^{-1}}}} (\chi_{\rho_y}(\omega_1,\omega_2))^m (\chi_{\rho_z}(\omega_1,\omega_2))^n J_{x,y}(\omega_1,\omega_2)^2,
\end{align*}
and we obtain:

\begin{Thm} \label{thm:B2_k-measureT2}
The joint spectral measure of $\mathcal{A}^{\rho_y}_k(SO(5))$, $\mathcal{A}^{\rho_z}_k(SO(5))$, (over $\mathbb{T}^2$) is given by
\begin{align}
\mathrm{d}\varepsilon_{y,z} & = \frac{1}{128 \pi^4} J_{x,y}(\omega_1,\omega_2)^2 \, \mathrm{d}_{2\kappa} \times \mathrm{d}_{2\kappa} \nonumber \\
& \qquad + \frac{1}{64 \kappa^2 \pi^4} J_{x,y}(\omega_1,\omega_2)^2 \sum_{j=1}^{k-1} \mathrm{d}^{(j/2\kappa,(\kappa-j)/2\kappa)},
\end{align}
where $\kappa = k+3$, $\mathrm{d}^{(\theta_1,\theta_2)}$ is as in Definition \ref{def:B2measure} and $\mathrm{d}_m$ is the uniform Dirac measure over the $m^{\mathrm{th}}$ roots of unity.
\end{Thm}

In fact, $\varepsilon_{SO(5)}: = \varepsilon_{y,z}$ is the joint spectral measure over $\mathbb{T}^2$ for any $\mathcal{A}^{\lambda}_k(SO(5))$, $\mathcal{A}^{\mu}_k(SO(5))$.

We can now easily deduce the joint spectral measures (over $\mathbb{T}^2$) for $\mathcal{A}_{\infty}(Sp(2))$, $\mathcal{A}_{\infty}(SO(5))$ claimed in Section \ref{sect:measureAinftyB2-T2}. Letting $k \rightarrow \infty$ in Theorems \ref{thm:C2_k-measureT2}, \ref{thm:B2_k-measureT2}, we obtain:

\begin{Thm} \label{thm:C2_infty-measureT2}
The joint spectral measure of any pair of infinite $Sp(2)$ graphs $\mathcal{A}^{\lambda}_{\infty}(Sp(2))$, $\mathcal{A}^{\mu}_{\infty}(Sp(2))$ (over $\mathbb{T}^2$), and the joint spectral measure of any pair of infinite $SO(5)$ graphs $\mathcal{A}^{\lambda}_{\infty}(SO(5))$, $\mathcal{A}^{\mu}_{\infty}(SO(5))$ (over $\mathbb{T}^2$), are identical and are both given by
\begin{equation} \label{eqn:jsm-A_inftyC2}
\mathrm{d}\varepsilon(\omega_1,\omega_2) = \frac{1}{128 \pi^4} J_{x,y}(\omega_1,\omega_2)^2 \mathrm{d}\omega_1 \, \mathrm{d}\omega_2,
\end{equation}
where $\mathrm{d}\omega$ is the uniform Lebesgue measure over $\mathbb{T}$.
\end{Thm}


\renewcommand{\arraystretch}{1.5}

\subsection{Graphs $\mathcal{D}_k(Sp(2))$, $k \leq \infty$} \label{sect:measures_DkB2}

\begin{figure}[tb]
\begin{center}
  \includegraphics[width=35mm]{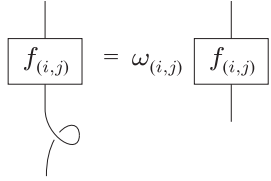}
 \caption{Statistical phase $\omega_{(i,j)}$} \label{Fig-stat_phase_wij}
\end{center}
\end{figure}

The centre of $Sp(2)$ is $\mathbb{Z}_2$. The graphs $\mathcal{D}^{\rho_u}_k(Sp(2))$, $u=x,y$, are associated to the orbifold inclusion $N \rightarrow N \rtimes_{\tau} \mathbb{Z}_2$, where $\tau = \lambda_{(0,k)}$ is a non-trivial simple current of order 2.
For such an orbifold inclusion to exist, one needs an automorphism $\tau_0$ such that $[\tau_0]=[\tau]$ and $\tau_0^2=\mbox{id}$ \cite[$\S3$]{bockenhauer/evans:1999i}, which exists precisely when the statistics phase $\omega_{\tau}$ of $\tau$ satisfies $\omega_{\tau}^2 = 1$ \cite[Lemma 4.4]{rehren:1990}.
Kuperberg's $Sp(2)$ spider \cite{kuperberg:1996} involves two types of strands, $Sp(2)$ and $SO(5)$. Using this, one can construct a semisimple braided modular tensor category whose simple objects are generalised Jones-Wenzl projections $f_{(i,j)}$, $(i,j) \in P_+^{k,Sp(2)}$ (see \cite{wenzl:1987} for ($SU(2)$) Jones-Wenzl projections and \cite{suciu:1997, kuperberg:1996, ohtsuki/yamada:1997} for generalised $SU(3)$ Jones-Wenzl projections) and whose morphisms are intertwiners between these projections (see \cite{turaev:2010, yamagami:2003, cooper:2007, evans/pugh:2010ii} for a similar construction in the case of $SU(2)$ and \cite{cooper:2007, evans/pugh:2010ii} for $SU(3)$). The statistics phase $\omega_{(i,j)} := \omega_{\lambda_{(i,j)}}$ is obtained by evaluating the twist applied to the generalised Jones-Wenzl projection $f_{(i,j)}$ (see Figure \ref{Fig-stat_phase_wij}, where the single strand drawn here represents $i$ $Sp(2)$-strands and $j$ $SO(5)$-strands). Then we see that $\omega_{(0,k)} = (-1)^k$, thus the orbifold inclusion exists. Further details will be given in a future publication. See \cite[Chapter XII]{turaev:2010} for a similar discussion in the case of $SO(3)$ and its double cover $SU(2)$.

\begin{figure}[tb]
\begin{center}
  \includegraphics[width=150mm]{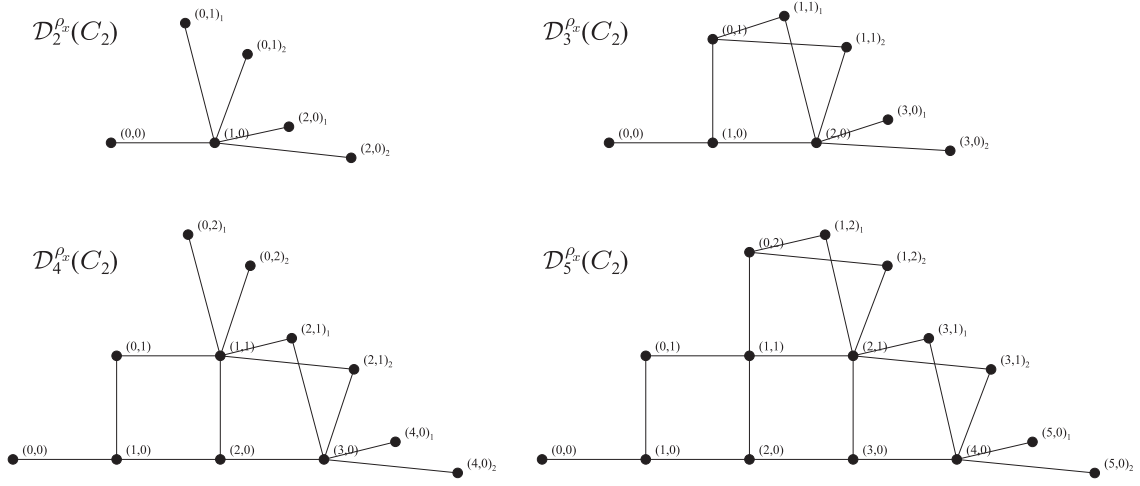}
 \caption{Orbifold graph $\mathcal{D}^{\rho_x}_k(Sp(2))$ for $k=2,3,4,5$} \label{fig-D_infty(C2)1}
\end{center}
\end{figure}

\begin{figure}[tb]
\begin{center}
  \includegraphics[width=150mm]{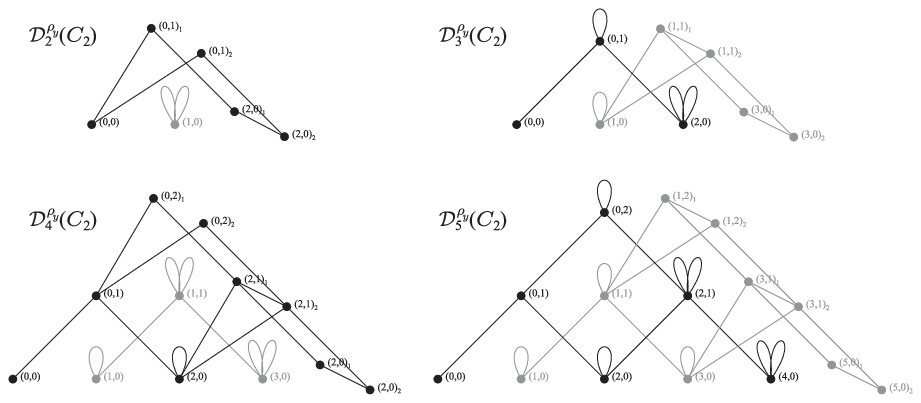}
 \caption{Orbifold graph $\mathcal{D}^{\rho_y}_k(Sp(2))$ for $k=2,3,4,5$} \label{fig-D_infty(C2)2}
\end{center}
\end{figure}

Following a similar method to \cite[$\S5.2$]{bockenhauer/evans/kawahigashi:2000}, one finds with $[\theta] = [\lambda_{(0,0)}] \oplus [\lambda_{(0,k)}]$ that $\mathcal{D}^{\rho_u}_k(Sp(2))$, $u=x,y$, are the nimrep graphs associated to the orbifold modular invariant
\begin{align*}
Z_{\mathcal{D}_{2l}} & = \sum_{\stackrel{(m,n) \in P_+^{2l,Sp(2)}(0)}{\scriptscriptstyle{m+2n<2l}}} |\chi_{(m,n)} + \chi_{(2l-m-n,2l-m-n)}|^2 + 2\sum_{j=0}^l |\chi_{(2l-2j,j)}|^2, \\
Z_{\mathcal{D}_{2l+1}} & = \sum_{(m,n) \in P_+^{2l+1,Sp(2)}(0)} |\chi_{(m,n)}|^2 + \sum_{\stackrel{(m,n) \in P_+^{2l+1,Sp(2)}}{\scriptscriptstyle{m \mbox{\tiny{ odd}}}}}\chi_{(m,n)}\chi_{(m,2l+1-m-n)}^{\ast},
\end{align*}
where $P_+^{k,Sp(2)}(0) = \{(m,n) \in P_+^{k,Sp(2)} | \, m \mbox{ even} \}$. This modular invariant appeared in \cite{bernard:1987}. The graphs $\mathcal{D}^{\rho_u}_k(Sp(2))$ are are $\mathbb{Z}_2$-orbifolds of the graphs $\mathcal{D}^{\rho_u}_k(Sp(2))$, and are illustrated in Figures \ref{fig-D_infty(C2)1}, \ref{fig-D_infty(C2)2}, where we have labeled the vertices by the corresponding Dynkin labels from the $\mathcal{A}_k(Sp(2))$ graphs.

The exponents of $\mathcal{D}_k(Sp(2))$ are given by $\mathrm{Exp}(\mathcal{D}_{2l}(Sp(2))) = \{ (m,n) \in P_+^{2l,Sp(2)}(0) | \, m \neq 2l-2n \} \cup \{ \mbox{twice } (2l-2j,j) | \, j=0,1,\ldots, l \}$ for $k=2l$ even, whilst $\mathrm{Exp}(\mathcal{D}_{2l+1}(Sp(2))) = P_+^{2l+1,Sp(2)}(0) \cup \{ (2l+1-2j,j) | \, j=0,1,\ldots, l \}$ for $k=2l+1$ odd. For $\lambda \in P_+^{k,Sp(2)}$ (which label the vertices of $\mathcal{A}_k(Sp(2))$) not a fixed point under the $\mathbb{Z}_2$-action, i.e. $\lambda \not \in \{ (k-2j,j) | \, j=0,1,\ldots, \lfloor k/2 \rfloor \}$ where $\lfloor x \rfloor$ denotes the integer part of $x$, the normalized eigenvector satisfies $|\psi_{\ast}^{\lambda}|^2 = 2 S_{\ast,\lambda}^2$. However for $\lambda \in \{ (k-2j,j) | \, j=0,1,\ldots, \lfloor k/2 \rfloor \}$, $|\psi_{\ast}^{\lambda_1}| = |\psi_{\ast}^{\lambda_2}| = \sqrt{2} S_{\ast,\lambda}/2$, where $\lambda_j$, $j=1,2$, denote the two copies of the fixed point in the orbifold graph $\mathcal{D}_k(Sp(2))$, so that $|\psi_{\ast}^{\lambda_1}|^2 + |\psi_{\ast}^{\lambda_2}|^2 = S_{\ast,\lambda}^2$.

With $\theta_1, \theta_2$ as in (\ref{eqn:theta-lambda_C2}), summing over all $\lambda = (\lambda_1,\lambda_2) \in P_+^{k,Sp(2)}(0)$ corresponds to summing over all $(\omega_1,\omega_2) \in \mathcal{C}_k^{Sp(2)}$ such that $\omega_1\omega_2 = e^{2\pi i (2m+1)/2\kappa}$ for $m \in \mathbb{Z}$, where $\mathcal{C}_k^{Sp(2)}$ is as in Section \ref{sect:measures_AkB2}.
Then from (\ref{eqn:moments-nimrep-B2}) and (\ref{eqn:PF-evector=J-C2}), with $\zeta_k = 1$ for $k$ odd and $\zeta_k = 3/2$ for $k$ even, we obtain
\begin{align*}
\int_{\mathbb{T}^2} & (\chi_{\rho_x}(\omega_1,\omega_2))^m (\chi_{\rho_y}(\omega_1,\omega_2))^n \mathrm{d}\varepsilon(\omega_1,\omega_2) \\
& \quad = \frac{2}{64 \kappa^2 \pi^4} \sum_{\stackrel{(\omega_1,\omega_2) \in \mathcal{C}_k^{Sp(2)}:}{\scriptscriptstyle{\omega_1\omega_2 = e^{2\pi i (2m+1)/2\kappa}}}} (\chi_{\rho_x}(\omega_1,\omega_2))^m (\chi_{\rho_y}(\omega_1,\omega_2))^n J_{x,y}(\omega_1,\omega_2)^2 \\
& \qquad \qquad + \frac{\zeta_k}{64 \kappa^2 \pi^4} \sum_{\stackrel{(\omega_1,\omega_2) \in \mathcal{C}_k^{Sp(2)}:}{\scriptscriptstyle{\omega_1=-\omega_2 \mbox{\tiny{ or }} \omega_1=-\omega_2^{-1}}}} (\chi_{\rho_x}(\omega_1,\omega_2))^m (\chi_{\rho_y}(\omega_1,\omega_2))^n J_{x,y}(\omega_1,\omega_2)^2 \\
& \quad = \frac{1}{8} \frac{1}{32 \kappa^2 \pi^4} \sum_{\stackrel{(\omega_1,\omega_2) \in \mathcal{C}_k^{W,Sp(2)}:}{\scriptscriptstyle{\omega_1\omega_2 = e^{2\pi i (2m+1)/2\kappa}}}} (\chi_{\rho_x}(\omega_1,\omega_2))^m (\chi_{\rho_y}(\omega_1,\omega_2))^n J_{x,y}(\omega_1,\omega_2)^2 \\
& \qquad \qquad + \frac{1}{8} \frac{\zeta_k}{64 \kappa^2 \pi^4} \sum_{\stackrel{(\omega_1,\omega_2) \in \mathcal{C}_k^{W,Sp(2)}:}{\scriptscriptstyle{\omega_1=-\omega_2 \mbox{\tiny{ or }} \omega_1=-\omega_2^{-1}}}} (\chi_{\rho_x}(\omega_1,\omega_2))^m (\chi_{\rho_y}(\omega_1,\omega_2))^n J_{x,y}(\omega_1,\omega_2)^2.
\end{align*}

Thus
\begin{Thm} \label{thm:D-measure}
The joint spectral measure of $\mathcal{D}^{\rho_x}_k(Sp(2))$, $\mathcal{A}^{\rho_y}_k(Sp(2))$, (over $\mathbb{T}^2$) is given by
\begin{align}
\mathrm{d}\varepsilon & = \frac{1}{128 \pi^4} J_{x,y}(\omega_1,\omega_2)^2 \, \left( \mathrm{d}_{\kappa} \times \left( \mathrm{d}_{2\kappa} - \mathrm{d}_{\kappa} \right) + \left( \mathrm{d}_{2\kappa} - \mathrm{d}_{\kappa} \right) \times \mathrm{d}_{\kappa} \right) \nonumber \\
& \qquad + \frac{\zeta_k}{64 \kappa^2 \pi^4} J_{x,y}(\omega_1,\omega_2)^2 \sum_{j=1}^{k-1} \mathrm{d}^{(j/2\kappa,(\kappa-j)/2\kappa)},
\end{align}
where $\kappa = k+3$, $\zeta_k = 1$ for $k$ odd and $\zeta_k = 3/2$ for $k$ even, $\mathrm{d}^{(\theta_1,\theta_2)}$ is as in Definition \ref{def:B2measure} and $\mathrm{d}_m$ is the uniform Dirac measure over the $m^{\mathrm{th}}$ roots of unity.
\end{Thm}

Letting $k \rightarrow \infty$ we easily obtain the following corollary:

\begin{Cor} \label{cor:D_infty-measure}
The joint spectral measure of $\mathcal{D}^{\rho_x}_{\infty}(Sp(2))$, $\mathcal{D}^{\rho_y}_{\infty}(Sp(2))$, (over $\mathbb{T}^2$) is precisely the joint spectral measure of the infinite $Sp(2)$ graphs $\mathcal{A}^{\rho_x}_{\infty}(Sp(2))$, $\mathcal{A}^{\rho_y}_{\infty}(Sp(2))$, given in Theorem \ref{thm:C2_infty-measureT2}.
\end{Cor}

\subsection{Exceptional Graph $\mathcal{E}_3(Sp(2))$: $(Sp(2))_3 \rightarrow (SO(10))_1$} \label{sect:measures_E3C2}

\begin{figure}[tb]
\begin{minipage}[t]{7.5cm}
\begin{center}
  \includegraphics[width=25mm]{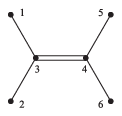} \\
 \caption{Exceptional Graph $\mathcal{E}_3^{\rho_1}(Sp(2))$} \label{fig-Graph_E3_B2-1}
\end{center}
\end{minipage}
\hfill
\begin{minipage}[t]{7.5cm}
\begin{center}
  \includegraphics[width=60mm]{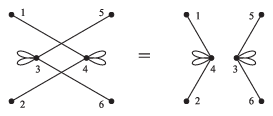} \\
 \caption{Exceptional Graph $\mathcal{E}_3^{\rho_2}(Sp(2))$} \label{fig-Graph_E3_B2-2}
\end{center}
\end{minipage}
\end{figure}

The graphs $\mathcal{E}_3(Sp(2))$ are associated to the conformal embedding $(Sp(2))_3 \rightarrow (SO(10))_1$ and are one of two nimreps associated to the modular invariant
$$Z_{\mathcal{E}_3} = |\chi_{(0,0)}+\chi_{(2,1)}|^2 + |\chi_{(2,0)}+\chi_{(0,3)}|^2 + 2|\chi_{(1,1)}|^2$$
which is at level 3 and has exponents $\mathrm{Exp}(\mathcal{E}_3(Sp(2))) = \{ (0,0), (2,1), (2,0), (0,3), \textrm{ and } (1,1) \textrm{ twice } \}$.
The other family $\mathcal{E}_3^M(Sp(2))$ are considered in the next section. The graphs $\mathcal{E}_3^{\rho_j}(Sp(2))$ are illustrated in Figures \ref{fig-Graph_E3_B2-1}, \ref{fig-Graph_E3_B2-2}. Note that $\mathcal{E}_3^{\rho_2}(Sp(2))$ has two connected components.

\begin{figure}[tb]
\begin{center}
  \includegraphics[width=75mm]{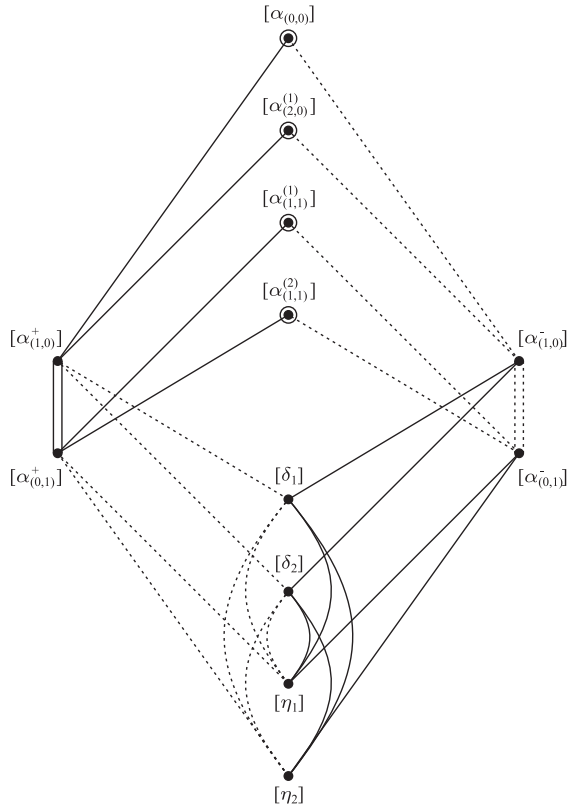} \\
 \caption{$\mathcal{E}_3(Sp(2))$: Multiplication by $[\alpha_{(1,0)}^+]$ (solid lines) and $[\alpha_{(1,0)}^-]$ (dashed lines)} \label{Fig-full_E3_C2}
\end{center}
\end{figure}

Following \cite[$\S$6]{bockenhauer/evans:1999ii} we can compute the principal graph and dual principal graph of the inclusion $(Sp(2))_3 \rightarrow (SO(10))_1$. The chiral induced sector bases ${}_M \mathcal{X}_M^{\pm} \subset \mbox{Sect}(M)$ and full induced sector basis ${}_M \mathcal{X}_M \subset \mbox{Sect}(M)$, the sector bases given by all irreducible subsectors of $[\alpha_{\lambda}^{\pm}]$ and $[\alpha_{\lambda}^+ \circ \alpha_{\lambda'}^-]$ respectively, for $\lambda, \lambda' \in {}_N \mathcal{X}_N$, are given by
\begin{align*}
{}_M \mathcal{X}_M^{\pm} &= \{ [\alpha_{(0,0)}], [\alpha_{(1,0)}^{\pm}], [\alpha_{(0,1)}^{\pm}], [\alpha_{(2,0)}^{(1)}], [\alpha_{(1,1)}^{(1)}], [\alpha_{(1,1)}^{(2)}] \},\\
{}_M \mathcal{X}_M &= \{ [\alpha_{(0,0)}], [\alpha_{(1,0)}^+], [\alpha_{(1,0)}^-], [\alpha_{(0,1)}^+], [\alpha_{(0,1)}^-], [\alpha_{(2,0)}^{(1)}], [\alpha_{(1,1)}^{(1)}], [\alpha_{(1,1)}^{(2)}], [\delta_1], [\delta_2], [\eta_1], [\eta_2] \},
\end{align*}
where $[\alpha_{(2,0)}^{\pm}] = [\alpha_{(0,1)}^{\pm}] \oplus [\alpha_{(2,0)}^{(1)}]$, $[\alpha_{(1,1)}^{\pm}] = [\alpha_{(1,0)}^{\pm}] \oplus [\alpha_{(1,1)}^{(1)}] \oplus [\alpha_{(1,1)}^{(2)}]$, $[\alpha_{(1,0)}^+ \alpha_{(1,0)}^-] = [\delta_1] \oplus [\delta_2]$, and $[\alpha_{(1,0)}^+ \alpha_{(0,1)}^-] = [\eta_1] \oplus [\eta_2]$, for $\alpha_{(i,j)} \equiv \alpha_{\lambda_{(i,j)}}$.
The fusion graphs of $[\alpha_{(1,0)}^+]$ (solid lines) and $[\alpha_{(1,0)}^-]$ (dashed lines) are given in Figure \ref{Fig-full_E3_C2}, see also \cite[Figure 7(a)]{coquereaux/rais/tahri:2010}. The marked vertices corresponding to sectors in the neutral system ${}_M \mathcal{X}_M^0 = {}_M \mathcal{X}_M^+ \cap {}_M \mathcal{X}_M^-$ have been circled.
Note that multiplication by $[\alpha_{(1,0)}^+]$ (or $[\alpha_{(1,0)}^-]$) does not give two copies of the nimrep graph $\mathcal{E}_3(Sp(2))$ as one might expect, but rather one copy each of $\mathcal{E}_3^{\rho_1}(Sp(2))$ and $\mathcal{E}_3^{\rho_1,M}(Sp(2))$. This is similar to the situation for the $SU(3)$ conformal embedding $SU(3)_9 \rightarrow (E_6)_1$ \cite[$\S$5.2]{evans:2002}.

\begin{figure}[tb]
\begin{center}
  \includegraphics[width=62.5mm]{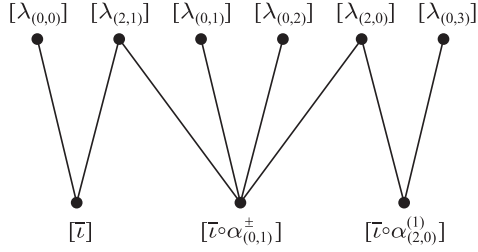} \\
 \caption{$\mathcal{E}_3(Sp(2))$: Principal graph of $(Sp(2))_3 \rightarrow (SO(10))_1$} \label{Fig-GHJ_Graph_E3_C2}
\end{center}
\end{figure}

Let $\iota:N \hookrightarrow M$ denote the injection map $\iota(n)=n \in M$, $n \in N$ and $\overline{\iota}$ its conjugate. The dual canonical endomorphism $\theta = \overline{\iota} \iota$ for the conformal embedding can be read from the vacuum block of the modular invariant: $[\theta] = [\lambda_{(0,0)}] \oplus [\lambda_{(2,1)}]$.
By \cite[Corollary 3.19]{bockenhauer/evans:1999ii} and the fact that $\langle \gamma, \gamma \rangle_M = \langle \theta, \theta \rangle_N = 2$, the canonical endomorphism $\gamma = \iota \overline{\iota}$ is given by
\begin{equation}
[\gamma] = [\alpha_{(0,0)}] \oplus [\delta_1].
\end{equation}
Then by \cite[Theorem 4.2]{bockenhauer/evans:1999ii}, the principal graph of the inclusion $(Sp(2))_3 \rightarrow (SO(10))_1$ of index $3+\sqrt{3} \approx 4.73$ is given by the connected component of $[\lambda_{(0,0)}] \in {}_N \mathcal{X}_N$ of the induction-restriction graph, and the dual principal graph is given by the connected component of $[\alpha_{(0,0)}] \in {}_M \mathcal{X}_M$ of the $\gamma$-multiplication graph. The principal graph and dual principal graph are the same, and we illustrate the principal graph in Figure \ref{Fig-GHJ_Graph_E3_C2}. These are the principal graphs for the 3311 Goodman-de la Harpe-Jones subfactor \cite{goodman/de_la_harpe/jones:1989}. The principal graph in Figure \ref{Fig-GHJ_Graph_E3_C2} appears as the intertwiner for the quantum subgroup $\mathcal{E}_3(Sp(2))$ in \cite[Figure 9]{coquereaux/rais/tahri:2010}.

One can also construct a subfactor $\alpha_{(1,0)}^{\pm}(M) \subset M$ with index $(1+\sqrt{3})^2 = 2(2+\sqrt{3}) \approx 7.46$, where $M$ is a type III factor.
Its principal graph is the nimrep graph $\mathcal{E}_3^{\rho_1}(Sp(2))$ illustrated in Figure \ref{fig-Graph_E3_B2-1}. The dual principal graph is isomorphic to the principal graph as abstract graphs \cite[Corollary 3.7]{xu:1998}.

We now determine the joint spectral measure of $\mathcal{E}_3^{\rho_1}(Sp(2))$, $\mathcal{E}_3^{\rho_2}(Sp(2))$.
With $\theta_1, \theta_2$ as in (\ref{eqn:theta-lambda_C2}) for $\lambda = (\lambda_1,\lambda_2) \in \mathrm{Exp}(\mathcal{E}_3(Sp(2)))$, we have the following values:
\begin{center}
\begin{tabular}{|c|c|c|c|} \hline
$\lambda \in \mathrm{Exp}$ & $(\theta_1,\theta_2) \in [0,1]^2$ & $|\psi^{\lambda}_{\ast}|^2$ & $\frac{1}{8\pi^2}|J(\theta_1,\theta_2)|$ \\
\hline $(0,0)$ & $\left(\frac{1}{12},\frac{2}{12}\right)$ & $\frac{3-\sqrt{3}}{24}$ & $\frac{3-\sqrt{3}}{2}$ \\
\hline $(2,1)$ & $\left(\frac{2}{12},\frac{5}{12}\right)$ & $\frac{3+\sqrt{3}}{24}$ & $\frac{3+\sqrt{3}}{2}$ \\
\hline $(2,0)$ & $\left(\frac{1}{12},\frac{4}{12}\right)$ & $\frac{3+\sqrt{3}}{24}$ & $\frac{3+\sqrt{3}}{2}$ \\
\hline $(0,3)$ & $\left(\frac{4}{12},\frac{5}{12}\right)$ & $\frac{3-\sqrt{3}}{24}$ & $\frac{3-\sqrt{3}}{2}$ \\
\hline $(1,1)$ & $\left(\frac{2}{12},\frac{4}{12}\right)$ & $\frac{1}{2}$ & $3$ \\
\hline
\end{tabular}
\end{center}
where the eigenvectors $\psi^{\lambda}$ have been normalized so that $||\psi^{\lambda}|| = 1$, and for the exponent $(1,1)$ which has multiplicity two, the value listed in the table for $|\psi^{(1,1)}_{\ast}|^2$ is $|\psi^{(1,1)_1}_{\ast}|^2 + |\psi^{(1,1)_2}_{\ast}|^2$.
Note that
\begin{equation} \label{eqn:psi=zetaJ-E3}
|\psi^{\lambda}_{\ast}|^2 = \zeta_{\lambda} \frac{1}{12} \frac{1}{8\pi^2} |J|
\end{equation}
where $\zeta_{\lambda} = 1$ for $\lambda \in \{ (0,0), (2,1), (2,0), (0,3) \}$ and $\zeta_{(1,1)} = 2$.

\begin{figure}[tb]
\begin{minipage}[t]{7.5cm}
\begin{center}
  \includegraphics[width=45mm]{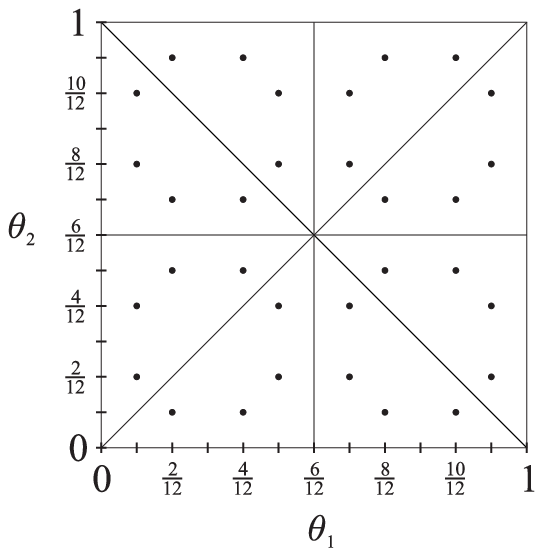}\\
 \caption{Orbit of $(\theta_1,\theta_2) \neq \left(\frac{2}{12},\frac{4}{12}\right)$.} \label{Fig-E3pointsB2}
\end{center}
\end{minipage}
\hfill
\begin{minipage}[t]{7.5cm}
\begin{center}
  \includegraphics[width=45mm]{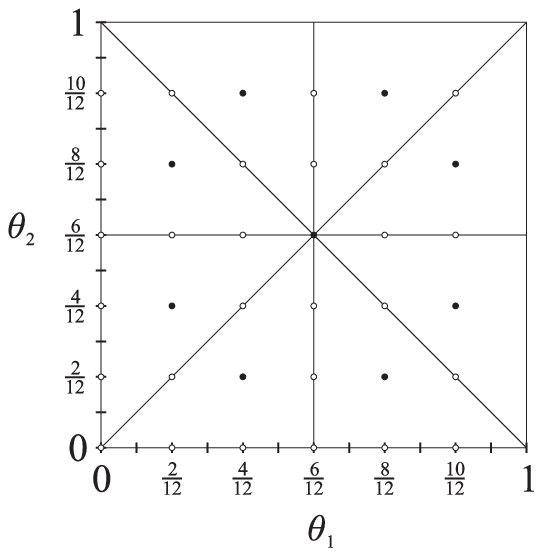}\\
 \caption{Orbit of $(\theta_1,\theta_2) = \left(\frac{2}{12},\frac{4}{12}\right)$.} \label{Fig-E3pointsB2-2}
\end{center}
\end{minipage}
\end{figure}

The orbit under $D_8$ of the points $(\theta_1,\theta_2) \in \left\{ \left(\frac{1}{12},\frac{2}{12}\right), \left(\frac{2}{12},\frac{5}{12}\right), \left(\frac{1}{12},\frac{4}{12}\right), \left(\frac{4}{12},\frac{5}{12}\right) \right\},$ are illustrated in Figure \ref{Fig-E3pointsB2}, whilst the orbit of $\left(\frac{2}{12},\frac{4}{12}\right)$ is illustrated by the black points in Figure \ref{Fig-E3pointsB2-2}.
The orbits of the first four points support the measure $\mathrm{d}^{(1/12,2/12)}+\mathrm{d}^{(1/12,4/12)}$, where $\mathrm{d}^{(\theta_1,\theta_2)}$ is the discrete uniform measure given in Definition \ref{def:B2measure}. Since the hollow points in Figure \ref{Fig-E3pointsB2-2} lie on the boundary of the orbit of fundamental domain, $J=0$ at these points, thus we see that the orbit of $(2/12,4/12)$ supports the measure $|J| \, \mathrm{d}_6 \times \mathrm{d}_6$, where $\mathrm{d}_n$ is the uniform Dirac measure on the $n^{\mathrm{th}}$ roots of unity. Note that when taking the orbit under $D_8$, the associated weight in (\ref{eqn:psi=zetaJ-E3}) is now counted 8 times, thus we must divide (\ref{eqn:psi=zetaJ-E3}) by 8.
Thus the joint spectral measure for $\mathcal{E}_3(Sp(2))$ is
$$\mathrm{d}\varepsilon = 16 \, \frac{1}{8} \, \frac{1}{12} \, \frac{1}{8\pi^2} |J| \, \left( \mathrm{d}^{(1/12,2/12)} + \mathrm{d}^{(1/12,4/12)} \right) + 36 \, \frac{1}{8} \, \frac{2}{12} \, \frac{1}{8\pi^2} |J| \, \mathrm{d}_6 \times \mathrm{d}_6.$$
Then we have obtained the following result:

\begin{Thm}
The joint spectral measure of $\mathcal{E}_3^{\rho_1}(Sp(2))$, $\mathcal{E}_3^{\rho_2}(Sp(2))$ (over $\mathbb{T}^2$) is
\begin{equation}
\mathrm{d}\varepsilon = \frac{1}{48\pi^2} |J| \, \mathrm{d}^{(1/12,2/12)} + \frac{1}{48\pi^2} |J| \, \mathrm{d}^{(1/12,4/12)} + \frac{1}{384\pi^2} |J| \, \mathrm{d}_6 \times \mathrm{d}_6,
\end{equation}
where $\mathrm{d}^{(\theta_1,\theta_2)}$ is as in Definition \ref{def:B2measure} and $\mathrm{d}_6$ is the uniform Dirac measure on the $6^{\mathrm{th}}$ roots of unity.
\end{Thm}

\subsection{Exceptional Graph $\mathcal{E}_3^M(Sp(2))$: $(Sp(2))_3 \rightarrow (SO(10))_1 \rtimes \mathbb{Z}_2$} \label{sect:measures_E3MC2}

\begin{figure}[tb]
\begin{minipage}[t]{5.5cm}
\begin{center}
  \includegraphics[width=25mm]{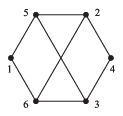} \\
 \caption{Graph $\mathcal{E}_3^{M,\rho_1}(Sp(2))$} \label{fig-Graph_E3M_B2-1}
\end{center}
\end{minipage}
\hfill
\begin{minipage}[t]{9.5cm}
\begin{center}
  \includegraphics[width=80mm]{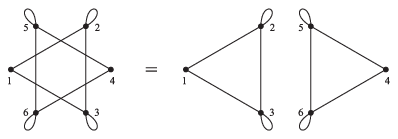} \\
 \caption{Graph $\mathcal{E}_3^{M,\rho_2}(Sp(2))$} \label{fig-Graph_E3M_B2-2}
\end{center}
\end{minipage}
\end{figure}

The graphs $\mathcal{E}_3^{M,\rho_j}(Sp(2))$, illustrated in Figures \ref{fig-Graph_E3M_B2-1}, \ref{fig-Graph_E3M_B2-2}, are the nimrep graphs for the type II inclusion $(Sp(2))_3 \rightarrow (SO(10))_1 \rtimes_{\tau} \mathbb{Z}_2$ with index $2(3+\sqrt{3}) \approx 9.46$, where $\tau = \alpha_{(2,0)}^{(1)}$ is a non-trivial simple current of order 2 in the ambichiral system ${}_M \mathcal{X}_M^0$, see Figure \ref{Fig-full_E3_C2}.
Now $\omega_{(2,0)} = -1$ \cite{coquereaux:webpage}, thus the orbifold inclusion exists (c.f. Section \ref{sect:measures_DkB2}). Note that $[\tau']=[\alpha_{(1,1)}^{(j)}] \in {}_M \mathcal{X}_M^0$ is a subsector of $[\alpha_{(1,1)}^{\pm}]$, for which $\omega_{(1,1)} = e^{7\pi i/4}$ \cite{coquereaux:webpage}. Then $\omega_{(1,1)}^2 = e^{7\pi i/2} \neq 1$, and hence the orbifold inclusion $(Sp(2))_4 \rightarrow (SO(10))_1 \rtimes_{\tau'} \mathbb{Z}_2$ does not exist.

\begin{figure}[tb]
\begin{center}
  \includegraphics[width=60mm]{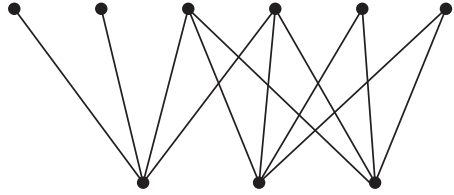} \\
 \caption{$\mathcal{E}_3^M(Sp(2))$: Principal graph of $(Sp(2))_3 \rightarrow (SO(10))_1 \rtimes \mathbb{Z}_2$} \label{Fig-GHJ_Graph_E3M_C2}
\end{center}
\end{figure}

The principal graph for this inclusion is illustrated in Figure \ref{Fig-GHJ_Graph_E3M_C2}. This will be discussed in a future publication using a generalised Goodman-de la Harpe-Jones construction analogous to that for the $D_{\mathrm{odd}}$ and $E_7$ modular invariants for $SU(2)$ \cite[$\S5.2, 5.3$]{bockenhauer/evans/kawahigashi:2000} and the type II inclusions for $SU(3)$ \cite[$\S5$]{evans/pugh:2009ii}.
It is not clear what the dual principal graph is in this case.

The associated modular invariant is again $Z_{\mathcal{E}_3}$ and the graphs are isospectral to $\mathcal{E}_3(Sp(2))$.
The eigenvectors $\psi^{\lambda}$ are not identical to those for $\mathcal{E}_3(Sp(2))$, however, as seen in the following table, the values of $|\psi^{\lambda}_{\ast}|^2$ are equal (up to a factor 2) to those for $\mathcal{E}_3(Sp(2))$, for $\lambda \neq (1,1)$.
With $\theta_1, \theta_2$ as in (\ref{eqn:theta-lambda_C2}) for $\lambda = (\lambda_1,\lambda_2) \in \mathrm{Exp}$, we have:
\begin{center}
\begin{tabular}{|c|c|c|c|} \hline
$\lambda \in \mathrm{Exp}$ & $(\theta_1,\theta_2) \in [0,1]^2$ & $|\psi^{\lambda}_{\ast}|^2$ & $\frac{1}{8\pi^2}|J(\theta_1,\theta_2)|$ \\
\hline $(0,0)$ & $\left(\frac{1}{12},\frac{2}{12}\right)$ & $\frac{3-\sqrt{3}}{12}$ & $\frac{3-\sqrt{3}}{2}$ \\
\hline $(2,1)$ & $\left(\frac{2}{12},\frac{5}{12}\right)$ & $\frac{3+\sqrt{3}}{12}$ & $\frac{3+\sqrt{3}}{2}$ \\
\hline $(2,0)$ & $\left(\frac{1}{12},\frac{4}{12}\right)$ & $\frac{3+\sqrt{3}}{12}$ & $\frac{3+\sqrt{3}}{2}$ \\
\hline $(0,3)$ & $\left(\frac{4}{12},\frac{5}{12}\right)$ & $\frac{3-\sqrt{3}}{12}$ & $\frac{3-\sqrt{3}}{2}$ \\
\hline $(1,1)$ & $\left(\frac{2}{12},\frac{4}{12}\right)$ & $\frac{1}{2}$ & $3$ \\
\hline
\end{tabular}
\end{center}
where the eigenvectors $\psi^{\lambda}$ have been normalized so that $||\psi^{\lambda}|| = 1$.
Then (\ref{eqn:psi=zetaJ-E3}) becomes $|\psi^{\lambda}_{\ast}|^2 = \zeta_{\lambda} \frac{1}{6} \frac{1}{8\pi^2} |J|$, where $\zeta_{\lambda}$ is as for $\mathcal{E}_3(Sp(2))$.
Thus we have the following result:

\begin{Thm}
The joint spectral measure of $\mathcal{E}_3^{M,\rho_1}(Sp(2))$, $\mathcal{E}_3^{M,\rho_2}(Sp(2))$ (over $\mathbb{T}^2$) is
\begin{equation}
\mathrm{d}\varepsilon = \frac{1}{24\pi^2} |J| \, \mathrm{d}^{(1/12,2/12)} + \frac{1}{24\pi^2} |J| \, \mathrm{d}^{(1/12,4/12)} + \frac{1}{384\pi^2} |J| \, \mathrm{d}_6 \times \mathrm{d}_6,
\end{equation}
where $\mathrm{d}^{(\theta_1,\theta_2)}$ is as in Definition \ref{def:B2measure} and $\mathrm{d}_6$ is the uniform Dirac measure on the $6^{\mathrm{th}}$ roots of unity.
\end{Thm}

\subsection{Exceptional Graph $\mathcal{E}_7(Sp(2))$: $(Sp(2))_7 \rightarrow (SO(14))_1$} \label{sect:measures_E7C2}

\begin{figure}[tb]
\begin{minipage}[t]{5cm}
\begin{center}
  \includegraphics[width=40mm]{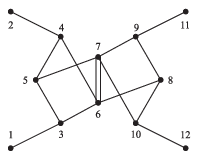} \\
 \caption{Graph $\mathcal{E}_7^{\rho_1}(Sp(2))$} \label{fig-Graph_E7_B2-1}
\end{center}
\end{minipage}
\hfill
\begin{minipage}[t]{10cm}
\begin{center}
  \includegraphics[width=90mm]{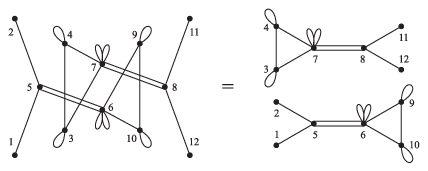} \\
 \caption{Graph $\mathcal{E}_7^{\rho_2}(Sp(2))$} \label{fig-Graph_E7_B2-2}
\end{center}
\end{minipage}
\end{figure}

The graphs $\mathcal{E}_7^{\rho_j}(G_2)$, illustrated in Figures \ref{fig-Graph_E7_B2-1}, \ref{fig-Graph_E7_B2-2}, are the nimrep graphs associated to the conformal embedding $(Sp(2))_7 \rightarrow (SO(14))_1$ and are one of two nimreps associated to the modular invariant
$$Z_{\mathcal{E}_7} = |\chi_{(0,0)}+\chi_{(6,1)}+\chi_{(2,2)}+\chi_{(0,5)}|^2 + |\chi_{(6,0)}+\chi_{(0,2)}+\chi_{(2,3)}+\chi_{(0,7)}|^2 + 2|\chi_{(3,1)}+\chi_{(3,3)}|^2$$
which is at level 7 and has 12 exponents
$$\mathrm{Exp}(\mathcal{E}_7(Sp(2))) = \{ (0,0), (6,1), (2,2), (0,5), (6,0), (0,2), (2,3), (0,7) \textrm{ and twice } (3,1), (3,3). \}$$
Note again that for the second fundamental representation $\rho_2$, the graph (Figure \ref{fig-Graph_E7_B2-2}) has two connected components.

As in Section \ref{sect:measures_E3C2}, we can compute the principal graph and dual principal graph of the inclusion $(Sp(2))_7 \rightarrow (SO(14))_1$. The chiral induced sector bases ${}_M \mathcal{X}_M^{\pm}$ and full induced sector basis ${}_M \mathcal{X}_M$ are given by
\begin{align*}
{}_M \mathcal{X}_M^{\pm} &= \{ [\alpha_{(0,0)}], [\alpha_{(1,0)}^{\pm}], [\alpha_{(0,1)}^{\pm}], [\alpha_{(2,0)}^{\pm}], [\alpha_{(1,1)}^{(j)\pm}], [\alpha_{(0,2)}^{(1)}], [\alpha_{(3,0)}^{(1)\pm}], [\alpha_{(2,1)}^{(j)\pm}], [\alpha_{(3,1)}^{(j)\pm}], \mbox{ for } j=1,2 \},\\
{}_M \mathcal{X}_M &= \{ [\alpha_{(0,0)}], [\alpha_{(1,0)}^{\varepsilon}], [\alpha_{(0,1)}^{\varepsilon}], [\alpha_{(2,0)}^{\varepsilon}], [\alpha_{(1,1)}^{(j)\varepsilon}], [\alpha_{(0,2)}^{(1)}], [\alpha_{(3,0)}^{(1)\varepsilon}], [\alpha_{(2,1)}^{(j)\varepsilon}], [\alpha_{(3,1)}^{(j)\varepsilon}], [\alpha_{(1,0)}^+\alpha_{(1,0)}^-], \\
& \qquad [\alpha_{(1,0)}^+\alpha_{(0,1)}^-], [\alpha_{(0,1)}^+\alpha_{(1,0)}^-], [\alpha_{(1,0)}^+\alpha_{(2,0)}^-], [\alpha_{(2,0)}^+\alpha_{(1,0)}^-], [\eta_j], [\zeta_j], [\psi_j], [\psi_j'], [\xi], [\xi'], [\varphi], \\
& \qquad [\gamma_j], [\gamma_j'], [\delta_j], [\lambda], [\lambda'], [\omega_j], [\phi_j], \mbox{ for } \varepsilon=+,- \mbox{ and } j=1,2 \},
\end{align*}
where $[\alpha_{(1,1)}^{\pm}] = [\alpha_{(1,1)}^{(1)\pm}] \oplus [\alpha_{(1,1)}^{(2)\pm}]$,
$[\alpha_{(0,2)}^{\pm}] = [\alpha_{(2,0)}^{\pm}] \oplus [\alpha_{(0,2)}^{(1)}]$,
$[\alpha_{(3,0)}^{\pm}] = [\alpha_{(1,1)}^{(1)\pm}] \oplus [\alpha_{(3,0)}^{(1)\pm}]$,
$[\alpha_{(2,1)}^{\pm}] = [\alpha_{(0,1)}^{\pm}] \oplus [\alpha_{(2,0)}^{\pm}] \oplus [\alpha_{(2,1)}^{(1)\pm}] \oplus [\alpha_{(2,1)}^{(2)\pm}]$,
$[\alpha_{(3,1)}^{\pm}] = [\alpha_{(1,0)}^{\pm}] \oplus 2[\alpha_{(1,1)}^{(1)\pm}] \oplus [\alpha_{(3,0)}^{(1)\pm}] \oplus [\alpha_{(3,1)}^{(1)}] \oplus [\alpha_{(3,1)}^{(2)}]$,
$[\alpha_{(0,1)}^+ \alpha_{(0,1)}^-] = [\eta_1] \oplus [\eta_2]$,
$[\alpha_{(2,0)}^+ \alpha_{(2,0)}^-] = [\zeta_1] \oplus [\zeta_2]$,
$[\alpha_{(2,0)}^+ \alpha_{(0,1)}^-] = [\psi_1] \oplus [\psi_2]$,
$[\alpha_{(0,1)}^+ \alpha_{(2,0)}^-] = [\psi_1'] \oplus [\psi_2']$,
$[\alpha_{(1,1)}^+ \alpha_{(1,0)}^-] = [\xi] \oplus [\varphi]$,
$[\alpha_{(1,0)}^+ \alpha_{(1,1)}^-] = [\xi'] \oplus [\varphi]$,
$[\alpha_{(1,1)}^+ \alpha_{(0,1)}^-] = [\alpha_{(1,0)}^+\alpha_{(0,1)}^-] \oplus [\gamma_1] \oplus [\gamma_2]$,
$[\alpha_{(0,1)}^+ \alpha_{(1,1)}^-] = [\alpha_{(0,1)}^+\alpha_{(1,0)}^-] \oplus [\gamma_1'] \oplus [\gamma_2']$,
$[\alpha_{(1,1)}^+ \alpha_{(2,0)}^-] = [\alpha_{(1,0)}^+\alpha_{(2,0)}^-] \oplus [\delta_1] \oplus [\delta_2]$,
$[\alpha_{(3,0)}^+ \alpha_{(1,0)}^-] = [\xi] \oplus [\lambda]$,
$[\alpha_{(1,0)}^+ \alpha_{(3,0)}^-] = [\xi'] \oplus [\lambda']$,
$[\alpha_{(3,0)}^+ \alpha_{(0,1)}^-] = [\gamma_1] \oplus [\gamma_2] \oplus [\omega_1] \oplus [\omega_2]$ and
$[\alpha_{(2,1)}^+ \alpha_{(1,0)}^-] = [\alpha_{(0,1)}^+\alpha_{(1,0)}^-] \oplus [\alpha_{(2,0)}^+\alpha_{(1,0)}^-] \oplus [\phi_1] \oplus [\phi_2]$.

\begin{figure}[tb]
\begin{center}
  \includegraphics[width=130mm]{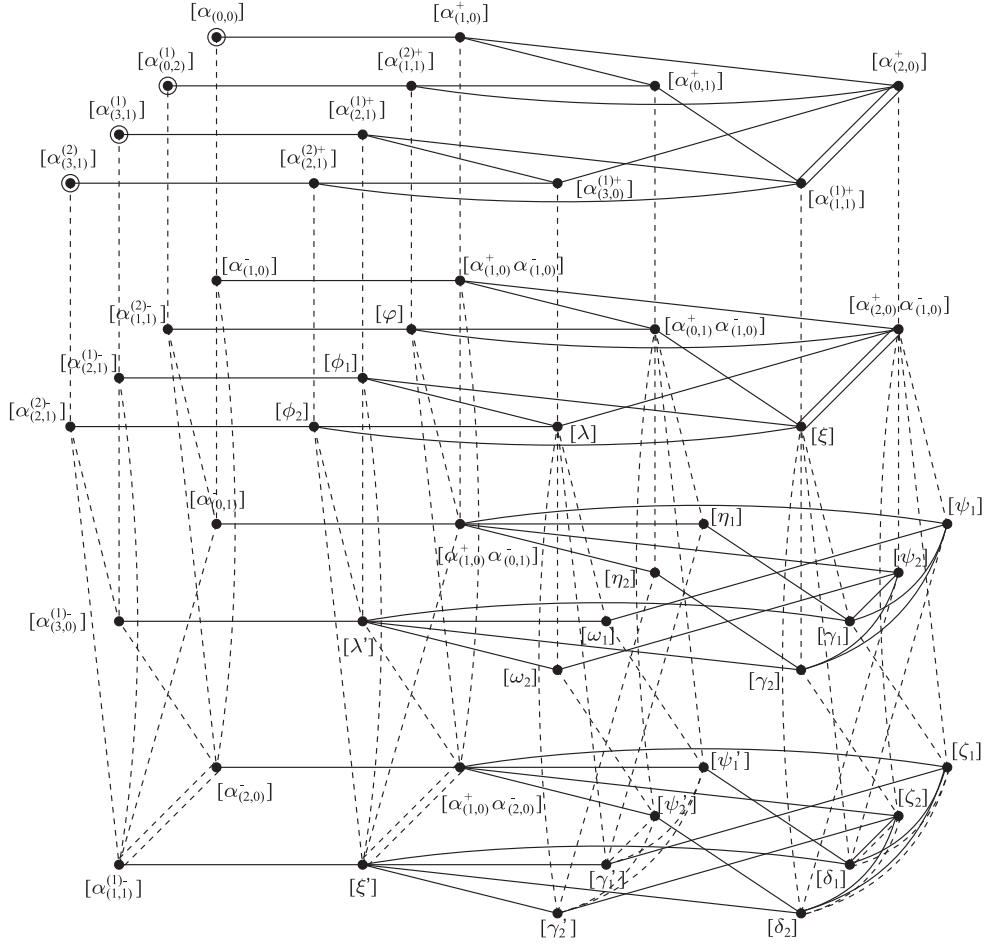} \\
 \caption{$\mathcal{E}_7(Sp(2))$: Multiplication by $[\alpha_{(1,0)}^+]$ (solid lines) and $[\alpha_{(1,0)}^-]$ (dashed lines)} \label{Fig-full_E7_C2}
\end{center}
\end{figure}

The fusion graphs of $[\alpha_{(1,0)}^+]$ (solid lines) and $[\alpha_{(1,0)}^-]$ (dashed lines) are given in Figure \ref{Fig-full_E7_C2}, where we have circled the marked vertices. Here multiplication by $[\alpha_{(1,0)}^+]$ (or $[\alpha_{(1,0)}^-]$) gives two copies each of $\mathcal{E}_7(Sp(2))$ and $\mathcal{E}_7^M(Sp(2))$. The ambichiral part ${}_M \mathcal{X}_M^0$ obeys $\mathbb{Z}_2 \times \mathbb{Z}_2$ fusion rules, corresponding to $SO(14)$ at level 1.

\begin{figure}[tb]
\begin{center}
  \includegraphics[width=160mm]{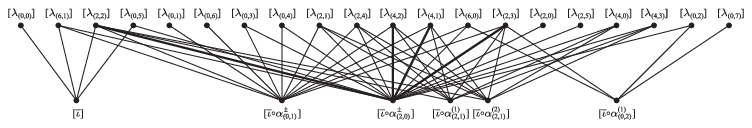} \\
 \caption{$\mathcal{E}_7(Sp(2))$: Principal graph of $(Sp(2))_7 \rightarrow (SO(14))_1$} \label{Fig-GHJ_Graph_E7_C2}
\end{center}
\end{figure}

We find
\begin{equation}
[\gamma] = [\alpha_{(0,0)}] \oplus [\alpha_{(1,0)}^+\alpha_{(1,0)}^-] \oplus [\eta_1] \oplus [\zeta_1],
\end{equation}
and the principal graph of the inclusion $(Sp(2))_7 \rightarrow (SO(14))_1$ of index $5(3 + \sqrt{5}) + \sqrt{250 + 110 \sqrt{5}} \approx 48.45$ is illustrated in Figure \ref{Fig-GHJ_Graph_E7_C2}, where the thick lines denote double edges.

Again, we can construct a subfactor $\alpha_{(1,0)}^{\pm}(M) \subset M$ of index $4 + \sqrt{5} + 2\sqrt{5 + 2\sqrt{5}} \approx 12.39$, where $M$ is a type III factor.
Its principal graph is the nimrep graph $\mathcal{E}_7^{\rho_1}(Sp(2))$ illustrated in Figure \ref{fig-Graph_E7_B2-1}. The dual principal graph is again isomorphic to the principal graph as abstract graphs.

We now determine the joint spectral measure of $\mathcal{E}_7^{\rho_1}(Sp(2))$, $\mathcal{E}_7^{\rho_2}(Sp(2))$.
With $\theta_1, \theta_2$ as in (\ref{eqn:theta-lambda_C2}) for $\lambda = (\lambda_1,\lambda_2) \in \mathrm{Exp}(\mathcal{E}_7(Sp(2)))$, we have the following values:
\begin{center}
\begin{tabular}{|c|c|c|c|} \hline
$\lambda \in \mathrm{Exp}$ & $(\theta_1,\theta_2) \in [0,1]^2$ & $|\psi^{\lambda}_{\ast}|^2$ & $\frac{1}{8\pi^2}|J(\theta_1,\theta_2)|$ \\
\hline $(0,0)$, $(0,7)$ & $\left(\frac{1}{20},\frac{2}{20}\right)$, $\left(\frac{8}{20},\frac{9}{20}\right)$ & $\frac{5-\sqrt{5}-\sqrt{10-2\sqrt{5}}}{80}$ & $\frac{5-\sqrt{5}-\sqrt{10-2\sqrt{5}}}{4}$ \\
\hline $(6,1)$, $(6,0)$ & $\left(\frac{2}{20},\frac{9}{20}\right)$, $\left(\frac{1}{20},\frac{8}{20}\right)$ & $\frac{5-\sqrt{5}+\sqrt{10-2\sqrt{5}}}{80}$ & $\frac{5-\sqrt{5}+\sqrt{10-2\sqrt{5}}}{4}$ \\
\hline $(2,2)$, $(2,3)$ & $\left(\frac{3}{20},\frac{6}{20}\right)$, $\left(\frac{4}{20},\frac{7}{20}\right)$ & $\frac{5+\sqrt{5}+\sqrt{10+2\sqrt{5}}}{80}$ & $\frac{5+\sqrt{5}+\sqrt{10+2\sqrt{5}}}{4}$ \\
\hline $(0,5)$, $(0,2)$ & $\left(\frac{6}{20},\frac{7}{20}\right)$, $\left(\frac{3}{20},\frac{4}{20}\right)$ & $\frac{5+\sqrt{5}-\sqrt{10+2\sqrt{5}}}{80}$ & $\frac{5+\sqrt{5}-\sqrt{10+2\sqrt{5}}}{4}$ \\
\hline $(3,1)$, $(3,3)$ & $\left(\frac{2}{20},\frac{6}{20}\right)$, $\left(\frac{4}{20},\frac{8}{20}\right)$ & $\frac{1}{4}$ & $\frac{5}{2}$ \\
\hline
\end{tabular}
\end{center}
where the eigenvectors $\psi^{\lambda}$ have been normalized so that $||\psi^{\lambda}|| = 1$.
Note that
\begin{equation} \label{eqn:psi=zetaJ-E7}
|\psi^{\lambda}_{\ast}|^2 = \zeta_{\lambda} \frac{1}{20} \frac{1}{8\pi^2} |J|,
\end{equation}
where $\zeta_{\lambda} = 1$ for $\lambda \in \{ (0,0), (6,1), (2,2), (0,5), (6,0), (0,2), (2,3), (0,7) \}$ and $\zeta_{(3,1)} = \zeta_{(3,3)} = 2$.

\begin{figure}[tb]
\begin{center}
  \includegraphics[width=45mm]{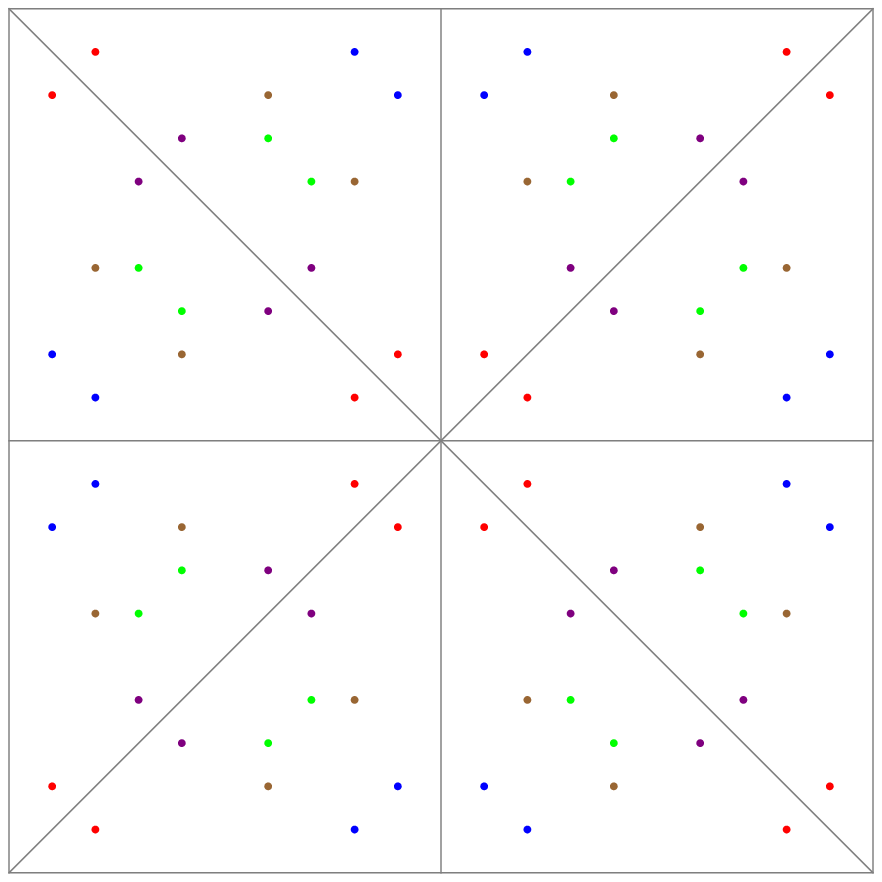}\\
 \caption{Orbit of $(\theta_1,\theta_2)$ for $\lambda \in \mathrm{Exp}(\mathcal{E}_7(Sp(2)))$.} \label{Fig-E7pointsB2}
\end{center}
\end{figure}

The orbits under $D_8$ of the points $(\theta_1,\theta_2) \in \left\{ \left(\frac{1}{20},\frac{2}{20}\right), \left(\frac{8}{20},\frac{9}{20}\right)\textcolor{red}{\cdot}, \left(\frac{2}{20},\frac{9}{20}\right), \left(\frac{1}{20},\frac{8}{20}\right)\textcolor{blue}{\cdot}, \right.$ $\left. \left(\frac{3}{20},\frac{6}{20}\right), \left(\frac{4}{20},\frac{7}{20}\right)\textcolor{green}{\cdot}, \left(\frac{6}{20},\frac{7}{20}\right), \left(\frac{3}{20},\frac{4}{20}\right)\textcolor{purple}{\cdot}, \left(\frac{2}{20},\frac{6}{20}\right), \left(\frac{4}{20},\frac{8}{20}\right)\textcolor{brown}{\cdot}, \right\}$ are illustrated in Figure \ref{Fig-E7pointsB2}.
The orbits of each successive pair of points support the measures $\mathrm{d}^{(1/20,2/20)}$, $\mathrm{d}^{(1/20,8/20)}$, $\mathrm{d}^{(3/20,6/20)}$, $\mathrm{d}^{(3/20,4/20)}$ and $\mathrm{d}^{(2/20,6/20)}$ respectively.
Thus, using (\ref{eqn:psi=zetaJ-E7}), we see that the joint spectral measure for $\mathcal{E}_7(Sp(2))$ is
$$\mathrm{d}\varepsilon = 16 \, \frac{1}{8} \, \frac{1}{20} \, \frac{1}{8\pi^2} |J| \, \left( \mathrm{d}^{(1/20,2/20)} + \mathrm{d}^{(1/20,2/5)} + \mathrm{d}^{(3/20,3/10)} + \mathrm{d}^{(3/20,1/5)} + 2\mathrm{d}^{(1/10,3/10)} \right).$$
Then we have obtained the following result:

\begin{Thm}
The joint spectral measure of $\mathcal{E}_7^{\rho_1}(Sp(2))$, $\mathcal{E}_7^{\rho_2}(Sp(2))$ (over $\mathbb{T}^2$) is
\begin{equation}
\mathrm{d}\varepsilon = \frac{1}{80\pi^2} |J| \, \left( \mathrm{d}^{(1/20,2/20)} + \mathrm{d}^{(1/20,2/5)} + \mathrm{d}^{(3/20,3/10)} + \mathrm{d}^{(3/20,1/5)} + 2\mathrm{d}^{(1/10,3/10)} \right),
\end{equation}
where $\mathrm{d}^{(\theta_1,\theta_2)}$ is as in Definition \ref{def:B2measure}.
\end{Thm}

\subsection{Exceptional Graph $\mathcal{E}_7^M(Sp(2))$: $(Sp(2))_7 \rightarrow (SO(14))_1 \rtimes \mathbb{Z}_2$}

\begin{figure}[tb]
\begin{minipage}[t]{5cm}
\begin{center}
  \includegraphics[width=40mm]{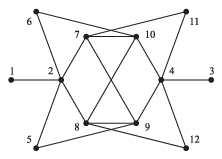} \\
 \caption{Graph $\mathcal{E}_7^{M,\rho_1}(Sp(2))$} \label{fig-Graph_E7M_B2-1}
\end{center}
\end{minipage}
\hfill
\begin{minipage}[t]{10cm}
\begin{center}
  \includegraphics[width=100mm]{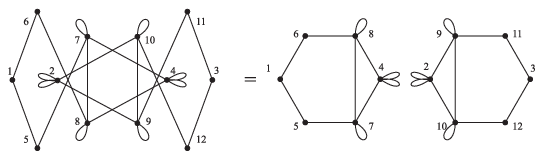} \\
 \caption{Graph $\mathcal{E}_7^{M,\rho_2}(Sp(2))$} \label{fig-Graph_E7M_B2-2}
\end{center}
\end{minipage}
\end{figure}

The graphs $\mathcal{E}_7^{M,\rho_j}(Sp(2))$, illustrated in Figures \ref{fig-Graph_E7M_B2-1} and \ref{fig-Graph_E7M_B2-2} are the nimrep graphs for the type II inclusion $(Sp(2))_7 \rightarrow (SO(14))_1 \rtimes_{\tau} \mathbb{Z}_2$ with index $10(3 + \sqrt{5}) + 2\sqrt{250 + 110 \sqrt{5}} \approx 96.90$, where $\tau = \alpha_{(0,2)}^{(1)}$ is a non-trivial simple current of order 2 in the ambichiral system ${}_M \mathcal{X}_M^0$, see Figure \ref{Fig-full_E7_C2}. From Section \ref{sect:measures_E7C2}, $[\tau]$ is a subsector of $[\alpha_{(0,2)}^{\pm}]$. Now $\omega_{(0,2)} = -1$ \cite{coquereaux:webpage}, which satisfies $\omega_{(0,2)}^2 = 1$, thus the orbifold inclusion exists (c.f. Section \ref{sect:measures_DkB2}). Note that $[\tau']=[\alpha_{(3,1)}^{(j)}] \in {}_M \mathcal{X}_M^0$, $j=1,2$, is a subsector of $[\alpha_{(3,1)}^{\pm}]$, for which $\omega_{(3,1)} = e^{7\pi i/4}$ \cite{coquereaux:webpage}. Then $\omega_{(3,1)}^2 = e^{7\pi i/2} \neq 1$, and hence the orbifold inclusion $(Sp(2))_7 \rightarrow (SO(14))_1 \rtimes_{\tau'} \mathbb{Z}_2$ does not exist.

\begin{figure}[tb]
\begin{center}
  \includegraphics[width=90mm]{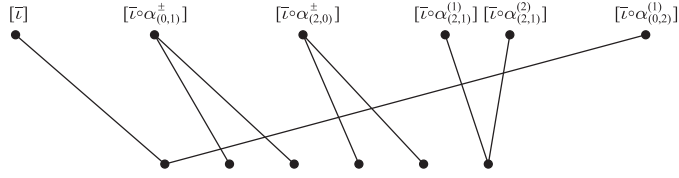} \\
 \caption{$\mathcal{E}_7^M(Sp(2))$: $\mathbb{Z}_2$-action on $(SO(14))_1$} \label{Fig-GHJ_Graph_factor_E7M_C2}
\end{center}
\end{figure}

The principal graph for this inclusion would be the graph obtained by composing the principal graph for $(Sp(2))_7 \rightarrow (SO(14))_1$, illustrated in Figure \ref{Fig-GHJ_Graph_E7_C2} with the graph for the $\mathbb{Z}_2$-action, illustrated in Figure \ref{Fig-GHJ_Graph_factor_E7M_C2}. This will be discussed in a future publication using a generalised Goodman-de la Harpe-Jones construction (c.f. the comments in Section \ref{sect:measures_E3MC2}). Again, it is not clear what the dual principal graph is in this case.

The associated modular invariant is again $Z_{\mathcal{E}_7}$ and the graphs are isospectral to $\mathcal{E}_7(Sp(2))$.
The eigenvectors $\psi^{\lambda}$ are not identical to those for $\mathcal{E}_3(Sp(2))$. However, as seen in the following table, the values of $|\psi^{\lambda}_{\ast}|^2$ are equal (up to a factor 2) to those for $\mathcal{E}_3(Sp(2))$, for $\lambda \neq (3,1), (3,3)$.
With $\theta_1, \theta_2$ as in (\ref{eqn:theta-lambda_C2}) for $\lambda = (\lambda_1,\lambda_2) \in \mathrm{Exp}$, we have:
\begin{center}
\begin{tabular}{|c|c|c|c|} \hline
$\lambda \in \mathrm{Exp}$ & $(\theta_1,\theta_2) \in [0,1]^2$ & $|\psi^{\lambda}_{\ast}|^2$ & $\frac{1}{8\pi^2}|J(\theta_1,\theta_2)|$ \\
\hline $(0,0)$, $(0,7)$ & $\left(\frac{1}{20},\frac{2}{20}\right)$, $\left(\frac{8}{20},\frac{9}{20}\right)$ & $\frac{5-\sqrt{5}-\sqrt{10-2\sqrt{5}}}{40}$ & $\frac{5-\sqrt{5}-\sqrt{10-2\sqrt{5}}}{4}$ \\
\hline $(6,1)$, $(6,0)$ & $\left(\frac{2}{20},\frac{9}{20}\right)$, $\left(\frac{1}{20},\frac{8}{20}\right)$ & $\frac{5-\sqrt{5}+\sqrt{10-2\sqrt{5}}}{40}$ & $\frac{5-\sqrt{5}+\sqrt{10-2\sqrt{5}}}{4}$ \\
\hline $(2,2)$, $(2,3)$ & $\left(\frac{3}{20},\frac{6}{20}\right)$, $\left(\frac{4}{20},\frac{7}{20}\right)$ & $\frac{5+\sqrt{5}+\sqrt{10+2\sqrt{5}}}{40}$ & $\frac{5+\sqrt{5}+\sqrt{10+2\sqrt{5}}}{4}$ \\
\hline $(0,5)$, $(0,2)$ & $\left(\frac{6}{20},\frac{7}{20}\right)$, $\left(\frac{3}{20},\frac{4}{20}\right)$ & $\frac{5+\sqrt{5}-\sqrt{10+2\sqrt{5}}}{40}$ & $\frac{5+\sqrt{5}-\sqrt{10+2\sqrt{5}}}{4}$ \\
\hline $(3,1)$, $(3,3)$ & $\left(\frac{2}{20},\frac{6}{20}\right)$, $\left(\frac{4}{20},\frac{8}{20}\right)$ & $0$ & $\frac{5}{2}$ \\
\hline
\end{tabular}
\end{center}
where the eigenvectors $\psi^{\lambda}$ have been normalized so that $||\psi^{\lambda}|| = 1$.
Then (\ref{eqn:psi=zetaJ-E7}) becomes $|\psi^{\lambda}_{\ast}|^2 = \zeta_{\lambda} \frac{1}{10} \frac{1}{8\pi^2} |J|$, where $\zeta_{\lambda} = 1$ for $\lambda \in \{ (0,0), (6,1), (2,2), (0,5), (6,0), (0,2), (2,3), (0,7) \}$ as for $\mathcal{E}_7(Sp(2))$, and $\zeta_{(3,1)} = \zeta_{(3,3)} = 0$.
Thus we have the following result:

\begin{Thm}
The joint spectral measure of $\mathcal{E}_7^{M,\rho_1}(Sp(2))$, $\mathcal{E}_7^{M,\rho_2}(Sp(2))$ (over $\mathbb{T}^2$) is
\begin{equation}
\mathrm{d}\varepsilon = \frac{1}{80\pi^2} |J| \, \left( \mathrm{d}^{(1/20,2/20)} + \mathrm{d}^{(1/20,2/5)} + \mathrm{d}^{(3/20,3/10)} + \mathrm{d}^{(3/20,1/5)} \right),
\end{equation}
where $\mathrm{d}^{(\theta_1,\theta_2)}$ is as in Definition \ref{def:B2measure}.
\end{Thm}

\subsection{Exceptional Graph $\mathcal{E}_8(Sp(2))$}

\begin{figure}[tb]
\begin{minipage}[t]{5.2cm}
\begin{center}
  \includegraphics[width=50mm]{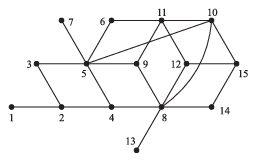} \\
 \caption{Graph $\mathcal{E}_8^{\rho_1}(Sp(2))$} \label{fig-Graph_E8_B2-1}
\end{center}
\end{minipage}
\hfill
\begin{minipage}[t]{10.2cm}
\begin{center}
  \includegraphics[width=100mm]{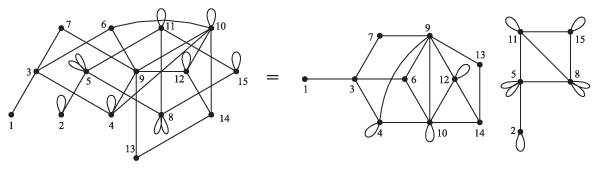} \\
 \caption{Graph $\mathcal{E}_8^{\rho_2}(Sp(2))$} \label{fig-Graph_E8_B2-2}
\end{center}
\end{minipage}
\end{figure}

The graphs $\mathcal{E}_8^{\rho_j}(G_2)$ are illustrated in Figures \ref{fig-Graph_E8_B2-1} and \ref{fig-Graph_E8_B2-2}.
To our knowledge neither graph has appeared in the literature before in the context of nimrep graphs or subfactors.
The associated modular invariant is \cite[(3.1)]{verstegen:1990}
\begin{align*}
Z_{\mathcal{E}_8} &= |\chi_{(0,0)}+\chi_{(0,8)}|^2 + |\chi_{(0,2)}+\chi_{(0,6)}|^2 + |\chi_{(4,0)}+\chi_{(4,4)}|^2 + |\chi_{(4,1)}+\chi_{(4,3)}|^2 \\
& \qquad + |\chi_{(2,2)}+\chi_{(2,4)}|^2 + (\chi_{(8,0)},\chi_{(0,1)},\chi_{(0,7)}) + (\chi_{(6,1)},\chi_{(0,3)},\chi_{(0,5)}) + (\chi_{(4,2)},\chi_{(2,1)},\chi_{(2,5)}) \\
& \qquad + (\chi_{(2,3)},\chi_{(6,0)},\chi_{(6,2)}) + (\chi_{(0,4)},\chi_{(2,0)},\chi_{(2,6)}),
\end{align*}
where $(\chi_{\lambda},\chi_{\mu},\chi_{\nu}) := |\chi_{\lambda}|^2 + \chi_{\lambda}(\chi_{\mu}+\chi_{\nu})^{\ast} + (\chi_{\mu}+\chi_{\nu})\chi_{\lambda}^{\ast}$.
This modular invariant is at level 8 and has exponents
\begin{align*}
\mathrm{Exp}(\mathcal{E}_8(Sp(2))) &= \{ (0,0), (0,8), (0,2), (0,6), (4,0), (4,4), (4,1), (4,3), (2,2), (2,4), \\
& \qquad (8,0), (6,1), (4,2), (2,3), (0,4) \}.
\end{align*}
Since the modular invariant associated with this family of graphs does not come from a conformal embedding, it has not yet been shown that the graphs $\mathcal{E}_8{\rho_j}(Sp(2))$ arises from a braided subfactor.
This modular invariant is a twist of the $\mathcal{D}_8(Sp(2)) = \mathcal{A}_8(Sp(2))/\mathbb{Z}_2$ orbifold invariant discussed in Section \ref{sect:measures_DkB2}, and is analogous to the $E_7$ modular invariant for $SU(2)$ \cite[$\S5.3$]{bockenhauer/evans/kawahigashi:2000} and the Moore-Seiberg $\mathcal{E}_{MS}^{(12)}$ invariant for $SU(3)$ \cite[$\S5.4$]{evans/pugh:2009ii}. The realisation of this nimrep by a braided subfactor will be discussed in a future publication, using a generalised Goodman-de la Harpe-Jones construction analogous to that for $E_7$, $\mathcal{E}_{MS}^{(12)}$ in \cite{bockenhauer/evans/kawahigashi:2000, evans/pugh:2009ii}.  This construction produces $\mathcal{E}_8^{\rho_j}(G_2)$ as nimrep graphs.
It is expected that $\mathcal{E}_8{\rho_j}(Sp(2))$ does indeed arise as the nimrep for a type II inclusion with index $4 (\cos(\pi/11) + \cos(2\pi/11))^2 \approx 12.97$.

However, for our purposes it is sufficient to know the eigenvalues and corresponding eigenvectors for these graphs, and it is not necessary for the graph to be a nimrep graph.
With $\theta_1, \theta_2$ as in (\ref{eqn:theta-lambda_C2}) for $\lambda = (\lambda_1,\lambda_2) \in \mathrm{Exp}(\mathcal{E}_8(Sp(2)))$, we have the following values:
\begin{center}
\begin{tabular}{|c|c|c|c|} \hline
$\lambda \in \mathrm{Exp}$ & $(\theta_1,\theta_2) \in [0,1]^2$ & $|\psi^{\lambda}_{\ast}|^2$ & $\frac{1}{64\pi^2}J_{y,z}(\theta_1,\theta_2)^2$ \\
\hline $(0,0)$, $(0,8)$ & $\left(\frac{1}{22},\frac{2}{22}\right)$, $\left(\frac{9}{22},\frac{10}{22}\right)$ & $a_5$ & $\frac{121}{2}a_5$ \\
\hline $(0,2)$, $(0,6)$ & $\left(\frac{3}{22},\frac{4}{22}\right)$, $\left(\frac{7}{22},\frac{8}{22}\right)$ & $a_4$ & $\frac{121}{2}a_4$ \\
\hline $(4,0)$, $(4,4)$ & $\left(\frac{1}{22},\frac{6}{22}\right)$, $\left(\frac{5}{22},\frac{10}{22}\right)$ & $a_3$ & $\frac{121}{2}a_3$ \\
\hline $(4,1)$, $(4,3)$ & $\left(\frac{2}{22},\frac{7}{22}\right)$, $\left(\frac{4}{22},\frac{9}{22}\right)$ & $a_2$ & $\frac{121}{2}a_2$ \\
\hline $(2,2)$, $(2,4)$ & $\left(\frac{3}{22},\frac{6}{22}\right)$, $\left(\frac{5}{22},\frac{8}{22}\right)$ & $a_1$ & $\frac{121}{2}a_1$ \\
\hline $(8,0)$ & $\left(\frac{1}{22},\frac{10}{22}\right)$ & $11b_1$ & $0$ \\
\hline $(6,1)$ & $\left(\frac{2}{22},\frac{9}{22}\right)$ & $11b_5$ & $0$ \\
\hline $(4,2)$ & $\left(\frac{3}{22},\frac{8}{22}\right)$ & $11b_2$ & $0$ \\
\hline $(2,3)$ & $\left(\frac{4}{22},\frac{7}{22}\right)$ & $11b_4$ & $0$ \\
\hline $(0,4)$ & $\left(\frac{5}{22},\frac{6}{22}\right)$ & $11b_3$ & $0$ \\
\hline
\end{tabular}
\end{center}
where $a_i$ is the $i^{\mathrm{th}}$ largest root of $56689952x^5 -15460896x^4 +1522664x^3 -63888x^2 +968x - 1$, $b_i$ is the $i^{\mathrm{th}}$ largest root of $x^5-11x^4+44x^3-77x^2+55x-11$, and the eigenvectors $\psi^{\lambda}$ have been normalized so that $||\psi^{\lambda}|| = 1$.

\begin{figure}[tb]
\begin{center}
  \includegraphics[width=45mm]{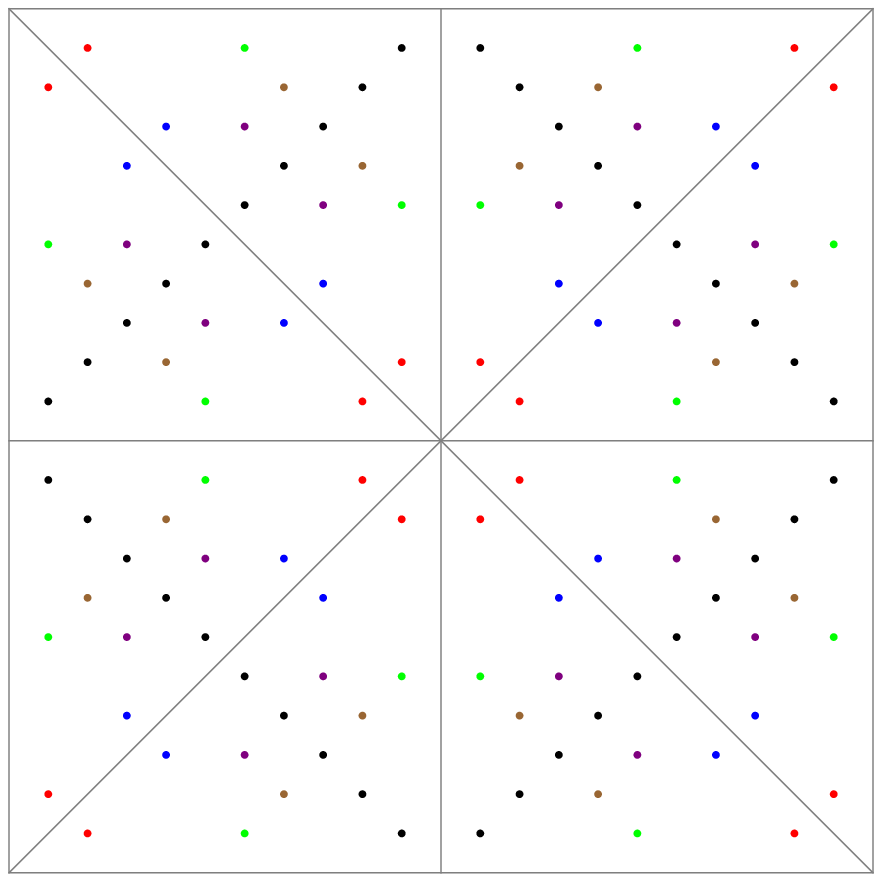}\\
 \caption{Orbit of $(\theta_1,\theta_2)$ for $\lambda \in \mathrm{Exp}(\mathcal{E}_8(Sp(2)))$.} \label{Fig-E8pointsB2}
\end{center}
\end{figure}

The orbits under $D_8$ of the points $(\theta_1,\theta_2) \in \left\{ \left(\frac{1}{22},\frac{2}{22}\right), \left(\frac{9}{22},\frac{10}{22}\right)\textcolor{red}{\cdot}, \left(\frac{3}{22},\frac{4}{22}\right), \left(\frac{7}{22},\frac{8}{22}\right)\textcolor{blue}{\cdot}, \right.$ $\left. \left(\frac{1}{22},\frac{6}{22}\right), \left(\frac{5}{22},\frac{10}{22}\right)\textcolor{green}{\cdot}, \left(\frac{2}{22},\frac{7}{22}\right), \left(\frac{4}{22},\frac{9}{22}\right)\textcolor{brown}{\cdot}, \left(\frac{3}{22},\frac{6}{22}\right), \left(\frac{5}{22},\frac{8}{22}\right)\textcolor{purple}{\cdot}, \right\}$ are illustrated in Figure \ref{Fig-E8pointsB2}.
The orbits of each successive pair of points support the measures $\mathrm{d}^{(1/20,2/22)}$, $\mathrm{d}^{(3/22,4/22)}$, $\mathrm{d}^{(1/22,6/22)}$, $\mathrm{d}^{(2/22,7/22)}$ and $\mathrm{d}^{(3/22,6/22)}$ respectively. The orbits under $D_8$ of the points $(\theta_1,\theta_2) \in \left\{ \left(\frac{1}{22},\frac{10}{22}\right), \left(\frac{2}{22},\frac{9}{22}\right), \left(\frac{3}{22},\frac{8}{22}\right), \left(\frac{4}{22},\frac{7}{22}\right), \left(\frac{5}{22},\frac{6}{22}\right) \right\}$ are the black points illustrated in Figure \ref{Fig-E8pointsB2}.
Thus we see that the joint spectral measure for $\mathcal{E}_8(Sp(2))$ is
\begin{align*}
&16 \, \frac{1}{8} \, \frac{2}{121} \, \frac{1}{64\pi^4} J_{y,z}^2 \, \left( \mathrm{d}^{(1/22,2/22)} + \mathrm{d}^{(3/22,4/22)} + \mathrm{d}^{(1/22,6/22)} + \mathrm{d}^{(3/22,6/22)} + 2\mathrm{d}^{(2/22,7/22)} \right) \\
& \quad + 8 \, \frac{1}{8} \, 11 \left( b_1 \, \mathrm{d}^{(1/22,10/22)} + b_2 \, \mathrm{d}^{(3/22,8/22)} + b_3 \, \mathrm{d}^{(5/22,6/22)} + b_4 \, \mathrm{d}^{(4/22,7/22)} + b_5 \, \mathrm{d}^{(2/22,9/22)} \right).
\end{align*}
Then we have the following result:

\begin{Thm}
The joint spectral measure of $\mathcal{E}_8^{\rho_1}(Sp(2))$, $\mathcal{E}_8^{\rho_2}(Sp(2))$ (over $\mathbb{T}^2$) is
\begin{align}
\mathrm{d}\varepsilon &= \frac{1}{1936\pi^4} J_{y,z}^2 \, \left( \mathrm{d}^{(1/22,2/22)} + \mathrm{d}^{(3/22,4/22)} + \mathrm{d}^{(1/22,6/22)} + \mathrm{d}^{(3/22,6/22)} + 2\mathrm{d}^{(2/22,7/22)} \right) \\
& \qquad + 11 \left( b_1 \, \mathrm{d}^{(1/22,10/22)} + b_2 \, \mathrm{d}^{(3/22,8/22)} + b_3 \, \mathrm{d}^{(5/22,6/22)} + b_4 \, \mathrm{d}^{(4/22,7/22)} + b_5 \, \mathrm{d}^{(2/22,9/22)} \right),
\end{align}
where $b_i$ is the $i^{\mathrm{th}}$ largest root of $x^5-11x^4+44x^3-77x^2+55x-11$, and $\mathrm{d}^{(\theta_1,\theta_2)}$ is as in Definition \ref{def:B2measure}.
\end{Thm}

\subsection{Exceptional Graph $\mathcal{E}_{12}(Sp(2))$: $(Sp(2))_{12} \rightarrow (E_8)_1$}

\begin{figure}[tb]
\begin{minipage}[t]{5cm}
\begin{center}
  \includegraphics[width=35mm]{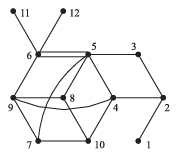} \\
 \caption{Graph $\mathcal{E}_{12}^{\rho_1}(Sp(2))$} \label{fig-Graph_E12_B2-1}
\end{center}
\end{minipage}
\hfill
\begin{minipage}[t]{10cm}
\begin{center}
  \includegraphics[width=100mm]{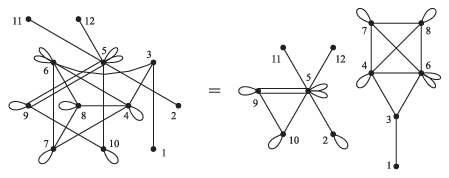} \\
 \caption{Graph $\mathcal{E}_{12}^{\rho_2}(Sp(2))$} \label{fig-Graph_E12_B2-2}
\end{center}
\end{minipage}
\end{figure}

The graphs $\mathcal{E}_{12}^{\rho_j}(Sp(2))$, illustrated in Figures \ref{fig-Graph_E12_B2-1}, \ref{fig-Graph_E12_B2-2}, are the nimrep graphs associated to the conformal embedding $(Sp(2))_{12} \rightarrow (E_8)_1$ and are the graphs associated to the modular invariant
$$Z_{\mathcal{E}_{12}} = |\chi_{(0,0)}+\chi_{(6,0)}+\chi_{(8,1)}+\chi_{(2,3)}+\chi_{(8,3)}+\chi_{(6,6)}+\chi_{(2,7)}+\chi_{(0,12)}+2\chi_{(4,4)}|^2$$
which is at level 12 and has exponents
$$\mathrm{Exp}(\mathcal{E}_{12}(Sp(2))) = \{ (0,0), (6,0), (8,1), (2,3), (8,3), (6,6), (2,7), (0,12), \textrm{ and four times } (4,4). \}$$
The graphs $\mathcal{E}_{12}^{\rho_j}(Sp(2))$ are illustrated in Figures \ref{fig-Graph_E12_B2-1}, \ref{fig-Graph_E12_B2-2}. Note again that $\mathcal{E}_{12}^{\rho_2}(Sp(2))$ has two connected components.

The chiral induced sector bases ${}_M \mathcal{X}_M^{\pm}$ are given by
\begin{align*}
{}_M \mathcal{X}_M^{\pm} &= \{ [\alpha_{(0,0)}], [\alpha_{(1,0)}^{\pm}], [\alpha_{(0,1)}^{\pm}], [\alpha_{(2,0)}^{\pm}], [\alpha_{(1,1)}^{\pm}], [\alpha_{(0,2)}^{\pm}], [\alpha_{(3,0)}^{(j)\pm}], [\alpha_{(2,1)}^{(j)\pm}], [\alpha_{(1,2)}^{(j)\pm}], \mbox{ for } j=1,2 \},\\
\end{align*}
where $[\alpha_{(3,0)}^{\pm}] = [\alpha_{(3,0)}^{(1)\pm}] \oplus [\alpha_{(3,0)}^{(2)\pm}]$,
$[\alpha_{(2,1)}^{\pm}] = [\alpha_{(0,2)}^{\pm}] \oplus [\alpha_{(2,1)}^{(1)\pm}] \oplus [\alpha_{(2,1)}^{(2)\pm}]$, and
$[\alpha_{(1,2)}^{\pm}] = [\alpha_{(1,1)}^{\pm}] \oplus [\alpha_{(3,0)}^{(1)\pm}] \oplus [\alpha_{(1,2)}^{(1)\pm}] \oplus [\alpha_{(1,2)}^{(2)\pm}]$.

One can in principle compute the principal graph and dual principal graph of the inclusion $(Sp(2))_{12} \rightarrow (E_8)_1$, as in Section \ref{sect:measures_E3C2}, but we do not do that here due to their size (the principal graph for instance has 55 vertices).
It is only possible to determine the edge set of the pair of vertices $[\overline{\iota} \circ \alpha_{(2,1)}^{(1)\pm}]$ and $[\overline{\iota} \circ \alpha_{(2,1)}^{(1)\pm}]$ together, but not which edges are attached to either vertex individually.
However, the correct choice could be verified by the generalised Goodman-de la Harpe-Jones method referred to in Section \ref{sect:measures_E3C2}, where the principal graph appears as the intertwining graph.

The subfactor $\alpha_{(1,0)}^{\pm}(M) \subset M$ of index $\frac{1}{2}(10 + 3\sqrt{5} + \sqrt{75 + 30\sqrt{5}}) \approx 14.31$, where $M$ is a type III factor, has principal graph the nimrep graph $\mathcal{E}_{12}^{\rho_1}(Sp(2))$ illustrated in Figure \ref{fig-Graph_E12_B2-1}, and the dual principal graph is again isomorphic to the principal graph as abstract graphs.

We now determine the joint spectral measure of $\mathcal{E}_{12}^{\rho_1}(Sp(2))$, $\mathcal{E}_{12}^{\rho_2}(Sp(2))$.
With $\theta_1, \theta_2$ as in (\ref{eqn:theta-lambda_C2}) for $\lambda = (\lambda_1,\lambda_2) \in \mathrm{Exp}(\mathcal{E}_3(Sp(2)))$, we have the following values:
\begin{center}
\begin{tabular}{|c|c|c|c|} \hline
$\lambda \in \mathrm{Exp}$ & $(\theta_1,\theta_2) \in [0,1]^2$ & $|\psi^{\lambda}_{\ast}|^2$ & $\frac{1}{8\pi^2}|J(\theta_1,\theta_2)|$ \\
\hline $(0,0)$, $(0,12)$ & $\left(\frac{1}{30},\frac{2}{30}\right)$, $\left(\frac{13}{30},\frac{14}{30}\right)$ & $\frac{9-\sqrt{5}-\sqrt{30+6\sqrt{5}}}{120}$ & $\frac{9-\sqrt{5}-\sqrt{30+6\sqrt{5}}}{8}$ \\
\hline $(6,0)$, $(6,6)$ & $\left(\frac{1}{30},\frac{8}{30}\right)$, $\left(\frac{7}{30},\frac{14}{30}\right)$ & $\frac{9+\sqrt{5}-\sqrt{30-6\sqrt{5}}}{120}$ & $\frac{9+\sqrt{5}-\sqrt{30-6\sqrt{5}}}{8}$ \\
\hline $(8,1)$, $(8,3)$ & $\left(\frac{2}{30},\frac{11}{30}\right)$, $\left(\frac{4}{30},\frac{13}{30}\right)$ & $\frac{9+\sqrt{5}+\sqrt{30-6\sqrt{5}}}{120}$ & $\frac{9+\sqrt{5}+\sqrt{30-6\sqrt{5}}}{8}$ \\
\hline $(2,3)$, $(2,7)$ & $\left(\frac{4}{30},\frac{7}{30}\right)$, $\left(\frac{8}{30},\frac{11}{30}\right)$ & $\frac{9-\sqrt{5}+\sqrt{30+6\sqrt{5}}}{120}$ & $\frac{9-\sqrt{5}+\sqrt{30+6\sqrt{5}}}{8}$ \\
\hline $(4,4)$ & $\left(\frac{5}{30},\frac{10}{30}\right)$ & $\frac{2}{3}$ & $3$ \\
\hline
\end{tabular}
\end{center}
where the eigenvectors $\psi^{\lambda}$ have been normalized so that $||\psi^{\lambda}|| = 1$, and for the exponent $(4,4)$ which has multiplicity four, the value listed in the table for $|\psi^{(4,4)}_{\ast}|^2$ is $\sum_{j=1}^4 |\psi^{(1,1)_j}_{\ast}|^2$.
Note that
\begin{equation} \label{eqn:psi=zetaJ-E12}
|\psi^{\lambda}_{\ast}|^2 = \zeta_{\lambda} \frac{1}{15} \frac{1}{8\pi^2} |J|
\end{equation}
where $\zeta_{\lambda} = 1$ for $\lambda \in \mathrm{Exp}(\mathcal{E}_3(Sp(2))) \setminus \{ (4,4) \}$ and $\zeta_{(4,4)} = 10/3$.

\begin{figure}[tb]
\begin{center}
  \includegraphics[width=45mm]{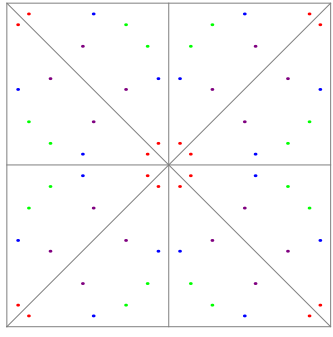}\\
 \caption{Orbit of $(\theta_1,\theta_2)$ for $\lambda \in \mathrm{Exp}(\mathcal{E}_{12}(Sp(2))) \setminus \{ (4,4) \}$.} \label{Fig-E12pointsB2}
\end{center}
\end{figure}

The orbits under $D_8$ of the points $(\theta_1,\theta_2) \in \left\{ \left(\frac{1}{30},\frac{2}{30}\right), \left(\frac{13}{30},\frac{14}{30}\right)\textcolor{red}{\cdot}, \left(\frac{1}{30},\frac{8}{30}\right), \left(\frac{7}{30},\frac{14}{30}\right)\textcolor{blue}{\cdot}, \right.$ $\left. \left(\frac{2}{30},\frac{11}{30}\right), \left(\frac{4}{30},\frac{13}{30}\right)\textcolor{green}{\cdot}, \left(\frac{4}{30},\frac{7}{30}\right), \left(\frac{8}{30},\frac{11}{30}\right)\textcolor{purple}{\cdot} \right\}$ are illustrated in Figure \ref{Fig-E12pointsB2}.
The orbits of the first four pairs of points support the measures $\mathrm{d}^{(1/30,2/30)}$, $\mathrm{d}^{(1/30,8/30)}$, $\mathrm{d}^{(2/30,11/30)}$ and $\mathrm{d}^{(4/30,7/30)}$ respectively. The orbit of the last pair has appeared before and supports the measure $|J| \, \mathrm{d}_6 \times \mathrm{d}_6$.
Thus, using (\ref{eqn:psi=zetaJ-E12}), we see that the joint spectral measure for $\mathcal{E}_{12}(Sp(2))$ is
\begin{align*}
\mathrm{d}\varepsilon &= 16 \, \frac{1}{8} \, \frac{1}{15} \, \frac{1}{8\pi^2} |J| \, \left( \mathrm{d}^{(1/30,1/15)} + \mathrm{d}^{(1/30,4/15)} + \mathrm{d}^{(1/15,11/30)} + \mathrm{d}^{(2/15,7/30)} \right) \\
& \qquad + 36 \, \frac{1}{8} \, \frac{10}{3} \, \frac{1}{15} \, \frac{1}{8\pi^2} |J| \, \mathrm{d}_6 \times \mathrm{d}_6.
\end{align*}
Then we have obtained the following result:

\begin{Thm}
The joint spectral measure of $\mathcal{E}_{12}^{\rho_1}(Sp(2))$, $\mathcal{E}_{12}^{\rho_2}(Sp(2))$ (over $\mathbb{T}^2$) is
\begin{equation}
\mathrm{d}\varepsilon = \frac{1}{60\pi^2} |J| \, \left( \mathrm{d}^{(1/30,1/15)} + \mathrm{d}^{(1/30,4/15)} + \mathrm{d}^{(1/15,11/30)} + \mathrm{d}^{(2/15,7/30)} \right) + \frac{1}{8\pi^2} |J| \, \mathrm{d}_6 \times \mathrm{d}_6,
\end{equation}
where $\mathrm{d}^{(\theta_1,\theta_2)}$ is as in Definition \ref{def:B2measure} and $\mathrm{d}_6$ is the uniform Dirac measure on the $6^{\mathrm{th}}$ roots of unity.
\end{Thm}

\bigskip \bigskip

\begin{footnotesize}
\noindent{\it Acknowledgement.}

The second author was supported by the Coleg Cymraeg Cenedlaethol.
\end{footnotesize}

\end{document}